\baselineskip=15pt plus 2pt
\magnification =1150
\def\sqr#1#2{{\vcenter{\vbox{\hrule height.#2pt\hbox{\vrule width.#2pt
height#1pt\kern#1pt \vrule width.#2pt}\hrule height.#2pt}}}}
\def\square{\mathchoice\sqr64\sqr64\sqr{2.1}3\sqr{1.5}3} 
\noindent {\bf Hill's spectral curves and the invariant measure of the periodic KdV
equation}\par
\vskip.05in
\centerline {Gordon Blower${}^a$, Caroline Brett${}^a$ and Ian
Doust${}^b$}\par
\vskip.05in
\noindent ${}^a$ Department of Mathematics and Statistics, Lancaster University, 
Lancaster, LA1 4YF  England, UK\par
\noindent ${}^a$ School of Mathematics, University of New South Wales, Sydney,
 Australia\par
\vskip.05in
\centerline {27th September 2014}\par
\vskip.05in
\hrule
\vskip.05in
\noindent {\bf Abstract} This paper analyses the periodic spectrum 
of Schr\"odinger's equation $-f''+qf=\lambda f$ when the 
potential is real, periodic, random and  subject to the 
invariant measure $\nu_N^\beta$ of the periodic KdV equation. This $\nu_N^\beta$ is
 the modified canonical ensemble, as given by Bourgain ({Comm. Math. Phys.} 
{166} (1994), 1--26), and $\nu_N^\beta$ satisfies a logarithmic Sobolev 
inequality. Associated concentration inequalities control the fluctuations of the 
periodic eigenvalues $(\lambda_n)$. For $\beta,
N>0$ small,  there exists a set of positive $\nu_N^\beta$ measure such that $(\pm \sqrt{2(\lambda_{2n}+\lambda_{2n-1})})_{n=0}^\infty$ gives a sampling 
sequence for Paley--Wiener space $PW(\pi )$ and the reproducing
kernels give a Riesz basis. Let $(\mu_j)_{j=1}^\infty$ be the tied spectrum; then 
$(2\sqrt{\mu_j}-j)$ belongs to a Hilbert cube in $\ell^2$ and is
distributed according to a measure that satisfies Gaussian concentration for Lipschitz
functions. The 
sampling sequence $(\sqrt{\mu_j})_{j=1}^\infty$ arises from a divisor on the spectral curve, which 
is hyperelliptic of infinite genus. The linear
 statistics $\sum_j g(\sqrt{\lambda_{2j}})$ with test function $g\in PW(\pi)$ satisfy Gaussian concentration 
inequalities.\par
\par
\noindent {\sl Classification:} 37L55\par
\vskip.05in
\hrule
\vskip.05in
\vfill
\eject
\noindent {\bf 1 Introduction}\par
\vskip.05in
\noindent A large class of systems can be modelled via differential equations of the form
$$   -f''(x) + q(x) f(x) = \lambda f(x) \qquad (x \in {\bf R})\eqno(1.1)$$
where $q$ is a periodic potential function. Here 
$q:{\bf R}\rightarrow {\bf R}$ is a $2\pi$-periodic and measurable function such that 
$N=\int_0^{2\pi} q(x)^2dx/(2\pi )$ is finite. Equation 
(1.1) is known variously as Hill's equation or the time-independent
 Schr{\"o}dinger equation for potential scattering. By classical 
results, [20], Equation (1.1) admits an infinite increasing sequence of real eigenvalues 
$$ \lambda_0 < \lambda_1 \leq \lambda_2 < \lambda_3 \leq \lambda_4 < \lambda_5 \leq \dots\eqno(1.2)$$
each corresponding to either a nontrivial $2\pi$-periodic solution giving the principal series
of eigenvalues, or else an antiperiodic
 solution satisfying $f(x+2\pi) = -f(x)$ giving the complementary series. 
The periodic spectrum $\{\lambda_j: j=0, 1,\dots \}$ partitions ${\bf R}$ 
into intervals of 
stability and instability for the nontrivial solutions for (1.1).\par
\indent Hill's equation is closely related to the periodic Korteweg-de~Vries (KdV) equation
$${{\partial u}\over{\partial t}} + {{\partial^3 u}\over{\partial x^3}} + \beta u 
{{\partial u}\over{\partial x}} = 0.\eqno(1.3)$$
\noindent Gardner, Greene, Kruskal and Miura [9] and Lax [17] noted that periodic spectrum of (1.1) is 
preserved if a time-dependent potential $q_t(x)$ evolves according to (1.3).\par 
\indent For suitable fixed $q$, much is known about the asymptotic behaviour 
of $\{\lambda_j\}$. Let
$$ \Omega_N = \Bigl\{\phi \in L^2({\bf T};{\bf R}) \,:\, \int_{\bf T} \phi(x)^2 
\, {{dx}\over{2\pi}} \leq N \Bigr\}\qquad (N>0).\eqno(1.4)$$
\noindent It is known, for example, that if $q \in \Omega_N$ then both 
$\lambda_{2j-1}$ and $\lambda_{2j}$ are asymptotically
$${{j^2}\over{4}}+\int_0^{2\pi} q(x) {{dx}\over{2\pi}}+
o(1)\qquad (j\rightarrow\infty ),\eqno(1.5)$$
\noindent and in particular that the $j^{th}$ intervals of instability 
$(\lambda_{2j-1},\lambda_{2j})$ has length
$d_j = \lambda_{2j} - \lambda_{2j-1}$ such that $(d_j)_{j=1}^\infty$ forms an $\ell^2$ sequence; see [20]
and [10]. Indeed 
the decay properties of the time-invariant sequence $\{d_j\}$ for a family of potentials
 solving (1.3) are closely related to regularity properties of those solutions. (See also
[13].)\par
\indent 
The central questions addressed in this paper concern the spectral properties of (1.1) for a
 random potential $q$. Given the link between (1.1) and (1.3) it is natural to choose $q$
 according to the Gibbs measures $\nu^\beta_N$ for the periodic KdV system which were
 introduced by Bourgain [7]. These make $(\Omega_N, \nu_N^\beta )$ into an inner regular
and Borel probability space. A typical $q$ in the support of $\nu^\beta_N$ is not differentiable, but by refining 
the classical spectral results from [20] we show that some, but not all, of the 
classical results apply.  In [3], we proved concentration inequalities and logarithmic Sobolev inequalities for the 
measures $\nu^\beta_N$. The general principle is that a
real Lipschitz function on $\Omega_N$ has very small average
oscillation with respect to $\nu_N^\beta$; see  [33, page 618]. For $q \in (\Omega_N,\nu^\beta_N)$,
 the sequence $(\lambda_j)$ is also random. Let ${\cal P}$ be the space of real periodic spectra. In
Propositions 2.3 and 3.1, we show that ${\cal P}$ is embedded in a Hilbert cube.\par
\indent We recall the definition of Hill's discriminant $\Delta$. 
Let $f_\lambda$ and $g_\lambda$ be the fundamental solutions of (1.1) such that 
$f_\lambda (0)=1$, $f'_\lambda (0)=0$; $g_\lambda (0)=0$ and $g'_\lambda (0)=1$. Then 
$$\Delta (\lambda )=f_\lambda (2\pi )+g_\lambda'(2\pi )\qquad (\lambda\in {\bf C})\eqno(1.6)$$
\noindent defines an entire function of order $1/2$, and the
 periodic spectrum is the set of zeros of $4-\Delta (\lambda )^2$. The set $\{ \lambda\in {\bf R}: \Delta (\lambda )^2\leq 4\}$ is typically an infinite union of 
closed and bounded intervals, called intervals of stability; whereas 
$\{ \lambda \in {\bf R}: \Delta (\lambda )^2>4\}$ is typically an infinite union of 
open intervals $(\lambda_{2j-1}, \lambda_{2j})$ such that all nontrivial
 solutions of (1.1) 
with $\lambda\in (\lambda_{2j-1}, \lambda_{2j})$ are unbounded. Nevertheless, the tied eigenvalues of (1.1) for the boundary conditions
$f(0)=f(2\pi )=0$ satisfy $\mu_j\in [\lambda_{2j-1}, \lambda_{2j}]$.\par
\indent In section 4 we obtain concentration
inequalities for the distribution of the tied eigenvalues $\mu_j$ and in section 5 for the 
periodic eigenvalues. Statistical information about the periodic spectrum $\{\lambda_j\}$ 
for random potentials is obtained by studying the distribution of scalar-valued random
 variables of the form $f(\{\lambda_j\})$ for suitable functions $f$, and in particular, 
linear eigenvalue statistics. This study is analogous to the results on the 
statistical properties for the eigenvalues of random unitary $n \times n$ matrices taken from the compact group $U(n)$ under 
normalized Haar measure (see, for example, [12], or [29]). The proofs in sections 4 and 5 of the current paper use similar
functional inequalities to their random matrix counterparts in [2].\par
\vskip.05in
\noindent {\bf Definition} (i) (Linear Statistics) Let $(\Omega ,\nu )$ be a probability space, and $t_j:\Omega\rightarrow
{\bf R}$ random variables for $j\in {\bf Z}$. Let $H$ be a real reproducing kernel
Hilbert space on ${\bf R}$, such that $h_t\in H$ satisfies $g(t)=\langle g,h_t\rangle_H$ for all $g\in H$ and
$t\in {\bf R}$, and $t\mapsto \langle g,h_t\rangle$ is continuous. Then for all $g\in H$, there is a sequence of
random variables $(g(t_j))_{j=-\infty}^\infty$ on $(\Omega ,\nu )$, and we define the 
corresponding linear statistics to be $\sum_{j=-m}^m 
g(t_j)$ or equivalently $\langle g, \sum_{j=-m}^m h_{t_j}\rangle$
 for $m=0, 1, \dots$.\par
\indent (ii) (Riesz basis) Let $(e_j)_{j=-\infty}^\infty$ be the usual orthonormal basis of
$\ell^2$, and let $U:\ell^2\rightarrow H$ be a bounded linear operator with
bounded inverse. Then $(Ue_j)_{j=-\infty}^\infty$ is called a Riesz basis for
$H$.\par
\vskip.05in
\indent We choose $H$ to be independent of $q$, and then select the 
$(t_j)$ depending upon $q$. The mutual dependence of the $(t_j)$ is
expressed in terms of the corresponding reproducing kernels
$h_{t_j}$ and $[\langle h_{t_j}, h_{t_k}\rangle ]$. We cannot expect that the $(h_{t_j})$ will be an orthogonal sequence;
nevertheless, we obtain conditions under which $(h_{t_j})$ is a Riesz basis. This terminology is familiar from the theory of 
wavelets and sampling theory [25].\par
\vskip.05in
\noindent {\bf Definition} (i) ($RPW(\pi )$) For $b>0$, the Paley--Wiener space  $PW(b)$ is the complex Hilbert space
 of entire functions $f$ of exponential type such that 
$$\lim\sup_{y\rightarrow\pm \infty}{{\log \vert f(iy)\vert }\over{\vert y\vert}}\leq b,
\eqno(1.7)$$
\noindent and $\int_{-\infty}^\infty \vert f(x)\vert^2 dx$ finite; so that $PW(b)$ is a closed linear 
subspace of $L^2({\bf R} ; {\bf C}).$ Let $RPW(b)$ be the closed linear subspace of $L^2({\bf R}; {\bf R})$ 
consisting of those $g\in PW(b)$ such that $\overline {g(z)} =g(\bar z)$.\par 
\indent (ii) (${\hbox{sinc}}$) The normalised cardinal sine function is ${\hbox{sinc}}(x)=(\sin \pi x)/(\pi x)$, so that
 ${\hbox{sinc}}\in RPW(\pi )$.
 Moreover, ${\hbox{sinc}}(t-s)$ is the reproducing kernel for $s\in {\bf R}$ and 
$PW(\pi )$ as in [26], so 
$$g(s)=\int_{-\infty}^\infty g(t)\,{\hbox{sinc}}(t-s)\, dt\qquad (g\in PW(\pi )).\eqno(1.8)$$
\indent (iii) (Sampling points) Introduce the points $(t_n)_{n=-\infty}^\infty$ by 
$$t_n=\cases{\sqrt{\lambda_{2n}-\lambda_0},&$ n=1,2, \dots ;$\cr 
 0,& $n=0;$\cr
-\sqrt{\lambda_{-2n}-\lambda_0},&$ n=-1,-2, \dots .$\cr} \eqno(1.9)$$
\vskip.05in
\noindent  Shannon's sampling theorem refers to
sampling on $Z=(j)_{j=-\infty}^\infty$. In the current paper, we sample on
$(t_j)_{J=-\infty}^\infty$, where $t_j$ are random variables on $(\Omega_N,
\nu_N^\beta)$. We can view the
differences $(\xi_j=t_j-j)_{j=-\infty}^\infty$ as random jitters in the
sampling points. For the linear statistics associated 
with the KdV measure, in section 5 we prove a concentration inequality which shows that the 
linear statistic $\sum_{j=-m}^m g(t_j)$ for $g\in RPW(\pi )$ is tightly concentrated about its mean value with 
Gaussian decay 
away from the mean. For typical $g\in RPW(\pi )$, the series 
$\sum_{n=-\infty}^\infty g(n)$ is not absolutely convergent, so we use the normalized 
series $\sum_{j=-\infty}^\infty (g(t_j) -g(j))$ and its partial sums. 
The mutual dependence of the $\lambda_j$ is described in terms of
the Gram matrix, 
and its generalized determinant, which is defined as follows.\par
\vskip.05in
\noindent {\bf Definition} ($\det_2$) Let $A$ be a Hilbert--Schmidt operator with norm $\Vert A\Vert_{HS}$ and
 eigenvalues 
$(\kappa_j)_{j=0}^\infty $ listed according to multiplicity. Then the Carleman determinant of $I+A$ is 
defined by the convergent infinite product
$$\det_2(I+A)=\prod_{j=0}^\infty (I+\kappa_j)e^{-\kappa_j}.\eqno(1.10)$$
\noindent We write $I+HS=\{ T\in B(\ell^2): T-I\in HS\}$ and observe that ${\cal G}=\{ T\in I+HS: \det_2 T\neq 0\}$ is an infinite
dimensional multiplicative group.\par
\vskip.05in
\indent Given $t=(t_j)$ as in (1.9) let $h_{t_j}(s) ={\hbox{sinc}}(s-t_j)$. Then 
the Gram matrix, defined by $G(t)=[{\hbox{sinc}}(t_j-t_k)]$, determines the 
properties of $(h_{t_j})_{j=-\infty}^\infty$ up to unitary equivalence. In section 6, we show that
$G(t)\in {\cal G}$ where $t=t(q)$ and $q$ belongs to a subset
of $\Omega_N$ of positive $\nu_N^\beta$ measure, and we use this fact to show that 
$(h_{t_j})_{j=-\infty}^\infty$ is a Riesz basis
for $RPW(\pi )$. In section 7, we interpret linear statistics as divisors on the spectral
curves.\par
\indent Hill's curve is the transcendental hyperelliptic curve
$${\cal E}=\bigl\{ (z, \lambda )\in {\bf C}^2: z^2=4-\Delta (\lambda
)^2\bigr\}\eqno(1.11)$$
\noindent which has branch points at the periodic eigenvalues. Choosing $\mu_j\in
(\lambda_{2j-1}, \lambda_{2j})$, we introduce a divisor $\delta$ and a differential form 
$\omega_\infty$ such that
$$\omega_\infty (\delta )=\sum_j\int_{\mu_j}^{\lambda_{2j}} {{\Delta'(\lambda )d\lambda}\over{\sqrt{\Delta (\lambda
)^2-4}}}.\eqno(1.12)$$
\noindent This formula requires careful interpretation in the case of random $q$. In section 7, we show that the
curve ${\cal E}$ has a real Jacobian ${\bf X}$, which is the real part of an infinite 
dimensional complex torus. In Proposition 7.2, we produce a map
$\Omega_N\rightarrow {\bf X}\times
{\cal P}$ which is one-to-one, up to translating the potential.\par
\indent 
 Analogously, a compact and connected Lie
group has a maximal torus
$T$ and set of cosets $G/T$ with fibration $T\times (G/T)\rightarrow G$. In the group $U(n)$ of 
complex unitary $n\times n$ matrices, the maximal torus ${\bf T}^n$ is realised as the space of diagonal unitary
 matrices. The eigenvalue map $U\mapsto {\hbox{Diag}}(e^{i\theta_j})$ induces a probability
 measure $\mu_{W}$ on ${\bf T}^n$ from Haar measure on $U(n)$, as expressed by the Weyl denominator 
formula [14].\par
\indent The table below expresses the analogy between linear statistics for the KdV measure and 
linear statistics for $U(n)$. Let $f\in L^1({\bf T}; {\bf C})$ have Fourier coefficients $\hat f_j$, and form the Toeplitz determinant 
$D_n(f)=\det  [\hat f_{j-k}]_{j,k=0}^{n-1}$. In particular, suppose that $g\in L^1({\bf T}; {\bf C})$ 
belongs to $H^{1/2}({\bf T})$ so that 
$\sum_{k\in {\bf Z}}\vert k\vert \vert \hat g(k)\vert^2$ finite and let $f=e^g$.
 Then Johansson [12] proved a strong Szeg\"o theorem, which describes the asymptotic growth 
of $D_n(e^g)$ as $n\rightarrow\infty$ in terms of a quadratic form in $g$, 
and established a central limit theorem for the linear statistics. Our concentration results are quite
analogous. 
$$\matrix{ {}& {\hbox{Hill's curve}}& {\hbox{Unitary group}}\cr
{\hbox{Probability space}}& \Omega_N=\{ q\in L^2({\bf T}); \Vert q\Vert^2_2\leq N\}& U(n)\cr
{\hbox{Probability measure}}& \nu_N^\beta& \mu_{W}\cr
{\hbox{Eigenvalues}}& {\bf X}=\{ (\mu_j)_{j=1}^\infty\}&
{\bf T}^n =\{ (e^{i\theta_j})_{j=1}^n\}\cr
{\hbox{Test functions}}& g(z)\in RPW(\pi )& g\in H^{1/2}({\bf T})\cr
{\hbox{Linear statistic}}& \sum_{j=-m}^m (g(t_j)-g(j))& \sum_{j=1}^n 
g(e^{i\theta_j})-n\int_{\bf T} g(e^{i\theta}) {{d\theta}\over{2\pi}}\cr
{}&{}&{}\cr
{\hbox{Quadratic form}}& 2^{-1}\sum_{j=-m}^m \vert g'(t_j)\vert^2 & 
\int\!\!\!\int_{{\bf T}^2} -\log\vert e^{i\theta}-e^{i\phi}\vert\, g'(e^{i\theta}) g'(e^{i\phi}) {{d\theta}\over{2\pi}}{{d\phi}\over {2\pi}}\cr}$$
\vfill
\eject
\noindent {\bf 2 Concentration of measures}\par
\vskip.05in
\indent In this section we prove some concentration inequalities for 
three measures. We introduce the infinite dimensional phase space of $u:{\bf T}\times (0, \infty )\rightarrow {\bf R}$ 
such that $u(\cdot ,t)\in L^2({\bf T})$ for all $t>0$, and the densely-defined Hamiltonian
$$H(u)={{1}\over{2}}\int_{\bf T} \Bigl( 
{{\partial u}\over{\partial x}}(x,t)\Bigr)^2 {{dx}\over{2\pi}}-{{\beta}\over{6}}\int_{\bf T}
u(x,t)^3 {{dx}\over{2\pi}},\eqno(2.1)$$
\noindent where $\beta >0$ is the inverse temperature. The canonical equation of motion reduces to the 
KdV equation (1.3) which gives a nonlinear evolution on $L^2({\bf T})$ such that 
$H(u)$ and $\int_{\bf T} u(x,t)^2{{dx}\over{2\pi}}$ are invariant with respect to $t$. 
The relevance to (1.1) is that the flow generated by the KdV equation 
${{\partial q}\over{\partial t}}=3q{{\partial q}\over{\partial x}}-{{1}\over{2}}{{\partial^3q}\over{\partial x^3}}$ 
preserves the periodic spectrum.\par
For some normalizing constant $Z_N(\beta )$, there exists a probability measure
$$\nu_N^\beta (d\phi )=Z_N(\beta )^{-1}
{\bf I}_{\Omega_N}(\phi )e^{-H(\phi )} \prod_{e^{ix}\in {\bf T}} d\phi (x)\eqno(2.2)$$
\noindent which is Radon, in the sense of being inner-regular and defined on the 
$\sigma$-algebra generated by the Borel sets in $L^2$. Bourgain [7] constructed this probability measure using random 
Fourier series, as follows. Let $(\gamma_n)_{n=-\infty}^\infty$ be mutually independent standard 
Gaussian random variables on some probability space $(\Omega , {\bf P})$, and let $\gamma$ be 
the Gaussian probability measure on $L^2({\bf T}; {\bf R})$ that is induced by the map
$$\omega\mapsto \phi_\omega (x) =\gamma_0+\sum_{n=-\infty}^{-1} \gamma_n {{\sin nx}\over{n}}+
\sum_{n=1}^\infty \gamma_n {{\cos nx}\over{n}}.\eqno(2.3)$$
\noindent This gives an interpretation to the product factor in (2.2).\par
\indent Let the Fourier expansion  of the potential be
 $q(x) =(1/2)a_0+\sum_{n=1}^\infty (a_n\cos nx +b_n\sin nx)$, 
where the Fourier coefficients are real random variables on  $(\Omega_N, \nu_N^\beta )$
 satisfying 
$(a_0/2)^2+2^{-1}\sum_{n=1}^\infty (a_n^2+b_n^2)\leq N$. \par
\vskip.05in 
\noindent {\bf Lemma 2.1} {\sl Let $A=[\alpha_{jk}]$ be a trace class and real symmetric matrix
 on $\ell^2$ with eigenvalues $(\sigma_j)_{j=0}^\infty$ which are listed
 according to multiplicity, and such that $(\vert \sigma_j\vert)_{j=0}^\infty$ is decreasing,
and let $\gamma_j$ be mutually independent normalized Gaussian random variables on $L^2$.
 Then for $0<\varepsilon <1/2$ and $K=\pi /2\sqrt {1-2\varepsilon}$, the random variable  $Q=\sum_{j,k}\alpha_{jk}\gamma_j\gamma_k$ satisfies}
$$\int_{\Omega_N}e^{sQ}\nu_N^0(dq)\leq 
(1-e^{-\varepsilon N}K)^{-1}\prod_{j=0}^\infty (1-2s\sigma_j)^{-1/2}\qquad (N>\varepsilon^{-1}\log K,
 \quad 2s\vert\sigma_0\vert <1).\eqno(2.4)$$ 
\vskip.05in
\noindent {\bf Proof.} The sequence $(\vert\sigma_j\vert)_{j=0}^\infty$ gives the eigenvalues of the positive square root $(A^2)^{1/2}$, hence gives a summable sequence, so the product on the right-hand side is absolutely convergent. There exists a real orthogonal matrix on $\ell^2$ such that $U^\dagger AU$ is the real diagonal matrix ${\hbox{diag}}(\sigma_j)_{j=0}^\infty$.\par
\indent By elementary results regarding Gaussian random variables, we can compute
$$\int_{L^2} e^{sQ}\gamma (dq) =
\prod_{j=0}^\infty \int_{-\infty}^\infty e^{s\sigma_j -x^2/2} {{dx}\over{\sqrt{2\pi}}}
=\prod_{j=0}^\infty (1-2s\sigma_j)^{-1/2}.\eqno(2.5)$$
\noindent When $\sigma_j=1/j^2$ we can compute this product explicitly, and so for $0<\varepsilon <1/2$ we obtain 
$$\int_{L^2} e^{\varepsilon \Vert q\Vert^2} \gamma (dq)= {{1}\over{\sqrt{1-2\varepsilon}}}{{\sqrt{\varepsilon \pi}}\over{\sin \sqrt{\varepsilon \pi}}}
=K_\varepsilon ;\eqno(2.6)$$
\noindent so by Chebyshev's inequality, $\gamma (\Omega_N^c)\leq K_\varepsilon e^{-\varepsilon N}$; hence  the normalizing constant for $\nu_N^0$ on $\Omega_N$ satisfies
$$\int_{\Omega_N} \gamma (dq)> \max\{1-K_\varepsilon e^{-\varepsilon N}, 0\}.\eqno(2.7)$$
\rightline{$\square$}\par
\indent The following transportation cost bounds the distance between $\nu^0_N$ and $\nu_N^\beta$ in the 
Wasserstein metric; see [32] for a general discussion.\par
\vskip.05in
\noindent {\bf  Lemma 2.2} {\sl For all $\beta, N>0$, there exists $\kappa (\beta ,N)$ such that
if $f:\Omega_N\rightarrow {\bf R}$ is any Lipschitz function such that 
$\vert f(p)-f(q)\vert \leq L\Vert p-q\Vert_{L^2}$ for all $p,q\in \Omega_N$, then} 
$$\Bigl\vert\int_{\Omega_N} f(q)\nu_N^0(dq)-\int_{\Omega_N} f(q)\nu_N^\beta (dq)\Bigr\vert 
\leq L\beta N\kappa.\eqno(2.8)$$
\noindent {\sl A valid choice of $\kappa$ is
$ \kappa =e^{C\beta^{5/2}N^{9/4}}$ where $C$ is some absolute constant 
taken from [3].}\par
\vskip.05in
\noindent {\bf Proof.} By the Kantorovich--Rubinstein duality theorem
[32, p 34], the quantity (2.8) in the Lemma may be expressed as
$$ \Bigl\vert \int_{\Omega_N} f(q)\nu_N^0(dq)-\int_{\Omega_N} f(q)\nu_N^\beta (dq)\Bigr\vert
=L\inf_\pi\Bigl\{ \int\!\!\!\int_{\Omega_N\times \Omega_N}\Vert p-q\Vert_{L^2} \pi (dpdq)\Bigr\},
\eqno(2.9)$$
\noindent where the infimum is taken over all the Radon probability measures on 
$\Omega_N\times \Omega_N$ that have marginals $\nu_N^0$ and $\nu_N^\beta$. 
The infimum on the right-hand side of (2.9) defines the transportation cost for the cost 
function 
$\Vert p-q\Vert_{L^2}$.\par
\indent Next we note that $\nu_N^0$ is symmetrical in distribution with respect to $a_n\leftrightarrow -a_n$ and
 $b_n\leftrightarrow -b_n$ so some expectations are easy to compute. Note also that
$F(q)=-\beta \int_{\bf T} q(x)^3dx/(2\pi)$ is Lipschitz on $\Omega_N$ with $L=\beta N$, 
and $\int_{\Omega_N} F(q)\nu_N^0(dq)=0$. 

We bound the left-hand side of (2.9) in terms of the relative
entropy, where the relative entropy of $\nu_N^0$ with respect to $\nu_N^\beta$ satisfies
$$\eqalignno{{\hbox{Ent}}(\nu_N^0\mid \nu_N^\beta )&=\int_{\Omega_N}\log {{d\nu_N^0}\over{d\nu_N^\beta}} 
d\nu_N^0\cr
&=\int_{\Omega_N} (-F(q)+\log Z_N(\beta ))\nu_N^0(dq)\cr
&=\log Z_N(\beta ).&(2.10)\cr}$$
\noindent By the concentration inequality for Gaussians, the normalizing constant for $\nu_N^\beta$ satisfies
$$Z_N^\beta =\int_{\Omega_N}e^{F(q)} \nu_N^0(dq)\leq e^{\beta^2N^2}.\eqno(2.11)$$
\noindent The concentration inequality for $\nu_N^\beta$ which was proved
 in [3] gives the upper bound 
$$(2.9)\leq L\kappa \bigl( 2{\hbox{Ent}}(\nu_N^0\mid \nu_N^\beta )\bigr)^{1/2}.
\eqno(2.12)$$
\indent Hence by (2.11) and (2.12) we have the bound $(2.9)\leq \sqrt{2} L\kappa \beta N.$ This estimate 
improves as $\beta\rightarrow 0$, 
since $\nu_N^\beta$ converges in Wasserstein metric to $\nu_N^0$.          {$\square$}\par
\vskip.05in
\noindent {\bf Proposition 2.3} {\sl (i) There exists $N_1>0$ such that for $0<N<N_1$ 
there exist positive constants $c(N, \beta )$ and $C(N, \beta )$ such that the mean 
length of the $j^{th}$ interval of instability satisfies}
$$ {{c(N, \beta )}\over{j}}\leq \int_{\Omega_N} (\lambda_{2j}-\lambda_{2j-1}) 
\nu_N^\beta (dq)\leq {{C(N, \beta )}\over{j}}\qquad (j=1, 2, \dots ).
\eqno(2.13)$$
\noindent {\sl (ii) With probability one with respect to $\nu_N^\beta$, all the
periodic eigenvalues are simple, so $\lambda_j<\lambda_{j+1}$ for all
$j$.}\par
\vskip.05in
\noindent {\bf Proof.} (i) Erdelyi [8] showed that there exists $N_1$ such that for all $q\in \Omega_N$ with $0<N<N_1$, the length of the $j^{th}$ spectral gap is equivalent to the $j^{th}$ Fourier coefficient of $q$, so
$$\lambda_{2j}-\lambda_{2j-1}\asymp cN\Bigl\vert \int_0^{2\pi} q(x)e^{-ijx} {{dx}\over{2\pi}}
\Bigr\vert\qquad (j\rightarrow\infty ).\eqno(2.14)$$
\noindent From (2.3), we deduce that 
$$\int_{\Omega_N} (\lambda_{2j}-\lambda_{2j-1})^2 \nu_N^0 (dq)\leq {{CN^2}\over{j^2}}\int_{\Omega_N} 
\gamma_j^2 \nu_N^0(dq ),\eqno(2.15)$$
\noindent and so we can use the Cauchy--Schwarz inequality to obtain an upper bound involving the measure $\nu_N^\beta$, namely
$$\eqalignno{\Bigl( \int_{\Omega_N} (\lambda_{2j}-\lambda_{2j-1})\nu_N^\beta (dq)\Bigr)^2&\leq \int_{\Omega_N} (\lambda_{2j}-\lambda_{2j-1})^2 \nu_N^0(dq) \int_{\Omega_N}\Bigl({{d\nu_N^{\beta}}\over{d\nu_N^0}}\Bigr)^2 \nu_N^0(dq) \cr
&\leq {{C(N, \beta )}\over{j^2}}& (2.16)}$$
\noindent where the constants involve the normalizing constant for the Gibbs measure. \par
\indent To prove the lower bound in (2.13), we introduce the events 
$A_N(j)=\{ \omega : \vert \hat q(j)\vert^2 \leq N/2\}$ and 
$B_N(j)=\{ \omega : \sum_{k:k\neq \pm j}\vert q\vert^2\leq N/2\}$ which are independent
 in $(L^2, \gamma )$ and satisfy $A_N(j)\cap B_N(j)\subseteq \Omega_N$. Hence there exists $c_N>0$ such that 
$$\eqalignno{ {{c_N}\over{\vert j\vert^{1/2}}}&\leq \int_{A_N(j)}\vert \hat q(j)\vert^{1/2}d\gamma 
\int_{B_N(j)} \Bigl( \sum_{k: k\neq \pm j} \vert \hat q(k)\vert^2\Bigr)^{1/4}d\gamma\cr
&\leq \int_{\Omega_N} \vert \hat q(j)\vert^{1/2}\Bigl( \sum_{k}\vert \hat q(k)\vert^2 
\Bigr)^{1/4}
 d\gamma ,&(2.17)}$$
\noindent wherein $(\sum_{j=-\infty}^\infty \vert \hat
q(j)\vert^2)^{1/4}\leq
N^{1/4}$. Now by (2.3) and (2.14) we have
$$\eqalignno{ (2.17)&\leq Z_N(\beta )\int_{\Omega_N} \sqrt{\lambda_{2j}-\lambda_{2j-1}} d\nu_N^0\cr
&\leq Z_N(\beta)\Bigl( \int_{\Omega_N} (\lambda_{2j}-\lambda_{2j-1})d\nu_N^\beta \Bigr)^{1/2}
\Bigl( \int_{\Omega_N}\Bigl( {{d\nu_N^0}\over{d\nu_N^\beta}}\Bigr)^{2} d\nu_N^0\Bigr)^{1/2}.&
(2.18)\cr}$$
\indent (ii) It suffices to show that
$\sum_{j=1}^\infty j^{-2}\vert \log(\lambda_{2j}-\lambda_{2j-1})\vert$
is integrable with respect to $\nu_N^\beta.$ The required
estimate is (2.16), with $\vert\log
(\lambda_{2j}-\lambda_{2j-1})\vert$ replacing
$\lambda_{2j}-\lambda_{2j-1}$. {$\square$}\par
\noindent {\bf Remark 2.4} An unwelcome consequence of (2.13) is that 
$\sum_{j=1}^\infty (\lambda_{2j}-\lambda_{2j-1})\lambda_{2j}^k$ diverges for $k\geq 0$ 
with probability one, so the polynomial 
approximation arguments from section 10 of [22] become inapplicable. In
Proposition 7.2, we consider the inverse spectral problem for $q\in (\Omega_N, \nu_N^\beta )$.\par 
\vskip.05in

\noindent {\bf 3 Fluctuations and sampling of the periodic eigenvalues}\par
\vskip.05in
\noindent In this section we analyse the discriminant of (1.6) and obtain detailed 
information about the fluctuations of the periodic eigenvalues by using the 
concentration results of section 2. Loosely speaking, we show
that the intervals of instability $[\lambda_{2j-1}, \lambda_{2j}]$
are random with mean length of order $1/j$ and centre of order
$j^2/4$. First we obtain a version of Borg's estimates from [20].\par
\vskip.05in
\noindent {\bf Proposition 3.1} {\sl For all 
$q\in \Omega_N$ and all $\delta >0$, the periodic eigenvalues of 
(1.1) satisfy}
$$\bigl\vert \sqrt{\lambda_{2n}}-{{n}\over{2}}\bigr\vert\leq 
{{\sqrt{\beta N}}\over{n^{1-\delta}}}\qquad (n>( \beta N )^{1/\delta}
 \kappa^{1/\delta})\eqno(3.1)$$
\noindent {\sl so $\{ (\sqrt{\lambda_{4n}-\lambda_0}-n)_{n=1}^\infty :q\in
\Omega_N\}$ is relatively compact in the norm topology of $L^2$.}\par 
\vskip.05in
\noindent {\bf Proof.} For $0<\eta\leq 2\pi$,
 let $V,W:L^\infty [0,\eta ]\rightarrow  L^\infty [0,\eta ]$, where $V=V(q)$ and $W=W(q)$, be the bounded linear operators 
$$Vf(x)=\int_0^x \cos (\sqrt{\lambda } (x-t))\, q(t)f(t)dt\qquad (f\in L^\infty [0,
\eta ]),\eqno(3.2)$$
$$Wg(x)=\int_0^x \sin (\sqrt{\lambda} (x-t))\, q(t)g(t)dt\qquad (g\in L^\infty [0,
\eta ]).\eqno(3.3)$$
\noindent These integral operators have operator norms that satisfy
$$\Vert V\Vert_{op}, \Vert W\Vert_{op}\leq \cosh (\eta Y)\int_0^\eta \vert q(t)\vert dt\leq \cosh (\eta Y)
 \sqrt{\eta N};\eqno(3.4)$$
\noindent for all $q\in \Omega_N$. This estimate is uniform over the probability space 
$\Omega_N$ and not random, and enables us to obtain non random error terms in the
 following. \par
\indent The fundamental solutions of Hill's equation satisfy
$$f_\lambda (x)=\cos\sqrt{\lambda} x +{{1}\over{\sqrt{\lambda}}}Wf_\lambda (x),\eqno(3.5)$$
$$g_\lambda (x)={{\sin\sqrt{\lambda}x}\over{\sqrt{\lambda}}}+{{1}\over{\sqrt{\lambda }}}Wg_\lambda
(x);\eqno(3.6)$$
\noindent so the discriminant may be expressed as the series
$$\eqalignno{ \Delta (\lambda )&=2\cos 2\sqrt{\lambda }\pi +{{1}\over{\sqrt{\lambda}}}\int_0^{2\pi}
 (q(t)+q(2\pi -t)) \sin \sqrt{\lambda} (2\pi -t)
 \cos \sqrt{\lambda} t \,dt\cr
&+ {{1}\over{\lambda}}\bigl( W^2(\cos\sqrt{\lambda} t) (2\pi )+VW(\sin \sqrt{\lambda}t)(2\pi )\bigr)\cr
&+{{1}\over{\lambda^{3/2}}}\bigl( W^3 (\cos\sqrt{\lambda} t)(2\pi )+VW^2(\sin \sqrt{\lambda }t)(2\pi )\bigr)
+O( {{1}\over{\lambda^2}} );&(3.7)\cr}$$
 \noindent in which the coefficient of $\lambda^{-k/2}$ consists of $2^k$ terms, each 
of which is a product of $k$ factors of $V$ and $W$; hence the series  converges
 for $\lambda >4N.$ We consider the terms on the right hand side of (3.7) in turn: $2\cos 2\sqrt{\lambda} \pi$ is the main term and is independent of $q$; the next term is random, and reduces by  Fourier analysis to $(\pi a_0\sin 2\sqrt{\lambda}\pi )/\sqrt{\lambda};$ then the coefficient of $1/\lambda$ is
$$Q=\int_0^{2\pi} q(t)\int_0^t \sin \sqrt{\lambda} (t-x)\,\sin \sqrt{\lambda}(2\pi -t+x)\, q(x)\, dxdt,\eqno(3.8)$$
\noindent which is a real quadratic form in $q$, and hence may be expressed as the symmetrical expression
involving the symmetric kernel
$$K(x,y)=\cases{ \sin \sqrt{\lambda} (y-x)\,\sin \sqrt{\lambda}(2\pi -y+x)& for $0\leq x\leq y\leq
2\pi$;\cr
\sin \sqrt{\lambda} (x-y)\,\sin \sqrt{\lambda}(2\pi -x+y)& for $0\leq y\leq x\leq 2\pi$,\cr}$$
\noindent this defines a trace class linear operator on $L^2[0,2\pi ]$, as one 
checks by considering the operation of $K$ on  $(e^{inx})_{n=-\infty}^\infty$. Passing to 
the Fourier coefficients, we deduce that 
$$Q={{1}\over{8}}a^2_0\Bigl(
{{\pi\sin 2\pi\sqrt{\lambda}}\over{\sqrt{\lambda}}}-2\pi^2\cos 2\pi
\sqrt{\lambda}\Bigr) -\pi\sqrt{\lambda} \sin
2\pi\sqrt{\lambda}\sum_{n=1}^\infty
{{a_n^2+b_n^2}\over{n^2-4\lambda}}.$$ 
\indent Hence the discriminant satisfies
$$2-\Delta (\lambda )=4\sin^2 \pi\sqrt{\lambda}-{{2\pi a_0\sin\pi\sqrt{\lambda}}\over{\sqrt{\lambda}}}
 \cos \pi\sqrt{\lambda } +O\Bigl( {{\beta N}\over{\lambda }}\Bigr)\eqno(3.9)$$
\noindent for large $\lambda$. From the product formula for $\sin$,
the factors are holomorphic functions for $\Re \lambda >1$. We let $\sqrt{\lambda}=2n+\tau+i\sigma$ where $\vert \tau+i\sigma\vert =\sqrt{\beta N}/n^{1-\delta}$  and apply Rouch\'e's theorem to 
compare $2-\Delta (\lambda )$ with $4\sin^2\pi\sqrt{\lambda}$. Note that $\sin^2\pi\sqrt\lambda$ has a
double zero at $\lambda =n^2$, and 
$$\bigl\vert 4\sin^2\pi\sqrt{\lambda}\bigr\vert=4\sin^2\pi \tau+4\sinh^2\pi\sigma 
>4(\tau^2+\sigma^2),
\eqno(3.10)$$
\noindent which is greater than $\beta N/\vert \lambda\vert$. This gives a pair of zeros of
$2-\Delta (\lambda )$, namely eigenvalues $\lambda_{4n-1}$
and $\lambda_{4n}$ of the principal series.\par
\indent In the analysis of periodic spectra in [20] page 24, the authors assume
that $q$ has mean zero. In our case, however, the value of $a_0(q)$ is random, we now consider its influence on
the spectrum. We have a coarse estimate on $a_0$, namely $\vert a_0\vert \leq \sqrt{N}$, so
$$\Bigl\vert {{\pi a_0}\over{\sqrt{\lambda}}}\Bigr\vert\leq  {{\pi \sqrt {N}}\over{2n}}\leq {{\pi}\over{n^\delta \sqrt{\beta}}} \vert \tau+i\sigma \vert.\eqno(3.11)$$
\noindent which is an appropriate estimate when $\beta>0$ is not small. \par
\indent  Also, note that $a_0=\int_0^{2\pi} q(x)dx/\pi$ is a Lipschitz function of $q\in \Omega_N$ with $L=1$, so by Lemma 2.2, we have
$$\Bigl\vert\int_{\Omega_N} a_0(q) \nu_N^\beta (dq)-\int_{\Omega_N} a_0(q)\nu_N^0 (dq)\Bigr\vert\leq N\beta\kappa ,\eqno(3.12)$$
\noindent \noindent where $ \kappa =e^{C\beta^{5/2}N^{9/4}}$; so 
this estimate is advantageous when $\beta>0$ is small. Then by the concentration inequality Corollary 2 of [3], we have
$$\nu_N^\beta \Bigl\{ q\in \Omega_N: \bigl\vert a_0(q)-\int a_0d\nu_N^\beta\bigr\vert
>t\Bigr\}\leq 2e^{-t^2/4\kappa}
\eqno(3.13)$$
\noindent and by choosing $0<\varepsilon <\delta $ and $t=2\beta\sqrt{\kappa }n^\varepsilon $ we deduce that
$$\nu_N^\beta \Bigl\{ q\in \Omega_N: \bigl\vert a_0(q)-\int a_0d\nu_N^\beta \bigr\vert >2\beta 
\sqrt{\kappa} n^\varepsilon\Bigr\} \leq 2e^{-\beta n^\varepsilon }.\eqno(3.14)$$
\noindent By the Borel--Cantelli Lemma, we deduce that 
$$\nu_N^\beta \Bigl\{ q\in \Omega_N: \vert a_0(q)\vert  \leq \beta N\kappa +2\beta \sqrt{\kappa }n^\varepsilon\,\,\,
{\hbox{for all but finitely many}}\,\,\, n\Bigr\} =1.\eqno( 3.15)$$
\noindent We deduce that, on a set of full $\nu_N^\beta$ measure, the inequality 
$${{\pi \vert a_0\vert}\over{\sqrt{ \vert \lambda\vert}}}\leq 
\Bigl( {{\kappa \sqrt{\beta N}}\over{n^\delta}}+ {{2\sqrt{\beta} \kappa}\over 
{n^{\delta -\varepsilon}}}\Bigr) \vert \tau +i\sigma\vert \eqno(3.16)$$
\noindent holds for all but finitely many $n$.\par
\indent It follows that on the curve ${\cal C}_n$ defined by $\sqrt{\lambda} =2n+\tau +i\sigma$, we have
$$ \bigl\vert 4\sin^2\pi\sqrt{\lambda}\bigr\vert>\bigl\vert (2-\Delta (\lambda ))-4\sin^2\pi\sqrt{\lambda}\bigr\vert\eqno(3.17)$$
\noindent so that $2-\Delta (\lambda )$ and $4\sin^2\pi \sqrt{\lambda}$ both have two 
zeros inside ${\cal C}_n$, and by the basic theory of Hill's equation, the zeros are 
real. Likewise, the discriminant satisfies
$$2+\Delta (\lambda )=4\cos^2 \pi\sqrt\lambda +{{2\pi a_0\sin\pi\sqrt{\lambda}}\over{\sqrt{\lambda}}} \cos \pi\sqrt{\lambda} +O\Bigl( {{\beta N}\over{\vert\lambda\vert}}\Bigr).\eqno(3.18)$$
\noindent hence $2+\Delta (\lambda )$ has a pair of real zeros $\lambda_{4n-3}$ and 
$\lambda_{4n-2}$ of the complementary series when 
$\sqrt{\lambda}$ is close to $(2n+1)/2$.\par
\indent Finally, observe that for all $C>0$ and $0<\delta <1/2$, the subset $\{
(\xi_n)_{n=1}^\infty : \vert\xi_n\vert\leq C/n^{1-\delta} , n=1, 2, \dots \}$  of
$\ell^2$ is a compact Hilbert cube. {$\square$}\par
\vskip.05in
\noindent {\bf Definition} (Cartwright class) The Cartwright class [15, 24] is the space of the entire functions $f$ such that 
$$\int_{-\infty}^\infty {{\log_+\vert f(x)\vert }\over{1+x^2}}\, dx<\infty .\eqno(3.19)$$
\vskip.05in
\noindent {\bf Proposition 3.2} {\sl For $q\in (\Omega_N, \nu_N^\beta )$, let $(t_j)_{j=-\infty}^\infty$ be as in (1.9), where the $\lambda_j$ are the
periodic eigenvalues of (1.1). Then for all $q$ in a subset of $\Omega_N$ with measure one with
respect to $\nu_N^\beta$ there exists $b_0>0$ such that, if $f$ is an entire functions 
of exponential type $b<b_0$ and 
$$\sum_{j\in {\bf Z}} {{\log_+\vert f(t_j)\vert }\over{1+t_j^2}}<\infty ,
\eqno(3.20)$$
\noindent then $f$ belongs to the Cartwright class.}
\vskip.05in
\noindent {\bf Proof.} In his study of almost periodic functions, Bohr called a 
discrete subset $\Lambda$ of
 ${\bf R}$ relatively dense if there exists $L>0$ such that $\Lambda$ intersects 
all intervals $[a, a+L]$ for $a\in {\bf R}$. Pedersen [27] refined this, by saying that a
 discrete subset $T$ of ${\bf R}$ is $h$-dense if there exists $h>0$ such that, outside of some bounded set, 
$T$ intersects all bounded intervals of length $h$. By Proposition 3.1, 
the set $T=\{t_j: j\in {\bf Z}\}$ is $3/2$ dense, symmetric and separated,
so that there exists $\delta >0$, depending on $q$,  such that $\vert t_k-t_m\vert \geq \delta$ for all $k\neq m$. 
 Hence by Theorem 5.1 of [27], all entire functions $f$ of exponential type $b$ for $0\leq b\leq b_0$ 
and such that (3.20) holds also have finite logarithmic integral (3.19), hence belong to 
Cartwright's class. {$\square$}\par
\indent The connection between the logarithmic integral and exponential 
sums is subtle, so we refer the reader to [15]. We summarize some standard results regarding (1) sampling and (2) interpolation for the
Paley--Wiener space. Suppose that $t=(t_j)_{j=-\infty}^\infty$ is separated (uniformly discrete in [26, p 219]),
 so that there exists $\delta >0$ such that $\vert t_k-t_m\vert \geq \delta$ for all $k\neq m$. 
Then the sampling map $S_t:PW(\pi )\rightarrow \ell^2:$ $f\mapsto (f(t_j))_{j=-\infty}^\infty$ 
defines a bounded linear operator, hence the adjoint map is also bounded, where
 $$S_t^\dagger :\ell^2\rightarrow PW(\pi ):\quad S_t^\dagger (a_j)=\sum_{j=-\infty}^
\infty  a_j {\hbox{sinc}}(s-t_j).\eqno(3.21)$$
\indent (1) If $S_t$ is an embedding, so there exists $B>0$ such that $\Vert S_tf\Vert\geq B\Vert f\Vert$ for all
 $f\in PW(\pi)$, then
$({\hbox{sinc}}(s-t_j))_{j=-\infty}^\infty$ is a frame in 
$PW(\pi )$, and we say that $(t_j)$ is a sampling sequence.\par
\indent (2) If $S_t^\dagger$ is an
 embedding, or equivalently $S_t: PW(\pi )\rightarrow \ell^2$ is surjective, 
 then $({\hbox{sinc}}(s-t_j))_{j=-\infty}^\infty$ is a Riesz basis for 
its closed linear span. Then we say that $(t_j)_{j=-\infty}^\infty$ is 
an interpolating sequence.\par
\indent Seip has produced an example of a frame $({\hbox{sinc}}(s-s_j))_{j=-\infty}^\infty$ 
for $PW(\pi )$ such that no subsequence gives a Riesz 
basis for $PW(\pi )$; see [25, 26].
We write 
$$Z+\ell^2 =\{ (t_j)_{j=-\infty}^\infty :
(t_j-j)_{j=-\infty}^\infty \in \ell^2\}\eqno(3.22)$$
\noindent  which is a complete metric space for the
norm $\ell^2$.\par
\vskip.05in
\noindent {\bf Lemma 3.3} {\sl For $t=(t_k)_{k=-\infty}^\infty$, let $U_t:L^2[-\pi ,\pi ]\rightarrow L^2[-\pi
,\pi ]$ be the linear operator $U_t:e^{ikx}\mapsto e^{it_kx}$ for all $k\in {\bf Z}$.
Then $U_t-I\in HS$ for all $t\in Z+\ell^2$, and} 
$$\Vert U_t-U_s\Vert_{HS}\leq {{2\pi}\over{\sqrt{3}}}\Vert
t-s\Vert_{\ell^2}\qquad (s,t\in Z+\ell^2).\eqno(3.23)$$   
\vskip.05in
\noindent {\bf Proof.} Clearly $U_t$ is determined by extending linearly its operation on the
 complete orthonormal basis $(e^{ikx})_{k=-\infty}^\infty$ of 
$L^2([-\pi
,\pi ]; dx/2\pi )$. By elementary estimates we have
$$\eqalignno{\bigl\Vert U_t(e^{ikx})-e^{ikx}\bigl\Vert_{L^2}^2&=\int_{-\pi}^\pi
\bigl\vert e^{it_kx} -e^{ikx}\bigr\vert^2 {{dx}\over{2\pi}}\cr
&\leq {{4\pi^2}\over{3}}(t_k-k)^2,&(3.24)\cr}$$
\noindent so $\sum_{k=-\infty}^\infty \Vert U_t(e^{ikx})-e^{ikx}\Vert^2_{L^2}$
converges; hence $U_t-I$ is Hilbert--Schmidt. Likewise we obtain $\Vert
U_t(e^{ikx})-U_s(e^{ikx})\Vert^2 \leq 4\pi^2 (t_k-s_k)^2/3$, which is a
summable sequence whenever $t-s\in \ell^2$; hence $t\mapsto U_t$ is
Lipschitz. {$\square$}\par
\vskip.05in
\noindent {\bf Proposition 3.4} {\sl Suppose that ${(2\pi
/\sqrt{3})}\Vert (t_k-k)_{k=-\infty}^\infty \Vert_{\ell^2}<1.$ Then} $({\hbox{sinc}}(x-t_j))_{j=-\infty}^\infty$ {\sl is a Riesz
basis for $RPW(\pi)$, so there exist $A,B>0$ such that for all for all $g\in RPW(\pi)$ 
there exists a unique $(a_j)_{j=-\infty}^\infty\in\ell^2$ such that} $g(x)=\sum_{j=-\infty}^\infty
a_j{\hbox{sinc}}(x-t_j)$ {\sl converges in $L^2$ and} 
$$A\int_{-\infty}^\infty \vert g(t)\vert^2 dt\leq \sum_{j=-\infty}^\infty 
\vert g(t_j)\vert^2 \leq B\int_{-\infty}^\infty \vert g(t)\vert^2 dt
\qquad (g\in PW(\pi )). \eqno(3.25)$$
\vskip.05in
\noindent {\bf Proof.} Let $F:\ell^2\rightarrow L^2[-\pi ,\pi ]$ be the
usual Fourier map $(a_j)_{j=-\infty}^\infty \mapsto \sum_{j=-\infty}^\infty
a_je^{ijx}$, and the Fourier transform ${\cal
F}f(x)=\int_{-\pi}^\pi e^{-ixu} f(u)\, du/(2\pi )$, so that $F$ and
$\sqrt{2\pi} {\cal F}$ are unitary operators. Then ${\cal F}U_t: e^{ikx}\mapsto
{\hbox{sinc}}(x-t_k)$, and extending this linearly, we deduce that $S_t^\dagger
={\cal F}U_tF$ is a bounded linear operator.\par
\indent In particular, with $Z=(j)_{j=-\infty}^\infty$ we obtain 
Shannon's map $S_Z:PW(\pi )\rightarrow \ell^2:$ $f\mapsto
(f(j))_{j=-\infty}^\infty$, which is unitary [26, p. 209]. By Lemma 3.3, we have
$$\eqalignno{\bigl\Vert S_t-S_Z\bigr\Vert_{op}&\leq \bigl\Vert S_t^\dagger
-S_Z^\dagger\bigr\Vert_{HS}\cr
&\leq {{{2\pi}\over{\sqrt{3}}}}\Bigl(\sum_{j=-\infty}^\infty
\bigl(t_j-j\bigr)^2\Bigr)^{1/2}&(3.26)\cr}$$
\noindent so that 
$$\bigl\Vert S_tf\bigr\Vert_{\ell^2}\geq \Vert f\Vert_{L^2}- (2\pi
/\sqrt{3})\bigl\Vert (t_k-k)_{k=-\infty}^\infty \bigr\Vert_{\ell^2}\Vert f\Vert_{L^2},\eqno(3.27)$$
\noindent and likewise
$$\bigl\Vert S^\dagger_t(a_j)\bigr\Vert_{\ell^2}\geq\bigl\Vert (a_j)_{j=-\infty}^\infty
 \bigr\Vert_{\ell^2}- (2\pi
/\sqrt{3})\bigl\Vert (t_k-k)_{k=-\infty}^\infty \bigr\Vert_{\ell^2}
\bigl\Vert (a_j)_{j=-\infty}^\infty \bigr\Vert_{\ell^2},\eqno(3.28)$$
\noindent so that for $A=1-(2\pi
/\sqrt{3})\Vert (t_k-k)\Vert_{\ell^2}>0$ both the linear operators $S_t$ and $S_t^\dagger$ are embeddings.
Hence $({\hbox{sinc}}(x-t_j))_{j=-\infty}^\infty$ is a Riesz basis for its closed linear span, which
is all of $PW(\pi )$. ${\square}$\par

\vskip.05in

\noindent {\bf 4. Concentration of measure on the tied spectrum}\par
\vskip.05in
\noindent The points of the tied spectrum consist of those $\mu\in {\bf R}$ such that there exists a 
nontrivial solution of $-f''+qf=\mu f$ with the Dirichlet boundary conditions  $f(0)=0=f(2\pi )$. These 
interlace the periodic spectrum, so we can assume that 
$\mu_j$ belongs to the $j^{th}$ interval of instability 
$(\lambda_{2j-1}, \lambda_{2j})$ for $j=1, 2, \dots $. For a typical 
$q\in (\Omega_N, \nu_N^\beta )$, as
$j\rightarrow\infty$ the length $\lambda_{2j}-\lambda_{2j-1}$ tends to zero, but the Lipschitz bounds on $q\mapsto
\mu_j$ becomes larger; so we need to balance these effects a specially 
formulated concentration theorem.\par  
\indent Let $\xi_j=2\sqrt{\mu_j}-j$, which is possibly 
complex for the first few $j$ and
assuredly real thereafter, and form $\xi=(\xi_j)_{j=1}^\infty$. Also let 
$$d_j=C_1N(N+1)/j, \quad \alpha_j=2^{-1}j^{-8}\exp \bigl(-C\beta^{5/2}N^{9/4}\bigr)\qquad 
(j=1, 2, \dots )\eqno(4.1)$$ 
\noindent for $C$ as in Lemma 2.2 and $C_1>0$ to be chosen.\par
\vskip.05in
\noindent {\bf Theorem 4.1} {\sl Let $F:\ell^2\rightarrow {\bf R}$ be $1$-Lipschitz.
Then 
$$\nu^\beta_N\Bigl\{ q\in \Omega_N: F(\xi (q))-\int_{\Omega_N} F(\xi )d\nu_N^\beta >s \Bigr\}
\leq e^{-c_F^*(s)}\qquad (s\geq 0),\eqno(4.2)$$
\noindent where $c_F^*:[0, \infty )\rightarrow [0, \infty ]$ is a convex 
function such that, for all $n=1, 2, \dots,$ }
$$c^*_F(s)\geq\cases{4^{-1}\alpha_n(s-d_n)^2,& for all $s\geq d_n$;\cr
4^{-1}(1+\alpha_nd_n^2)^{-1}\alpha_ns^2,& for all
$d_n+1/(\alpha_nd_n)\geq s\geq 0$.\cr}
\eqno(4.3)$$
\vskip.05in
\noindent Theorem 4.1 gives a Gaussian concentration inequality, except that the
constants change for the various ranges of $s$. In transportation theory,
 it is natural to have cost functions with a different shape for short distances, 
as in [33] page 593. The proof of Theorem 4.1 is split into three results in the
remainder of this section, and begins by considering finitely many $\mu_j(q)$, 
as functions of $q\in \Omega_N$. First we prove that each $\mu_j$ is a Lipschitz function of $q$.\par
\vskip.05in
\noindent {\bf Lemma 4.2} {\sl For all $N>0$ there exist $C(N), \kappa_2>0$ such
that $\Phi: q\mapsto (\mu_1(q), \mu_2(q), \dots ,\mu_n(q))$ is
Lipschitz from $(\Omega_N, L^2)$ to $(\ell_n^2, \Vert \,\cdot
\,\Vert_{\ell^2} )$ with constant $L(n,N)\leq \kappa_2C(N)n^4$.}\par
\vskip.05in
\noindent {\bf Proof.} 
Let $HS$ be the space of Hilbert--Schmidt operators on
$L^2[0, 2\pi ]$ with the usual norm. For $\zeta >N$, the operator $\zeta I-d^2/dx^2
+q$ with boundary conditions $f(0)=f(2\pi )=0$ is invertible and has inverse $G(q)$, and $G(q)$ is given by an
integral operator with bounded Greens function. 
 For these boundary conditions, we compute the Greens function for $f''=\zeta f$ and obtain
$$G_\zeta (x,y)=\cases{ {{\sinh x\sqrt{\zeta}\sinh (2\pi -y)\sqrt{\zeta}}\over{\sqrt{\zeta}\sinh 2\pi\sqrt{\zeta}}}&
 for $0\leq x\leq y\leq 2\pi$;\cr
{{\sinh (2\pi -x)\sqrt{\zeta }\sinh y\sqrt{\zeta}}\over{\sqrt{\zeta}\sinh 2\pi \sqrt{\zeta}}}& for 
$0\leq y\leq x\leq 2\pi$.\cr}\eqno(4.4)$$
\noindent Moreover, from the
identity
$G(q)-G(p)=G(q)(p-q)G(p)$ and bounds on the Greens function, we obtain a constant
$C(N)$ such that 
$$\bigl\Vert G(q)-G(p)\bigr\Vert_{HS}\leq C(N)\bigl\Vert
p-q\bigr\Vert_{L^2}.\eqno(4.5)$$
\indent Real Lipschitz functions operate on differences of
self-adjoint operators in the Hilbert--Schmidt norm, so for all Lipschitz functions
$\varphi :{\bf R}\rightarrow {\bf R}$ such that $\vert \varphi
(x)-\varphi (y)\vert \leq L\vert x-y\vert$, we have 
$$\bigl\Vert \varphi (G(q))-\varphi (G(p))\bigr\Vert_{HS}\leq
\kappa_2L\bigl\Vert
G(q)-G(p)\bigr\Vert_{HS}.\eqno(4.6)$$
\noindent By a result of Lidskii [28], 
the map $\Lambda :HS\rightarrow
\ell^2$ which associates to a
self-adjoint positive operator $A$ the decreasing list of eigenvalues
gives a Lipschitz function. In particular,
$$\bigl\Vert\Lambda (\varphi (G(q)))-\Lambda (\varphi
(G(p)))\bigr\Vert_{\ell^2}\leq 
\bigl\Vert \varphi (G(q))-\varphi (G(p))\bigr\Vert_{HS}.\eqno(4.7)$$
\indent The required Lipschitz function is partially specified by
$$ \varphi (y)=\cases{ 0, & if  $y\leq 1/2n^2$;\cr
(1/y)-\zeta, & if $1/n^2\leq y\leq 1/\zeta$;\cr
0, & if $y\geq 2/\zeta$,\cr}\eqno(4.8)$$
\noindent with straight line segments added to complete the graph and
make the function continuous. Note that $\varphi (1/(y+\zeta ))=y$ and $\vert\varphi'\vert \leq n^4$ in
the middle of the domain, so 
$$ \Lambda (\varphi (G(q))) =(\mu_n(q), \mu_{n-1}(q), \dots , \mu_1(q),
0, \dots ).\eqno(4.9)$$
\noindent By combining the Lipschitz maps in (4.5), (4.6) and (4.7), we obtain the stated result. {$\square$}\par
\vskip.05in
\noindent {\bf Definition} (Free energy) Let $(X, \mu )$ be a probability
space and $F$ a real random variable on $(X, \mu )$ such that $F$ has finite mean. Define the
normalized free energy of $F$ as $c_F:{\bf R}\rightarrow (-\infty ,\infty ]$, where
$$c_F(t)=\log\Bigl( \int_X e^{tF(\xi )}\mu (d\xi
)\Bigr)-t\int_X F(\xi )\mu (d\xi )\qquad (t\in{\bf R}).\eqno(4.10)$$
\noindent Then $c_F$ is convex and its Legendre--Fenchel transform is
defined by $$c_F^*(s)=\sup\{
st-c_F(t): t;
c_F(t)<\infty \};$$
\noindent see [32] page 23. Also, $c^*_F$ is known as the (normalized) rate function.\par
\vskip.05in
\indent The next step is to prove a concentration inequality for Lipschitz functions of
finitely many tied eigenvalues.\par
\vskip.05in
\noindent {\bf Lemma 4.3} {\sl Let 
$\alpha (N,\beta )=2^{-1}\exp(-C\beta^{5/2}N^{9/4})$,
$\xi_j=2\sqrt{\mu_j}-j$, and for a 
$1$-Lipschitz function $\Psi_n:({\bf C}^n, \ell^2) \rightarrow {\bf R}$, let 
$F_n(q)= \Psi_n(\xi_1(q), \dots
,\xi_n(q))$ on $(\Omega , \nu_N^\beta)$. Then $F_n$ satisfies the concentration
inequality}
$$c_{F_n}(t)\leq {{\kappa_2^2C(N)^2n^8t^2}\over{2\alpha (N, \beta )}}\qquad (t\in {\bf R}).
\eqno(4.11)$$
\vskip.05in
\noindent {\bf Proof.} By Lemma 4.2, $F_n
:\Omega_N\rightarrow {\bf R}$ is Lipschitz with constant $\kappa_2C(N)n^4$. By
Corollary 2 of [3], any Lipschitz function on $(\Omega_N, \nu_N^\beta )$
satisfies a Gaussian concentration of measure inequality as in
(4.11). Here $\alpha (N, \beta )$ is the logarithmic Sobolev constant of [3], Theorem 1. See [33, Theorem 22.10] for a detailed discussion. {$\square$}\par
\vskip.05in

\indent Let $(d_j)_{j=1}^\infty$ be a decreasing positive sequence with $d_j\rightarrow 0$ as
$j\rightarrow\infty$ and let 
$$X=\Bigl\{ (\xi_j)_{j=-\infty}^\infty\in \ell^2({\bf
C}): \vert \xi_j\vert \leq d_j/2; j=1, 2, \dots \Bigr\}.\eqno(4.12)$$
\vskip.05in
\noindent {\bf Proposition 4.4} {\sl Let $\mu$ be a Radon probability measure on
$X$, and suppose that there exists $\alpha_n>0$ such that $c_{F_n}(t)\leq t^2/(2\alpha_n)$ for all $1$-Lipschitz functions
$F_n:(X, \ell^2)\rightarrow {\bf R}$ that depend only on the first $n$ coordinates. 
Then for all $1$-Lipschitz functions
$F:(X, \ell^2)\rightarrow {\bf R}$, the normalized free energy satisfies
$$c_F(t)\leq {{t^2}\over{\alpha_n}}+{{1}\over{2}}\log\cosh (2d_nt)\qquad
(t\in {\bf R})\eqno(4.13)$$
\noindent and the Legendre transform satisfies (4.2)}
$$c^*_F(s)\geq\cases{4^{-1}\alpha_n(s-d_n)^2& for all $s\geq d_n$;\cr
4^{-1}(1+\alpha_nd_n^2)^{-1}\alpha_ns^2& for all $d_n+1/(\alpha_nd_n)\geq s\geq 0$.
\cr}\eqno(4.14)$$
\vskip.05in
\noindent {\bf Proof.} Let ${\cal F}_n$ be the $\sigma$-algebra of Borel sets
that is generated by the first $n$ coordinate functions on $X$, and for $F$
as above let $F_n=E(F\mid {\cal F}_n)$ be the conditional expectation with
respect to ${\cal F}_n$ in $L^2(\mu )$; note that $F_n$ is $1$-Lipschitz and $\int F_nd\mu
=\int Fd\mu $. Now let
$$Z_n(t)=\int_X e^{t(F-F_n)} d\mu\eqno(4.15)$$
\noindent which satisfies $Z_n(0)=1$, $Z_n'(0)=0$ and $Z_n''(t)=\int_X
(F-F_n)^2e^{t(F-F_n)} d\mu,$ so
$0\leq Z_n''(t)\leq d_n^2Z_n(t)$ and so $Z_n'(t)\geq 0$ for all $t>0$. Integrating this differential inequality,
we obtain $Z_n(t)\leq \cosh (d_nt)$. Then by the Cauchy--Schwarz inequality, we have
$$\int_X e^{t(F-\int Fd\mu)} d\mu\leq \Bigl(\int_X e^{2t(F-F_n)}
d\mu\Bigr)^{1/2}\Bigl(\int_X e^{t(F_n-\int F_nd\mu)} d\mu\Bigr)^{1/2},\eqno(4.16)$$
\noindent so that 
$$\eqalignno{c_F(t)&\leq (1/2)\log Z_n(2t)+(1/2)c_{F_n}(2t)\cr
&\leq (1/2)\log \cosh (2d_nt)+\alpha_n^{-1}t^2.&(4.17)\cr}$$
\indent To estimate $c_F^*(s)$, we first suppose that $0\leq s\leq d_n+1/\alpha_nd_n$, and use the estimate
$$\eqalignno{c_F^*(s)&\geq \sup\{ st-c(t): 0<t<1/2d_n\}\cr
&\geq \sup\Bigl\{ st-\alpha_n^{-1}t^2-d_n^2t^2:0<t<1/2d_n\Bigr\}\cr
&={{\alpha_ns^2}\over{4(1+\alpha_nd_n^2)}}.&(4.18)\cr}$$
\noindent For $s\geq d_n$, we have  
$$c^*_F(s)\geq \sup\bigl\{ st-\alpha_n^{-1}t^2-d_nt: t>0\bigr\}
=4^{-1}\alpha_n(s-d_n)^2.\eqno(4.19)$$        
\rightline{$\square$}\par
\vskip.05in   
\noindent {\bf Proof of Theorem 4.1} 
 As in (2.14), let $q\in \Omega_N$ be such that for some $C_1, C_2>0$ the periodic
spectrum satisfies $\lambda_{2j}-\lambda_{2j-1}\leq C_1N\vert \hat q(j)\vert$
and $C_2j^2\leq \lambda_{2j-1}$ for all $j=1, 2, \dots $. Also
$$\sum_{j=n}^\infty \vert\xi_j\vert^2=\sum_{j=n}^\infty \Bigl(
{{4\mu_j-j^2}\over{2\sqrt{\mu_j}+j}}\Bigr)^2 \leq
{{(CN(N+1))^2}\over{C_2n}}.\eqno(4.20)$$
\noindent Hence $\xi\in X$ for the above $X$. We can apply Proposition
4.4 to the measure $\mu$ that is induced on $X$ from $(\Omega_N,
\nu_N^\beta )$ by the map $q\mapsto\xi (q)$, since Lemma 4.3 implies that the
hypotheses of Proposition 4.4 hold with the constants as in (4.1). We can therefore introduce $c_F^*$ as the 
Legendre--Fenchel transform of $c_F$, where $c_F$ is continuously differentiable on
${\bf R}$, and hence $\{ s: c^*_F(s) <\infty\}$ contains the
range of $c'_F.$ Theorem 4.1 follows from Proposition 4.4 by Chebyshev's
inequality, since
$$\nu_{N}^\beta \Bigl\{ q\in \Omega_N: F(\xi (q))-\int F(\xi )d\nu_N^\beta >s\Bigr\}\leq
e^{-st+c_F(t)}\qquad (s,t>0),\eqno(4.21)$$
\noindent and we can optimize this inequality over $t>0$ for each fixed $n$, and use the bounds
(4.14). {${\square}$}\par 
\vskip.1in
\noindent {\bf 5. Transportation of measure and linear statistics}\par
\vskip.05in 
\indent Let $\tau =(\tau_j)_{j=-\infty}^\infty$ where
$$\tau_j =\cases{ \sqrt{2(\lambda_{2j}+\lambda_{2j-1})}& for $j=1,2,\dots $;\cr
0, & for $j=0$;\cr
-\sqrt{ 2(\lambda_{-2j}+\lambda_{-2j-1})}& for $j=-1,-2, \dots $,\cr}\eqno(5.1)$$
\noindent and for $g\in RPW(\pi )$ introduce the linear statistic
$$F(\tau )=\sum_{j=-\infty}^\infty (g(\tau_j)-g(j)).\eqno(5.2)$$
In this section we obtain bounds on $F(\tau )$ and its fluctuations. 
By Propositions 2.3 and 3.1, $\tau_j$ is close to $j$, and the new idea is that the fluctuations of
$\sum_{j=-m}^m (g(\tau_j)-g(j))$ are controlled, or at least attenuated, by
 $\sum_{j=-m}^m \vert g'(\tau_j)\vert^2$. To make this precise, we 
introduce the infimum convolution, which generalizes the
Legendre--Fenchel transform, and then we consider the special case of linear
statistics.\par 
\vskip.05in
\noindent {\bf Definition} (Infimum convolution) For a continuous and bounded $F:\ell^2\rightarrow {\bf R}$, the Hopf--Lax infimum
convolution is defined as in [5] by 
$$Q_sF(\eta )=\inf\Bigl\{ F(\xi )+{{1}\over{2s}}\Vert \xi-\eta\Vert^2: \xi \in 
\ell^2\Bigr\} \qquad (\eta \in
\ell^2, s>0).\eqno(5.3)$$
\vskip.05in
\indent The purpose of the time parameter $s$ in $Q_s$ is to produce a semigroup that solves
the Hamilton--Jacobi equation
$${{\partial}\over{\partial s}}Q_sF=-(1/2)\Vert\nabla Q_sF\Vert^2.\eqno(5.4)$$
\vskip.05in
\noindent {\bf Proposition 5.1} {\sl  For $g\in RPW(\pi )$ with
$\Vert g\Vert_{L^2}\leq 1$ and for $\xi =(\xi_j)_{j=-\infty}^\infty \in \ell^2$ let
$$F(\xi )=\sum_{j=-\infty}^\infty (g(\xi_j +j)-g(j)).\eqno(5.5)$$
\noindent (i) Then $F$ is bounded, $F\mapsto Q_sF(\eta )$ is concave, and 
$Q_sF(\eta )\leq F(\eta )$ for all $s>0$ and $\eta \in \ell^2$;\par
\noindent (ii)  the semigroup law $Q_{s+t}F=Q_tQ_sF$ holds for all $s,t>0$;\par
\noindent (iii) $Q_sF(\eta )\rightarrow F(\eta )$ as
$s\rightarrow 0+$; and}
$$\Bigl({{\partial }\over{\partial
s}}\Bigr)_{s=0}Q_sF(\eta )={{-1}\over{2}}\sum_{j=-\infty}^\infty g'(\eta_j +j)^2\qquad (\eta
=(\eta_j)_{j=-\infty}^\infty \in\ell^2).\leqno(iv)$$   
\vskip.05in
\noindent {\bf Proof.} (i) First we check that the series defining $F (\xi
)$ is convergent and $F$ is bounded. By Fourier
analysis, $g'\in RPW(\pi )$. Also, by Proposition 3.4, for all $\xi= (\xi_j)_{j=-\infty}^\infty \in \ell^2$, 
we can perturb $(j)_{j=-\infty}^\infty$ to obtain a sampling sequence 
$(\xi_j+j)_{j=-\infty}^\infty$ 
for $RPW(\pi )$ by [25]. Hence there exists $B>0$ such that
$$\sum_{j=-\infty}^\infty \vert g'(j+\xi_j)\vert^2\leq 
B\int_{-\infty}^\infty \vert g'(x)\vert^2 dx\eqno(5.6)$$
\noindent for all $g\in RPW(\pi )$ and all $\xi\in \ell^2$ such that 
$\Vert \xi\Vert_{\ell^2}\leq 1$. From the Cauchy--Schwarz inequality, we deduce that $\sum_{j=1}^\infty g'(j+\theta
\xi_j)\xi_j$ converges absolutely for all $0<\theta <1$, and hence
$\sum_{j=-\infty}^\infty g(\xi_j+j)-g(j))$ is absolutely convergent by the
mean value theorem. The remaining statements are now straightforward.\par
\indent (ii) This is a general property of the quadratic Hamilton--Jacobi
semigroup [33, p. 584].\par
\indent (iii) and (iv) We prove that for $0<s<\sqrt{5}/\pi^2$
there exists for each $\eta\in \ell^2$ a unique $\xi\in\ell^2$ in a
neighbourhood of $\eta$ such that 
$Q_sF(\eta )=F((\xi_j+j)_{j=-\infty}^\infty )+\Vert
(\eta_j-\xi_j)_{j=-\infty}^\infty\Vert^2/(2s)$. Using the Fourier representation, one can easily show that $g$ has bounded derivatives of order
$k$ such that $\Vert g^{(k)}\Vert_\infty \leq \pi^k\Vert g\Vert_{L^2}
/\sqrt{2k+1}$ for $k=0, 1, \dots.$  So for  all $\eta_j\in {\bf R}$,
$\sqrt{5}/\pi^2>s>0$ and $j\in {\bf Z}$, the real function $\xi_j\mapsto
g(j+\xi_j)+(\xi_j-\eta_j)^2/(2s)$ is differentiable and diverges to
infinity as $\xi_j\rightarrow\pm \infty$, hence attains its infimum at
a unique $\xi_j$ such that $sg'(\xi_j+j)+\xi_j=\eta_j$ and 
$\vert \eta_j-\xi_j\vert \leq s\pi \Vert
g\Vert_{L^2}/\sqrt{3}$. By the mean value theorem, 
we deduce that there exists $\zeta_j$ 
between $\xi_j$ and $\eta_j$
such that
$$g(\xi_j+j)-g(j)+{{1}\over{2s}}(\xi_j-\eta_j)^2=g(\eta_j+j)-g(j)-
{{s}\over{2}}
g'(\xi_j+j)^2-{{1}\over{2}}g''(\zeta_j+j)(\xi_j-\eta_j)^2.\eqno(5.7)$$
\indent The map $T:\ell^2\rightarrow\ell^2$ given by $T:(\xi_j)_{j=-\infty}^\infty \mapsto
(\xi_j+sg'(\xi_j+j))_{j=-\infty}^\infty$ is continuous and has Fr\'echet derivative
$\nabla T(\xi )={\hbox{diagonal}}\,
(1+sg''(\xi_j+j))_{j=-\infty}^\infty$; so for $s\pi^2\Vert
g\Vert_{L^2}/\sqrt{5}<1$, the operator $\nabla T$ is invertible on $\ell^2$, and hence $T$ is an open
mapping. By summing (5.5) over $j$, we deduce that
$$Q_sF(\eta )=F(\eta )-{{s}\over{2}}\sum_{j=-\infty}^\infty g'(\xi_j +j)^2
+O(s^2),\eqno(5.8)$$
\noindent where the series converges 
to $\sum_{j=-\infty}^\infty g'(j+\eta_j)^2$ as $s\rightarrow 0+$. 
{$\square$}\par
\vskip.05in
\noindent {\bf Theorem 5.2} {\sl There exists $\rho =\rho (N, \beta, m)>0$ such that for
all $g\in RPW(\pi )$ with $\Vert g\Vert_{L^2}\leq 1$, the linear statistic 
$F_m(\tau )=\sum_{j=-m}^m (g(\tau_j)-g(j))$ satisfies}
$$\int_{\Omega_N} \exp\bigl({\rho Q_1F_m(\tau (q) )}\bigr)\nu_N^\beta (dq) \leq 
 \exp\Bigl({\rho \int_{\Omega_N} F_m(\tau (q) )
\nu_N^\beta(dq)}\Bigr).\eqno(5.9)$$
\vskip.05in
\noindent {\bf Proof.} By Proposition 2.3(ii), we can assume that all the
periodic eigenvalues are simple. Then we note that $\lambda_j>0$ for all but finitely 
many indices $j$; whereas, in the exceptional cases where
$\lambda_{2j}<0$, we have $\tau_j$ purely imaginary so the sum
 $g(\tau_j)+g(\tau_{-j})$ unambiguously gives a real random variable for all $g\in RPW(\pi )$.\par
\indent We show that the solution operator of Hill's
equation and hence the characteristic function are Lipschitz functions of 
$q\in \Omega_n$, so we can apply known concentration theorems to control some linear
statistics. Let $p$ be an integer such that $m^2<p<(m+1)^2$, and let $S$ be the oriented square with vertices $\pm p\pm ip$ which is described 
 once in the positive sense. Taking $m$ to be large, we can use the nonrandom bound 
$\Vert q\Vert^2_{L^2}\leq N$ to bound the terms in the series (3.7) uniformly for all $q\in \Omega_N$. We write Hill's equation in the style
$$\eqalignno{\Bigl( \left[\matrix{ 0&1\cr -1 &0\cr}\right] 
{{d}\over{dx}}& -{{1}\over{2}}\left[\matrix{ 1+q& -i(1-q)\cr -i(1-q)& -(1+q)\cr}\right]
-{{\lambda}\over{2}}\left[\matrix{1&0\cr0&1\cr}\right]\Bigr)\Psi_\lambda (x)=0.&(5.10)\cr}$$ 
\noindent Let 
$$R={{1}\over{2}}\left[\matrix {1+q& iq-i\cr iq-i&-1-q\cr}\right].\eqno(5.11)$$
\noindent The integral equation
$$\Psi_\lambda (x)=\exp \Bigl( x\sqrt{\lambda}\left[\matrix{0&-1\cr 1&0\cr}\right]\Bigr)
\Psi_\lambda (0) +\int_0^x \exp\Bigl((x-s)\sqrt{\lambda}\left[\matrix{ 0&-1\cr 1&
 0\cr}\right]\Bigr) R(s)\Psi_\lambda (s)\,ds\eqno(5.12)$$
\noindent has variational equation
$$\eqalignno{\delta\Psi (x)&=\int_0^x \exp\Bigl((x-s)\sqrt{\lambda}\left[\matrix{0&-1\cr 1&0\cr}\right]\Bigr)
 \delta R(s)\Psi_\lambda (s)ds\cr
&\quad +\int_0^x \exp\Bigl((x-s)\sqrt{\lambda}\left[\matrix{ 0&-1\cr 1&0\cr}\right]\Bigr) R(s)\delta\Psi (s)\, ds.
&(5.13)\cr}$$
\noindent  Then $\Vert \Psi_\lambda\Vert\leq M$. We choose $\kappa>0$ such that 
$\sqrt{\kappa N}e^{p\kappa } <1/2$ and split $[0, 2\pi ]$ into consecutive subintervals of length $\kappa$; then we consider the supremum norm $\Vert\,.\,\Vert_\infty$ on the matrix functions to obtain the bound
$$\Vert\delta \Psi\Vert_\infty \leq \sqrt{\kappa}\Bigl( \int_0^\kappa \Vert R(s)\Vert^2ds\Bigr)^{1/2} e^{p\kappa} 
\Vert \delta \Psi\Vert_\infty +e^{p\kappa} \sqrt {\eta} \Big(\int_0^\kappa\Vert \delta R(s)\Vert^2ds\Bigr)^{1/2} 
\Vert \Psi_\lambda\Vert_\infty\eqno(5.14)$$
\noindent so that 
$$\Vert \delta \Psi\Vert_\infty \leq {{2}\over{\sqrt{N}}}\Bigl( \int_0^\kappa \Vert\delta R(s)\Vert^2ds\Bigr)^{1/2}
 \Vert \Psi\Vert_\infty .\eqno(5.15)$$
\noindent We repeat this bound for each successive interval and thus we obtain a Lipschitz constant
$$L_\kappa\leq    \Bigl( 1+{{M}\over {2\sqrt{N}}}\Bigr)^{2\pi/\kappa }.\eqno(5.16)$$
\noindent We deduce that $q\mapsto\Delta (\lambda )$ is Lipschitz continuous from
$(\Omega_N, L^2)$ to the space of holomorphic functions inside $S$ with
 the uniform norm. The function
 $\lambda \mapsto \Psi_\lambda (s)$ is holomorphic and hence we can differentiate the
 integral equations (5.9) with respect to $\lambda$. By Cauchy's estimates, $q \mapsto \Delta'(\lambda )$ is Lipschitz continuous.\par  
We observe that there exists $\delta >0$ such that $\vert \Delta (\lambda )^2-4\vert >\delta $ for all $\lambda$ on $S$.
 Note that the zeros of $\Delta (\lambda )^2-4$ all lie on the real axis, and their position are described in the proof of Proposition 3.1. So we can form the product
$$\Delta (\lambda )^2-4=c\prod_{j=0}^\infty 
\Bigl( 1-{{\lambda}\over{\lambda_j}}\Bigr)\eqno(5.17)$$
\noindent By estimating each factor on the compact set $S$, we obtain a lower bound $\delta >0$.\par
\indent It follows from these estimates that the functions 
$$q\mapsto {{1}\over {\pi i}}\int_S {{\Delta (\lambda )\Delta'(\lambda )\lambda^kd\lambda }\over {\Delta (\lambda )^2 
-4}}\qquad (k=1, \dots ,m)\eqno(5.18)$$
\noindent are Lipschitz continuous on $\Omega_N$. \par
\indent As in [4], we introduce the circles $C(n^2/4, 1/4)$ and apply Cauchy's integral formula to obtain
$${{\lambda_{2j}+\lambda_{2j-1}}\over{2}}={{1}\over{2\pi
i}}\int_{C(n^2/4, 1)} {{\lambda \Delta (\lambda
)\Delta'(\lambda ) \, d\lambda}\over{\Delta (\lambda )^2-4}}\eqno(5.19)$$
\noindent for the midpoint of the $j^{th}$ interval of instability. Hence $q\mapsto
(\lambda_{2j}+\lambda_{2j-1})/2$ is Lipschitz from $(\Omega_N, L^2)\rightarrow{\bf R}$ with constant
$L_j\leq \kappa_0(N, \beta ) j^2$ for some $\kappa_0(N, \beta )$ and all $j$. We deduce that $q\mapsto
\tau_j$ is also Lipschitz. (It follows likewise that the functions $q\mapsto \sum_{j=1}^{2m}\lambda_j^k$ are Lipschitz
 continuous for $k=1, \dots $, but we have not quite proved that $q\mapsto \lambda_j^k$ is Lipschitz.)\par
\indent Let $\mu_m$ be the probability measure that is induced from $\nu_N^\beta$ by the Lipschitz map 
$\varphi :q\mapsto (\tau_j)_{j=-m}^m$, where the Lipschitz constant $L_m$ is finite. The measure $\nu_N^\beta$ satisfies a logarithmic Sobolev inequality 
$$LSI(1/\kappa ):\qquad \int_{\Omega_N} f(q )^2 \log \Bigl( f(q) )^2/\int f^2d\nu_N^\beta \Bigr) \nu_N^\beta (dq )
\leq 2\kappa \int_{\Omega_N}\Vert \nabla f(q)\Vert^2\nu_N^\beta (dq )\eqno(5.20)$$
by Theorem 1 of [3]. Now $\varphi$ induces $\mu_m$ on ${\bf R}^{2m+1}$ from 
$\nu_N^\beta$ on $\Omega_N$. Also, by the chain rule applied to $f(q)=F_m(\varphi (q))$ we have  
$\Vert \nabla f(q)\Vert\leq L_m\Vert (\nabla F_m)\circ \varphi \Vert$, and hence we have the logarithmic Sobolev inequality     
such that 
$$\int_{{\bf R}^{2m+1}} F_m(\xi )^2 \log \Bigl( F_m(\xi )^2/\int F_m^2d\mu_m\Bigr) \mu_m(d\xi )
\leq 2L^2_m\kappa\int_{{\bf R}^{2m+1}}\Vert \nabla F_m(\xi )\Vert^2\mu_m(d\xi )\eqno(5.21)$$
\noindent for all $C^\infty$ functions $F_m:{\bf R}^{2m+1}\rightarrow {\bf R}$.
 In particular, this applies to the linear statistic $F_m$ of (5.7), since $\xi\mapsto Q_tF_m(\xi )$ is
Lipschitz continuous and $t\mapsto Q_tF_m(\xi )$ is differentiable for almost all $\xi$, by [5, p 673]. 
The logarithmic Sobolev inequality (5.21) implies 
the transportation inequality (5.9) for $\mu_m$ with constant $\rho>0$ depending upon 
$L^2_m\kappa$ by [32, p. 292]. The crucial point is that
$$\phi (t)={{L^2_m\kappa }\over{t}}\log\Bigl(\int_{{\bf
R}^{2m+1}}e^{tQ_t(F_m)(\xi )/\kappa L^2_m}\mu_m(d\xi )\Bigr)$$
\noindent is a decreasing function of $0<t<1$ by the Hamilton--Jacobi equation
and (5.21), so $\phi(1)\leq \lim_{t\rightarrow 0+}\phi (t)$. {$\square$}\par
\vskip.05in
\noindent {\bf Corollary 5.3} {\sl There exists $\kappa (m, \beta, N) >0$ such that for all 
$g\in RPW(\pi )$ with $\Vert g\Vert_{L^2}\leq 1$, the linear statistic 
$F_m(\tau )=\sum_{j=-m}^m (g(\tau_j)-g(j))$ satisfies}  
$$\nu_N^\beta \Bigl\{ q\in\Omega_N:  F_m(\tau (q) )-\int
F_md\nu_N^\beta >\varepsilon \Bigr\}\leq 
e^{-\kappa (m, \beta ,N)\varepsilon^2}\qquad
(\varepsilon >0),\eqno(5.22)$$   
\vskip.05in
\noindent {\bf Proof.} This follows from Theorem 5.2 by 
[33, Theorem 22.22]. {$\square$}\par
\vskip.05in

\noindent {\bf Proposition 5.4} {\sl For $1>\beta >0$ and $N_1>N>0$, there exist 
$\kappa_1(N,\beta )>0$ and $\kappa_2(N, \beta )>0$ such that, for all $g\in RPW(\pi )$ with $\Vert
g\Vert_{L^2}\leq 1$, the corresponding linear statistic satisfies}
$$\nu_N^\beta \Bigl\{ q\in \Omega_N: F(\tau (q))-\int Fd\nu_N^\beta >\varepsilon \Bigr\}\leq 
\exp\bigl(-\kappa_1(n,
\beta )\varepsilon^2\min \{1, \kappa_2 (N, \beta )\varepsilon^4\}\bigr)\qquad (\varepsilon >0 ).
\eqno(5.23)$$
\vskip.05in   
\noindent {\bf Proof.} Let $F_m(\tau )$ be the partial sum $F_m(\tau
)=\sum_{j=-m}^m (g(\tau_j)-g(j))$.   
Now we use the mean value theorem and the Cauchy--Schwarz inequality to bound
$$\eqalignno{ \vert F(\tau )-F_m(\tau )\vert&=\Bigl\vert\sum_{j: \vert j\vert \geq m+1}
\bigl(g(\tau_j)-g(j)\bigr)\Bigr\vert\cr
&\leq \Bigl( \sum_{j: \vert j\vert \geq m+1}\vert \tau_j-j\vert^2\Bigr)^{1/2}\Bigl( 
\sum_{j: \vert j\vert \geq m+1}\vert g'(\eta_j )\vert^2\Bigr)^{1/2}\cr
&\leq \Bigl(\sum_{j: \vert j\vert \geq
m+1}\bigl\vert{{\tau_j^2-j^2}\over{j+\tau_j}}\bigr\vert^2\Bigr)^{1/2}\Bigl(\sum_{j=-\infty}^\infty
\vert g'(\eta_j)\vert^2\Bigr)^{1/2},&(5.24)\cr}$$
\noindent for some $\eta_j$ between $j$ and $\tau_j$. Then we use (5.4) and (3.1) to obtain the bound
$$\vert F(\tau )-F_m(\tau )\vert\leq m^{-1}B\bigl( C(\beta )N(N+1)\bigr)^{1/2}\Vert
g'\Vert_{L^2}\eqno(5.25)$$
\noindent for some constant $C(\beta )$, and $\Vert g'\Vert_{L^2}\leq
\pi \Vert g\Vert_{L^2}$. \par
\indent Hence, in the decomposition 
$$\Bigl[ F-\int F>\varepsilon \Bigr]\subseteq  
\Bigl[ F_m-\int F_m>\varepsilon/3 \Bigr]\cup \Bigl[ \int F_m-\int F>\varepsilon/3 \Bigr]\cup
\Bigl[ F-\int F_m>\varepsilon /3 \Bigr],\eqno(5.26)$$
\noindent we can ensure that the final two events are empty by selecting $m\geq \max\{ 1,
\kappa/\varepsilon \}$.\par
\indent For such an $m$, the function $q\mapsto F_m(\tau )$ is Lipschitz with constant of order $m^2$,
so by Corollary 2 of [3], there exists $\alpha (N, \beta )$ such that
$$\nu_N^\beta\Bigl\{ q\in \Omega_N:F_m -\int F_md\nu_N^\beta >\eta \Bigr\} 
\leq \exp\bigl( -\alpha (N, \beta )\eta^2/m^4\bigr)\qquad
(\eta >0).\eqno(5.27)$$
\noindent We choose $\eta =\varepsilon/3$ to conclude the proof. $\square$\par
\vskip.05in
\noindent {\bf 6. Riesz bases}\par
\vskip.05in
\noindent The purpose of this section is to prove the following theorem 
concerning sampling from
Paley--Wiener space on a sequence given by the random eigenvalues, as in
(1.9).\par
\vskip.05in
\noindent {\bf Definition} (Gram matrix) For a real
sequence $t=(t_j)_{j=-\infty}^\infty$, the Gram matrix of the sequence
$({\hbox{sinc}}(x-t_j))_{j=-\infty}^\infty$ in $RPW(\pi )$ is 
$$\Gamma(t)=\bigl[{\hbox{sinc}}(t_j-t_k)\bigr]_{j,k=-\infty}^\infty.\eqno(6.1)$$
\vskip.05in
\noindent {\bf Theorem 6.1} {\sl Let $t=(t_j)_{j=-\infty}^\infty$ be the modified eigenvalues of Hill's equation as in 
(1.9). There exist $\beta>0$ and $N>0$ such that on a set of 
strictly positive $\nu_N^\beta$ measure,} \par
\noindent {\sl (i)} $({\hbox{sinc}}(t-t_j))_{j=-\infty}^\infty$ {\sl gives a Riesz basis 
for $RPW(\pi );$\par
\noindent (ii) there exists a corresponding system of biorthogonal functions
$(g_j)_{j=-\infty}^\infty$ such that for all $f\in RPW(\pi )$, the series}
$$f(x)=\sum_{j=-\infty}^\infty \langle f,g_j\rangle
{\hbox{sinc}}(x-t_j)\eqno(6.2)$$
\noindent {\sl converges in $L^2$, and $(\langle f,g_j\rangle_{L^2}
)_{j=-\infty}^\infty \in \ell^2$;}\par
\noindent {\sl (iii) the Gram matrix $\Gamma(t)$ is invertible, and belongs to
$I+HS$;}\par
\noindent {\sl (iv) the Carleman determinant satisfies}
$$\det_2 \Gamma(t)=\lim_{n\rightarrow\infty} \bigl\vert\det \bigl[
{\hbox{sinc}}(t_j-k)\bigr]_{j,k=-n}^n\bigr\vert^2>0.\eqno(6.3)$$

 \vskip.05in
\indent  The main idea of this section is to analyse $\Gamma(t)$ in the space $I+HS$, about which we record some facts. First, $I+HS$ is a complete metric space for the $HS$ norm. By elementary functional calculus, one can carry out
a type of polar decomposition in $I+HS$. The group ${\cal G}=\{ X\in I+HS: \det_2X\neq 0\}$
contains a subgroup ${\cal K}=\{ W\in {\cal G}: W^\dagger W=WW^\dagger =I\}$ and a convex
set ${\cal P}_+=\{ G\in {\cal G}: G=G^\dagger, G\geq 0\}$. 
The unitary group acts on each of ${\cal G}, {\cal K}$ and ${\cal P}_+$ 
by $(W,X)\mapsto WXW^\dagger$. In particular, for each $G\in {\cal P}_+$ 
there exists a unitary $W$ and a diagonal operator
$D={\hbox{diag}}(x_j)\in {\cal P}_+$ with eigenvalues $x_j>0$ and $\sum_j (x_j-1)^2$
finite such that $G=WDW^\dagger$. Also ${\cal K}\times {\cal K}$
acts on ${\cal G}$ by $(U,V): X\mapsto UXV^\dagger$ for $U,V\in {\cal
K}$ and $X\in {\cal G}$, and the map ${\cal K}\times 
{\cal K}\times {\cal P}_+\rightarrow {\cal G}:$ $(U,V,G)\mapsto UGV^\dagger $ is surjective. The space ${\cal P}_+$ is introduced to describe Gram
matrices.\par
\vskip.05in
\noindent {\bf Lemma 6.2} {\sl The Gram matrix $\Gamma(t)$ belongs to $I+HS$, 
and $\det_2 \Gamma(t)$ is a Lipschitz function of $t$ on bounded 
subsets of $Z+\ell^2$.}\par
\vskip.05in
\noindent {\bf Proof.} In the notation of Lemma 3.3, $U_t-I$ and hence $U_t^\dagger U_t-I$ are Hilbert--Schmidt. Taking the 
inverse Fourier transform ${\cal F}:L^2[-\pi , \pi ]\rightarrow RPW(\pi )$,
we deduce that 
$$\Gamma(t)=\bigl[ \langle   
 U_t(e^{ijx}),U_t(e^{ikx})\rangle\bigr]_{j,k=-\infty}^\infty\eqno(6.4)$$
\noindent belongs to $\{ G\in I+HS: G=G^\dagger, G\geq 0\}$ for all
$t\in Z+\ell^2$, and has a Carleman determinant.\par
\indent By Lemma 3.3, $(t_j)\mapsto U_t$ is Lipschitz on $Z+\ell^2$, with $\Vert
U_t-U_s\Vert_{HS}\leq 2\pi \Vert t-s\Vert_{\ell^2}/\sqrt{3}$, so  
$(t_j)\mapsto U_t^\dagger U_t$ is Lipschitz on bounded subsets of $Z+\ell^2$. 
On the convex bounded set $\{ A: \Vert A\Vert_{HS}\leq M\}$, the function $A\mapsto \det_2(I+A)$ is Lipschitz 
with constant $L\leq e^{c(2M+1)^2}$, where $c>0$ is some universal constant. Hence 
$t\mapsto \det_2 U_t^\dagger U_t$ is
Lipschitz continuous on bounded sets. {$\square$}\par
\vskip.05in
\noindent {\bf Proof of Theorem 6.1}  (i) and (iii). Then $\Gamma(t)$ represents $S_tS_t^\dagger$ with respect to the 
standard orthonormal basis of $\ell^2$. To prove that $S_t^\dagger$ is an embedding, it suffices to show that
 $\det_2 \Gamma(t) >0$ on a set of positive measure. 
To do this, we choose $0<\delta <1/2$ and then $n>(\beta N)^{1/\delta }\kappa^{1/\delta }$ as in 
Proposition 3.1, and reduce the
 analysis to a finite rank operator $A_n$. \par
\indent By Proposition 3.1, for all $N$, there exists $M$ such that $\Vert (t_j-j)\Vert_{\ell^2}\leq M$ for
all $(t_j(q))$ that arise as sampling sequences for $q\in \Omega_N$. By Lemma 6.2, 
there exists
$L>0$ such that $t\mapsto \det_2 \Gamma(t)$ is $L$-Lipschitz on $\{ t=(t_j):\Vert
(t_j-j)\Vert_{\ell^2}\leq M\}$.\par
\indent  We introduce the finite-rank operator 
$A_n:L^2[-\pi , \pi ]\rightarrow L^2[-\pi
,\pi ]$ such that
$$A_n(e^{ikx})=\cases{ e^{it_kx}-e^{i kx},& for $k=-n, \dots ,n$;\cr
                              0& else\cr}\eqno(6.5)$$
\noindent and write $U_t$ from Lemma 6.1 as $U_t=I+A_n+B_n$. Observe that 
$q\mapsto A_n$ gives a
Lipschitz map $\Omega_N\rightarrow HS$, such that $0\mapsto 0$. Note also that  
$[{\hbox{sinc}}(t_j-t_k)]_{j,k=-n}^n$ is a
block submatrix of $(I+A_n^\dagger )(I+A_n)$, and that 
$\langle (I+A_n^\dagger )(I+A_n)(e^{ikx}), e^{ikx}\rangle =1$ for all $k\in {\bf
Z}$, so 
$$\eqalignno{\det_2(I+A_n^\dagger )(I+A_n)&=\bigl(\det (I+A_n^\dagger
)(I+A_n)\bigr)e^{-{\hbox{trace}}(A_n+A_n^\dagger +A_n^\dagger A_n)}\cr
&=\bigl\vert\det (I+A_n)\bigr\vert^2\cr
&=\Bigl\vert\det
\bigl[{\hbox{sinc}}(t_j-k)\bigr]_{j,k=-n}^n\Bigr\vert^2,&(6.6)\cr}$$ 
\noindent a formula reminiscent of kernels from random matrix
theory [14, p 124]. Then by Lemma 6.2, 
$$\eqalignno{\int_{\Omega_N} \det_2 \Gamma(t)\,\nu_N^\beta (dq)&\geq \int_{\Omega_N} \det_2\bigl(
(I+A^\dagger_n)(I+A_n)\bigr) \nu_N^\beta (dq)\cr
&\quad -L\int_{\Omega_N}\bigl(2(1+M)\Vert
B_n\Vert_{HS} +\Vert B_n\Vert^2_{HS}\bigr)\nu_N^\beta (dq)&(6.7)}$$
\noindent where, as in (6.5),
$$\Vert B_n\Vert_{HS}\leq {{2\pi}\over{\sqrt{3}}}\Bigl( \sum_{k: \vert k\vert \geq n}
(t_k-k)^2\Bigr)^{1/2}\eqno(6.8)$$
\noindent and 
$$\Vert A_n\Vert_{HS}\leq {{2\pi}\over{\sqrt{3}}}\Bigl( \sum_{k=-n}^n
(t_k-k)^2\Bigr)^{1/2}.\eqno(6.9)$$
\noindent By Proposition 3.1, we have $\sum_{k=n}^\infty
(t_k-k)^2\leq C\beta Nn^{2\delta -1}/(1-2\delta )$ so we can choose $n$ so large that the final 
integral in (7.14) is less than $1/8$. Then we choose 
$N, \beta >0$ so small that 
$$\eqalignno{\int_{\Omega_N} \det_2\bigl(
(I+A^\dagger_n)(I+A_n)\bigr) \nu_N^\beta (dq)
&\geq 
1-{{2L\pi}\over{\sqrt{3}}}\int_{\Omega_N}\Bigl( \sum_{k=-n}^n
(t_k-k)^2\Bigr)^{1/2}\nu_N^\beta (dq)\cr
&>3/4;&(6.10)\cr}$$
\noindent hence (6.7) is greater than $1/8$ some some $N, \beta >0$. Consequently $\det_2\Gamma(t)>1/8$ on an open set that 
has strictly positive $\nu_N^\beta$ measure and there
$({\hbox{sinc}}(s-t_j))_{j=-\infty}^\infty$ gives a Riesz basis for its linear span. Since $\Vert
B_n\Vert\rightarrow 0$ as $n\rightarrow\infty$, we can arrange for $U_n$ to be invertible, so the
linear span is all of $RPW(\pi )$.\par
\indent (iv) As $n\rightarrow\infty$, we have $\Vert
B_n\Vert_{HS}\rightarrow 0$, and hence $\vert\det
(I+A_n)\vert^2\rightarrow \det_2\Gamma(t)$.\par
\indent (ii) There exists an 
invertible linear operator
$W_t:L^2[-\pi ,\pi ]\rightarrow RPW(\pi )$ such that $W_t(e^{ijx}) 
={\hbox{sinc}}(s-t_j)$ and then $(g_j)_{j=-\infty}^\infty$ where
$g_j=(W_t^\dagger )^{-1}(e^{ijx})$ gives a system of biorthogonal functions for 
$({\hbox{sinc}}(s-t_j))_{j=-\infty}^\infty$. Note that $(g_j)_{j=-\infty}^\infty$ is itself a Riesz
basis, with $\Vert g_j\Vert_{L^2}\leq M_1$ for some $M_1$ and all $j\in {\bf Z}$, 
and 
$[{\hbox{sinc}}(t_j-t_k)]$ has inverse matrix $[\langle g_j, g_k\rangle ].$\par
\rightline{$\square$}\par
\vskip.05in
\noindent {\bf Remarks 6.3} (i) In [26] page 124 , Katz and Sarnak considered a version of the Gram matrix $\Gamma(t)$ when the entries arise from
eigenvalues of random unitary operators. The determinant of (6.3) resembles some
expressions which appear in the representation theory of classical
groups as in [14] page 122, especially when written in the style
$$\det \bigl[ {\hbox{sinc}}
(t_j-k)\bigr]_{j,k=-n}^n$$
$$={{1}\over{(2n+1)!}}\int_{[-\pi, \pi ]^{2n+1}}  
{{\det [ e^{it_jx_k}]_{j,k=-n}^n}\over  
{\det [ e^{ijx_k}]_{j,k=-n}^n}}\prod_{-n\leq j<k\leq n} \vert
e^{ix_j}-e^{ix_k}\vert^2 {{dx_{-n}}\over{2\pi}}\dots 
{{dx_{n}}\over{2\pi}}.\eqno(6.11)$$
\noindent This integral formula follows from Andr\'eief's identity and the usual
Vandermonde determinant.\par
\indent (ii) In terms of Theorem E of [25], the phase function for the sampling sequence
$(t_j)_{j=-\infty}^\infty$ is, for some real $\alpha$,  
$$\varphi (x)=\alpha +{\hbox{sign}}(x)\int_{\{ \lambda \in [0, x^2]:
4-\Delta (\lambda )^2\geq 0\}} {{\vert \Delta'(\lambda )\vert \,
d\lambda}\over{\sqrt{4-\Delta (\lambda )^2}}},\eqno(6.12)$$
\noindent so that $\varphi$ is continuous, 
increasing, constant when $x^2$ belongs to
an interval of instability, and increases by $\pi$ as $x^2$ increases 
over each interval of stability. These properties follow from the fact 
that $\Delta' (\lambda
)$ is of constant sign on each interval of stability by Laguerre's
theorem.\par 
\vskip.05in
\noindent {\bf 7. The Jacobian and linear statistics}\par
\vskip.05in
\indent In this section we consider the set ${\cal M}_\lambda$ of potentials $q\in \Omega_N$ that have a given
periodic spectrum $\lambda =(\lambda_{j})_{j=0}^\infty$. Hochstadt [21, p. 219] observed that if only finitely many of the zeros of $\Delta (x)^2-4$
are simple, so
that the spectrum of $q$ has only finitely many gaps, then
$\Delta'(x)/\sqrt{\Delta (x)^2-4}$ is an algebraic function. However, for typical $q$ in
$(\Omega_N,\nu_N^\beta )$, Proposition 2.3(ii) shows that all the periodic eigenvalues are simple, so Hill's curve ${\cal E}$ of (1.11) 
is a hyperelliptic transcendental Riemann surface that has infinite genus. In Lemma 7.1 we define a 
suitable space of divisors on ${\cal E}$, and a map which associates to each $q\in {\cal M}_\lambda$ a
divisor $\delta$ on ${\cal E}$ which is determined by the tied spectrum $(\mu_j)_{j=1}^\infty$ of
$q$. Moreover, there is a pairing of divisors with the 
differential $\omega_\infty=(\Delta'(x)/\sqrt{\Delta (x)^2-4})dx$. The main task is to interpret (1.12) and the addition rule
on divisors. Clearly, translating $q(x)$ to $q(x+s)$
preserves the periodic spectrum, and the measure $\nu_N^\beta$, but changes the
tied spectrum. Suppose that $\lambda_2>0$, let 
$t_j=\sqrt{\lambda_{2j}}$ for $j=1, 2, \dots $, $t_0=0$ and $t_{-j}=-t_j$ and suppose
further that 
$({\hbox{sinc}}(s-t_j))$ is a Riesz basis for 
$RPW(\pi )$. By Theorem 6.1, this event has positive probability with respect to
$\nu_N^\beta$ for suitable $\beta , N>0$. Using the
classical language of divisors, we analyse the addition rule on the expressions
$\sum_{j=-\infty}^\infty x_jg(t_j)$.\par  
\indent The Jacobian of ${\cal E}$ is a complex torus of infinite 
dimension which has real part ${\bf X}$. To construct the map from the divisors on 
${\cal E}$ to 
${\bf X}$, McKean and Trubowitz [22] used sampling on a space of entire functions similar 
to $PW(\pi )$ and thus obtained a suitable family of holomorphic differentials. Their
theory requires smooth $q$, so we extend this to typical $q$ in $(\Omega_N, \nu_N^\beta )$. In this section we define ${\bf X}$ in terms of sampling on $RPW(\pi )$. 
We also show that the Jacobian map has a
Carleman determinant, with a rescaling argument to avoid the formal computations of 
section
12 of [22].\par
\indent We momentarily suppose that $\mu_j=\lambda_{2j}$ for all but finitely many $j$
and then in Theorem 7.3 allow more general assumptions. 
By analogy with classical examples [23], we introduce the generating function
$$S=\sum_{j}2\int_{\mu_j}^{\lambda_{2j}}\sqrt{\Delta (x)^2-4}\,
dx\eqno(7.1)$$
\noindent and introduce the Taylor coefficients of $\Delta
(x)^2-4=\sum_{j=0}^\infty \alpha_jx^j$ as new variables. The
corresponding phases are defined by
$$\varphi_k={{\partial S}\over{\partial \alpha_k}}=\sum_{j=1}^\infty 
\int_{\mu_j}^{\lambda_{2j}}{{x^{k}dx}\over{\sqrt{\Delta
(x)^2-4}}}\qquad (k=0, 1, \dots
).\eqno(7.2)$$
\noindent Then we associate with an entire function
$h(x)=\sum_{k=0}^\infty \beta_kx^k$ the series 
$$\sum_{k=0}^\infty\beta_k\varphi_k =\sum_{j=1}^\infty
\int_{\mu_j}^{\lambda_{2j}} {{h(x)\,dx}\over{\sqrt{\Delta
(x)^2-4}}}.\eqno(7.3)$$
\indent  A real point ${\bf q}_j$ on ${\cal E}$ has the form 
${\bf q}_j=(\mu_j, \varepsilon_j\sqrt{\Delta (\mu_j)^2-4})$ 
where $\lambda_{2j-1}\leq \mu_j\leq \lambda_{2j}$ is in the $j^{th}$ spectral gap, 
$\sqrt{\Delta (\mu_j)^2-4}\geq 0$ and $\varepsilon_j=\pm 1$ where the signs indicate the 
top and bottom of the cut in ${\cal E}$. In particular, let ${\bf p}_j=(\lambda_{2j}, 0)$, 
and, for any subset $M$ of ${\bf N}$, introduce the real divisor 
$\delta =\sum_{j\in M}({\bf p}_j-{\bf q}_j)$. 
The set of such $\delta$ generates a 
free abelian group ${\hbox{Div}}$ under formal addition and subtraction, which we regard as real
divisors of degree zero. When $M$ is finite, we say that $\delta$ has
finite support, and such $\delta$ generate a subgroup ${\hbox{Div}}_0$ 
of ${\hbox{Div}}$.\par
\vskip.05in

\noindent {\bf Lemma 7.1} {\sl Let 
$V=\{ g\in RPW(\pi ): z^2g(z)\in RPW(\pi )\}$ with $\Vert g\Vert^2_V=\Vert g\Vert^2_{L^2}+\Vert
z^2g(z)\Vert^2_{L^2}$ have dual space $V^*$ with respect to the pairing $\langle f,g\rangle_{L^2}$.} \par
\noindent {\sl (i) Then for all} $\delta\in {\hbox{Div}}$ {\sl there exists 
$\psi_\delta\in V^*$, given by}
$$\psi_\delta (g)= \sum_{j\in M}\varepsilon_j\int_{\mu_{j}}^{\lambda_{2j}} 
{{(g(\sqrt{x}) +g(-\sqrt{x}))dx}\over{\sqrt{\Delta (x)^2-4}}},\eqno(7.4)$$
\noindent {\sl such that $J:\delta\mapsto\psi_\delta$ is a group homomorphism} ${\hbox{Div}}\rightarrow V^*$.\par 
\noindent {\sl (ii) Let} $g_\varepsilon
(z)=2^{-1}\Delta'((1-6\varepsilon )z^2)({\hbox{sinc}}\,(\varepsilon
z))^3.$ {\sl Then
$g_\varepsilon\in V$ for all
$\varepsilon >0$ and the limit}
$$\psi_\delta (g_\varepsilon )\rightarrow
\omega_\infty (\delta )=\sum_{j\in M}\varepsilon_j\int^{\lambda_{2j}}_{\mu_j}{{\Delta'(x)\, dx}
\over{\sqrt{\Delta (x)^2-4}}}\qquad (\varepsilon\rightarrow 0+),\eqno(7.5)$$
\noindent {\sl exists for all} $\delta\in{\hbox{Div}}_0$, {\sl and defines a group homomorphism} $\omega_\infty : {\hbox{Div}}_0\rightarrow {\bf R}$.\par 
\noindent {\sl (iii) For any sequence of signs $\varepsilon_j=\pm 1$, there is a map} 
${\cal M}_\lambda\rightarrow {\hbox{Div}}$ {\sl given by $q\mapsto 
\sum_{j=1}^\infty ({\bf p}_j-{\bf q}_j)$ where $(\mu_j)_{j=1}^\infty$ is the tied spectrum of (1.1)
for $q$ and ${\bf q}_j=(\mu_j, \varepsilon_j\sqrt{\Delta (\mu_j)^2-4})$.}\par 
\vskip.05in
\noindent {\bf Proof.} (i) Note that
$g(\sqrt{x})+g(-\sqrt{x})=(2/\pi) \int_0^\pi \cos (s\sqrt {x})\Re
\hat g(s)\, ds$, so the numerator in (7.4) is entire of order $1/2$. Now ${\hbox{Div}}$ gives a subgroup of $V^*$ since $\psi_\delta$ defines a bounded linear functional on $V$; indeed, we can bound the $j^{th}$ summand by a
constant multiple of
$$\int_{\mu_{j}}^{\lambda_{2j}}
{{\vert g(\sqrt{x}) +g(-\sqrt{x})\vert dx}
 \over{\sqrt{(1-x/\lambda_{2j})(x/\lambda_{2j-1}-1)}}},\eqno(7.6)$$
\noindent where $\sum_{j=1}^\infty j^2\vert g(\sqrt{\lambda_{2j}})\vert$
 converges by (3.25) and the Cauchy--Schwarz inequality. Thus there exists a uniquely
determined sequence $(x_k(\delta
))_{k=-\infty}^\infty$ such that 
$(x_k(\delta )/(1+k^2))_{k=-\infty}^\infty \in \ell^2$ and
$$\sum_{k=-\infty}^\infty x_k(\delta )g(t_k)=\sum_{j=1}^\infty 
\varepsilon_j\int_{\mu_{j}}^{\lambda_{2j}}
 {{(g(\sqrt{x}) +g(-\sqrt{x}))dx}\over{\sqrt{\Delta (x)^2-4}}}
 \qquad (g\in V),\eqno(7.7)$$ 
\noindent so with respect to this basis  
$\psi_\delta =\sum_{j=-\infty}^\infty x_j(\delta )\,{\hbox{sinc}}(s-t_j)\in V^*$. This correspondence respects the group law, so that if $\delta\mapsto (x_j(\delta ))$ and
$\varepsilon\mapsto (y_j(\varepsilon ))$, then $\delta -\varepsilon\mapsto
(x_j(\delta )-y_j(\varepsilon ))$. For $\delta\in {\hbox{Div}}_0$, we can recover 
the coefficients by operating on the biorthogonal functions, obtaining
$\psi_\delta (g_j)=x_j(\delta )$.\par
\indent In particular, there exist bounded linear functionals
$$\psi_{(\lambda_{2j},0)-(\lambda_{2j-1},0)}:g\mapsto 2\int_{\lambda_{2j-1}}^{\lambda_{2j}}
{{g(\sqrt{x}) +g(-\sqrt{x})} \over{\sqrt{\Delta (x)^2-4}}}dx\qquad 
(g\in V)\eqno(7.8)$$
\noindent classically known as the real periods. By extension, 
 $\sum_{j\in M}\psi_{(\lambda_{2j},0)-(\lambda_{2j-1},0)}$ is also a bounded linear functional
on $V$ for all subsets $M$ of ${\bf N}$. Let $\Lambda$ be the lattice in $V^*$ that is generated by  
$\{\sum_{j\in M}\psi_{(\lambda_{2j},0)-(\lambda_{2j-1},0)}: 
M\subseteq {\bf N}\}$.\par
\indent (ii) By Theorem 2.1 of [20], all of the zeros of $4-\Delta (z)^2$ are real, so by
Laguerre's theorem, all the zeros of $\Delta'(z)$ are also real, and
separated by the zeros of $4-\Delta (z)^2$; see [10, p. 264].  Hence from the 
resulting product
representation,  $\Delta'(z^2)$ is
even and entire of exponential type; also $\Delta'(z)$ is real for all real $z$. 
By (3.9) and the Cauchy integral formula, we deduce that $\Delta'(z^2)$ is bounded on
${\bf R}$ as in [10, p
264], so $x^2g_\varepsilon (x), g_\varepsilon (x)\in L^2({\bf R})$; hence $zg_\varepsilon (z),
g_\varepsilon (z)\in RPW(\pi )$. Also, $g_\varepsilon (z)
\rightarrow 2^{-1}\Delta'(z^2)$ as
$\varepsilon \rightarrow 0+$, uniformly on compact subsets of ${\bf
C}$. We observe that, for all $\delta,\eta\in {\hbox{Div}}_0$, we have $\omega_\infty (\delta +\eta )=\omega_\infty
(\delta )+\omega_\infty (\eta )$.\par
\indent (iii) The tied spectrum interlaces the 
periodic spectrum so $\lambda_{2j-1}\leq \mu_j\leq \lambda_{2j}$ and we can
 apply (i). In the following result, we will deal with the ambiguity associated with the choice of
$\varepsilon_j=\mp 1$. {$\square$}\par 
\vskip.05in

\noindent {\bf Definition} (Jacobian) Regarding ${\hbox{Div}}$ as a subgroup of
$V^*$ under the map $J$ of Lemma 7.1(i), we define the real 
Jacobian of ${\cal E}$ to be the abelian group ${\bf X}={\hbox{Div}}/\Lambda$.
Then $J(\delta )$ has coordinates 
$J(\delta )=(x_j(\delta ))_{j=-\infty}^\infty$ with respect to the
sampling sequence. Also let ${\bf X}_0={\hbox{Div}}_0/(\Lambda\cap {\hbox{Div}}_0)$.\par 
\vskip.05in
\indent The inverse spectral problem involves recovering $q$ from the 
spectral data consisting of the periodic spectrum 
$\lambda =\{ \lambda_j: j=0, 1,\dots\}$ and the family of tied spectra 
$\{ \mu_j(s):j\in {\bf N};
s\in [0, 2\pi ]\}$ of the translated potentials $q(x+s)$ 
for $s\in [0, 2\pi ]$. This gives $({\bf p}_j-{\bf
q}_j(s))\in {\hbox{Div}}$, where 
${\bf q}_j(s)=(\mu_j(s), \varepsilon\sqrt{\Delta (\mu_j(s))^2-4})$ is a
real point on ${\cal E}$.\par
\vskip.05in
\noindent {\bf Proposition 7.2} {\sl (i) The real periods of $\omega_\infty$
vanish, so that $\omega_\infty ((\lambda_{2j},0)-(\lambda_{2j-1},0))=0$
for all $j$. Hence $\omega_\infty$ induces a 
homomorphism ${\bf X}_0\rightarrow {\bf R}.$}\par
\noindent {\sl (ii) The functions 
$q(x+s)\mapsto \omega_\infty ({\bf p}_j-{\bf q}_j(s))$ are continuous
$\Omega_N\rightarrow {\bf R}$. The mean value of $(d/ds)\omega_\infty ({\bf
p}_j-{\bf q}_j(s))$ with respect to $\nu_N^\beta$ is zero.}\par 
\noindent {\sl (iii) The map $q\mapsto (\omega_\infty ({\bf p}_j-{\bf
q}_j(s)):j=1, 2,\dots ; s\in [0, 2\pi ])$ is one-to-one on ${\cal M}_\lambda$.}\par  
\vskip.05in
\noindent {\bf Proof.}  (i) Observe that $\Delta (\lambda_{2j})=\Delta (\lambda_{2j-1})$, so
$\int_{\lambda_{2j-1}}^{\lambda_{2j}}\Delta'(\lambda )d\lambda /
\sqrt{\Delta (\lambda )^2-4}=0$. The rest follows from Lemma
7.1(ii).\par
\indent (ii) The map $q\mapsto \Delta$ is continuous from the norm
topology to the uniform topology on compact planar sets; also $q\mapsto
\mu_j$ is continuous by Lemma 4.1. Hence $q\mapsto \cosh^{-1}\Delta
(\mu_j)$ is continuous. By (2.2), the measure $\nu_N^\beta$ is invariant under the translation $q(x)\mapsto
q(x+s)$, and the translation preserves the periodic spectrum. Hence by linearity 
the mean 
of $\Delta'(\mu_j(s))\mu_j'(s)/\sqrt{\Delta (\mu_j(s))^2-4}$ with respect to
$\nu_N^\beta$ equals the mean of 
$$\int_0^{2\pi} {{\Delta'(\mu_j(s))\mu_j'(s)ds}\over{\sqrt{\Delta
(\mu_j(s))^2-4}}}.\eqno(7.9)$$
\noindent Now as $s$ ranges over $[0, 2\pi ]$, $\mu_j(s)$ describes
$[\lambda_{2j-1}, \lambda_{2j}]$ back and forth an integral number of times,
returning to $\mu_j(0)$. So by (i), (7.9) equals zero, hence the mean of $(d/ds)\omega_\infty ({\bf
p}_j-{\bf q}_j(s))$ equals zero.\par
\indent (iii) Let $\lambda_j'\in (\lambda_{2j-1}, \lambda_{2j})$ satisfy $\Delta'(\lambda_j')=0$. Then the graph of  
$s\mapsto\omega_\infty ({\bf p}_j-{\bf q}_j(s))$ crosses the axis when
$\mu_j=\lambda_{2j-1}$ and has height $\cosh^{-1}(\Delta'(\lambda_j')/2)$, known as
the spike height. In turn, the points $\mu_j(s)$ of the tied spectrum are determined
by $\cosh^{-1}(\Delta (\mu_j(s))/2)$.\par
\indent Given the periodic and tied spectra,  we note that by page 329 of [31] for a simple periodic eigenvalue
$\lambda_k$, there exists a corresponding periodic eigenfunction $f$ such that
$$-\Delta'(\lambda_k)f^2(s) =\prod_{j=1}^\infty
{{4(\mu_j(s)-\lambda_k)}\over{j^2}}.\eqno(7.10)$$
\noindent Finally, one can in principle recover $q$ from
$$4f^4(s)(q(s)-\lambda_k) =2f^2(s)(f^2(s))''-((f^2(s))')^2.\eqno(7.11)$$
\rightline {$\square$}\par
\vskip.05in
\noindent {\bf Remark}. Whereas typical $q\in (\Omega_N, \nu_N^\beta )$ are unbounded, one can
improve upon Proposition 7.2(iii) in the case in which $q$
is smooth. Then the differences 
$\lambda_{2j}-\lambda_{2j-1}$ are rapidly decreasing as $j\rightarrow\infty$, and
the spike heights are then summable. Trubowitz [31] showed that $q$ may be recovered from the
spectral data via
$$q(s)-\int_0^{2\pi}q(x) {{dx}\over{2\pi}}=-2{{d}\over{ds}}\sum_{j=1}^\infty
\omega_\infty({\bf p}_j-{\bf q}_j(s))\qquad (s\in [0, 2\pi ]).\eqno(7.12)$$
\noindent Equivalently, $q$ can be recovered from the periodic spectrum,
the tied spectrum and so-called norming constants.\par
\indent In the classical theory of finite genus [16, p.64], the Jacobian map has a 
nonzero determinant, under certain conditions on the divisors. 
The following Theorem 7.3 introduces a determinant and gives criterion for the linear statistic associated with sampling at
$(t_j)$ to arise from
a divisor on ${\cal E}$ in the sense of Lemma 7.1.\par
\indent By construction, ${\bf X}$ is a infinite-dimensional torus, and has a
system of real coordinates. We introduce $\eta_j^0\in [\lambda_{2j-1},
\lambda_{2j}]$ and then $\eta_j(s_j)\in [\lambda_{2j-1}, \lambda_{2j}]$
by the condition
$$\int_{\eta_j(s_j)}^{\lambda_{2j}}{{dx}\over{\sqrt{\Delta
(x)^2-4}}}=1-s_j\eqno(7.13)$$
\noindent where $\eta_j(0)=\eta_j^0$ and $s_j$ is a new real variable. Let  
$$L_0=\Bigl( \sum_{j=1}^\infty
\Bigl({{\lambda_{2j}-\eta_j^0}\over{\sqrt{\lambda_{2j}}
+\sqrt{\eta_{j}^0}}}\Bigr)^2\Bigr)^{1/2}.\eqno(7.14)$$
\noindent We recall from section 6 the biorthogonals $(g_k)_{k\in {\bf 
Z}}$ of the Riesz basis $({\hbox{sinc}}(s-t_k))_{k\in {\bf Z}}.$ For $\sigma =(s_j)_{j=1}^\infty$ 
let $X(\sigma )=J(\sum_j\delta_{\lambda_{2j}}-\delta_{\mu_j(s_j)})-(1)$, so
$X(\sigma )=(X_k(\sigma ))_{k=-\infty}^\infty$ where 
$$X_k(\sigma )=-1+\sum_{j=1}^\infty \int_{\eta_j(s_j)}^{\lambda_{2j}}
{{g_k(\sqrt{x})+g_k(-\sqrt{x})}\over{\sqrt{\Delta
(x)^2-4}}}dx.\eqno(7.15)$$
\vskip.05in
\noindent {\bf Theorem 7.3} {\sl (i) There exists $C_\infty$ such that $\Vert
X(0)\Vert_{\ell^\infty} \leq C_\infty\sqrt{\pi }M_1L_0$, where $\int L_0d\nu_N^\beta \rightarrow 0$ as $N\rightarrow 0$. The Fr\'echet derivative $X'(0)$ defines a
bounded linear operator on $\ell^\infty$.\par
\noindent (ii) On $\ell^2$, the operator $X'(0)-I$ is Hilbert--Schmidt and exists 
$C_2$ such that the Carleman 
determinant satisfies} 
$$\bigl\vert\det_2X'(0)-1\bigr\vert\leq C_2M_1L_0.\eqno(7.16)$$
\noindent {\sl (iii) Suppose there exists
$\sigma= (s_j)$ in
the unit ball of $\ell^\infty$ such that $X(\sigma )=0$. Then 
$$\sum_{j=-\infty}^\infty g(t_j)=\sum_{j=1}^\infty   
\int_{\eta_j(s_j)}^{\lambda_{2j}}
{{g(\sqrt{x})+g(-\sqrt{x})}\over{\sqrt{\Delta
(x)^2-4}}}dx\eqno(7.17)$$
\noindent for all $g$ in the linear span of the $(g_k)_{k\in {\bf 
Z}}$}.\par
\vskip.05in
\noindent {\bf Proof.} (i) We observe that
$$L_0\leq \Bigl( \sum_{j=1}^\infty
\Bigl({{\lambda_{2j}-\lambda_{2j-1}}\over{\sqrt{\lambda_{2j}}
+\sqrt{\lambda_{2j-1}}}}\Bigr)^2\Bigr)^{1/2},$$
\noindent which converges by Proposition 3.1. Also
$\lambda_{2j}-\lambda_{2j-1}\rightarrow 0$ as $N\rightarrow 0$, so
$L_0\rightarrow 0$. By the mean value theorem for integrals, there exists
$\nu_{j,k}\in (\eta_j^0, \lambda_{2j})$ such that 
$$X_k(0)=\sum_{j=1}^\infty \bigl( g_k(\sqrt{\nu_{j,k}})+  
g_k(-\sqrt{\nu_{j,k}})-
g_k(\sqrt{\lambda_{2j}})-g_k(-\sqrt{\lambda_{2j}})\bigr)
\int_{\eta_j^0}^{\lambda_{2j}} {{dx}\over{\sqrt{\Delta
(x)^2-4}}},\eqno(7.18)$$
\noindent where all of the integrals are equal to unity. Then by the
mean value theorem, there exist $\omega_{j,k}\in (\nu_{j,k},
\lambda_{2j})$ such that 
$$X_k(0)=\sum_{j=0}^\infty \bigl(
g_k'(\sqrt{\omega_{j,k}})-g_k'(-\sqrt{\omega_{j,k}})\bigr)
(\sqrt{\omega_{j,k}}-\sqrt{\lambda_{2j}}),\eqno(7.19)$$
\noindent so by the Cauchy--Schwarz inequality
$$\vert X_k(0)\vert \leq\Bigl(\sum_{j=1}^\infty
(g_k'(\sqrt{\omega_{j,k}})^2+g_k'(-\sqrt{\omega_{j,k}})^2\Bigr)^{1/2}
\Bigl( \sum_{j=1}^\infty (\sqrt{\lambda_{2j}}-
\sqrt{\omega_{j,k}})^2\Bigr)^{1/2},\eqno(7.20)$$
\noindent and since $(\pm \sqrt{\omega_{j,k}})$ gives a sampling
sequence for $RPW(\pi ) $ we can choose $C_\infty $ independent of $k$ such
that 
$$\vert X_k(0)\vert \leq C_\infty\Vert g_k'\Vert_{L^2}\Bigl( \sum_{j=1}^\infty
\Bigl(
{{\lambda_{2j}-\eta_j^0}\over{\sqrt{\lambda_{2j}}+\sqrt{\eta_j^0}}}
\Bigr)^2\Bigr)^{1/2};\eqno(7.21)$$
\noindent hence $X(0)=(X_k(0))$ is bounded and $\Vert
X(0)\Vert_{\ell^\infty} \leq C_\infty\sqrt{\pi }M_1L_0$.\par
\indent The function $X:{\hbox{Ball}}(\ell^\infty )\rightarrow
\ell^\infty$ is Fr\'echet differentiable near to $s=0$, and the 
derivative is expressed as an infinite matrix with respect to the usual
weak${}^*$ basis 
$$X'(s)=\Bigl[ {{\partial X_k}\over{\partial s_j}}\Bigr]_{j,k}=\bigl[
g_k(\sqrt{\eta_j(s_j)}+      
g_k(-\sqrt{\eta_j(s_j)}\bigr],\eqno(7.22)$$
\noindent so in particular 
$$X'(0)-[\delta_{j, \pm k}]=\bigl[ 
 g_k(\sqrt{\eta_j^0})+g_k(-\sqrt{\eta_j^0})-g_k(\sqrt{\lambda_{2j}})-
 g_k(-\sqrt{\lambda_{2j}})\bigr].\eqno(7.23)$$    
\noindent The rows of this matrix are absolutely summable with
uniformly bounded sums, so as in (7.20), $X'(0)$ defines a bounded linear operator on
$\ell^\infty$. \par 
\indent (ii) Furthermore, $(g_k)$ is itself a Riesz basis, and hence for any pair
of real sequences $(u_j)$ and $(v_j)$ such that $(u_j-v_j)\in \ell^2$, we have
$$\sum_{j,k}\vert
g_k(u_j)-g_k(v_j)\vert^2\leq M_1^2\sum_j \bigl( u_j-v_j\bigr)^2,\eqno(7.24)$$
\noindent as in (3.21), where the final series converges. In particular, we can
take $u_j=\sqrt{\eta_j^0}$ and $v_j=\sqrt{\lambda_{2j}}$, so that
$g_k(v_j)=\delta_{k,j}$. 
Hence $X'(0)-I$ is a Hilbert--Schmidt operator, with norm bounded by a
constant
multiple of $M_1L_0$. Hence $X'(0)$ has a
Carleman determinant, and $\det_2$ is a Lipschitz continuous
function on bounded subsets of $I+HS$. (We have not quite proved that $X'(0)$ has a Hill's 
determinant as in [20, p. 29], since (7.20) involves sums of squares.)\par

\indent  (iii) We observe that the biorthogonal system
satisfies
$$\sum_{j=-\infty}^\infty g_k(t_j)=\sum_{j=-\infty}^\infty \langle
g_k(x), {\hbox{sinc}} (x-t_j)\rangle_{L^2}=1,\eqno(7.25)$$
\noindent since only one term in the sum is nonzero. Hence the condition
$X_k(s)=0$ gives the identity (7.17) for $g=g_k$, and the general case
follows by linearity. {$\square$}\par
\vskip.05in
\noindent {\bf Corollary 7.4} {\sl (i) Suppose  
$C_\infty\sqrt{\pi }M_1L_0<1/3$. Then $X'(0)$ defines an invertible linear operator on
$\ell^\infty$.\par
\noindent (ii) Suppose further that} 
$$\Vert X(s)-X(0)-X'(0)s\Vert_{\ell^\infty }\leq (1/6)\Vert
s\Vert_{\ell^\infty }\qquad (s\in {\hbox{Ball}}_1(\ell^\infty
)).\eqno(7.26)$$ 
\noindent {\sl Then the sequence $(\sigma_n)_{n=1}^\infty$ produced 
by Newton's modified algorithm 
$\sigma_0=0$ and
$$\sigma_{n+1}=\sigma_n-X'(0)^{-1}X(\sigma_n),\qquad (n=0,
1, \dots )\eqno(7.27)$$
\noindent converges to $\sigma$ such that $X(\sigma )=0,$ so Theorem 7.3(iii) holds.}\par
\vskip.05in
\noindent {\bf Proof.} (i) When $C_\infty\sqrt{\pi}M_1L_0<1$, the 
operator
$X'(0)$ on $\ell^\infty$ satisfies $\Vert X'(0)-I\Vert <1$ by
Theorem 7.2(i), and hence $X'(0)$ is invertible with $\Vert
X'(0)^{-1}\Vert\leq (1-C_\infty\sqrt{\pi }M_1L_0)^{-1}$.\par
\indent (ii) This follows from (iv) by Corollary 2 of [1]. {$\square$}\par
\vskip.05in
\noindent {\bf Acknowledgements} The first-named author carried out part of this 
work during a visit to the University of New South Wales. He also thanks H.P. McKean for a
motivating discussion.  
The work of the second-named author was partially supported by an 
EPSRC research studentship. The work of the third-named author was partially supported by
an LMS Scheme 2 grant.\par
\vskip.05in
\noindent {\bf References}\par
\noindent [1] R.G. Bartle, Newton's method in Banach spaces, {Proc. 
Amer. Math. Soc.} {6} (1955), 827--831.\par
\noindent [2] G. Blower, Almost sure weak convergence for the
generalized orthogonal ensemble, {J. Statist. Phys.} {105}
(2001),  309--335.\par
\noindent [3] G. Blower, A logarithmic Sobolev inequality for the invariant 
measure of the periodic Korteweg--de Vries equation, {Stochastics} {84} (2012), 533-542.\par
\noindent [4] G. Blower, C. Brett and I. Doust,
 Logarithmic Sobolev inequalities and spectral concentration for the 
periodic Schr\"odinger equation, {Stochastics} (2014) {\sl to appear}.\par
\noindent [5] S.G. Bobkov, I. Gentil and M. Ledoux, Hypercontractivity of 
Hamilton--Jacobi equations, {J. Math. Pures Appl. (9)} {80} (2001), 669-696.\par
\noindent [6] J. Bourgain, Fourier transform restriction phenomena for certain lattice
 subsets and applications to nonlinear evolution equations II: 
The KdV equation, {Geom. Funct. Anal.} {3} (1993), 209-262.\par
\noindent [7] J. Bourgain, Periodic nonlinear Schr\"odinger equation and
 invariant measures, {Comm. Math. Phys.} {166} (1994), 1--26.\par
\noindent [8] A. Erdelyi, \"Uber die freien Schwingungen in Kondensatorkreisen mit periodisch\par
\noindent  ver\"andlicher Kapaziteier, {Ann. der Physik} {5} folge band 19 (1934), 585-622.\par
\noindent [9] C. S. Gardner, J. M. Greene, M. D. Kruskal, and 
R. M. Miura, Korteweg-de Vries equation and generalization VI.
Methods for exact solution. {Comm. Pure Appl. Math.} {27} (1974),
97--133.\par
\noindent [10] J. Garnett and E. Trubowitz, Gaps and bands of one-dimensional
periodic Schr\"odinger operators, {Comment. Math. Helvetici} {59}
(1984), 258--312.\par 
\noindent [11] J.R. Higgins, Five short stories about the cardinal series, {Bull. Amer. Math. Soc.
NS} {12} (1985), 45--89.\par  
\noindent [12] K. Johansson, On random matrices from the
compact classical groups, 
{\sl Ann. of Math. (2)} {145} (1997), 519--545.\par 
\noindent [13] T. Kappeler and B. Mityagin,
Gap estimates of the spectrum of Hill's equation and action variables for KdV, {Trans. Amer. Math. Soc.} 
{351} (1999), 619--646.\par
\noindent [14] N.M. Katz and P. Sarnak, {Random matrices, Frobenius eigenvalues and monodromy}, (American Mathematical Society, Providence, R.I. 1999).\par
\noindent [15] P. Koosis, {The logarithmic integral I}, (Cambridge University Press, 1988).\par
\noindent [16] S. Lang, Introduction to algebraic and abelian functions,
Second Edition, (Springer-Verlag, New York, 1982).\par
\noindent [17] P. D. Lax, Integrals of nonlinear equations of evolution and solitary waves,
{Comm. Pure Appl. Math.} {21} (1968), 467--490. 

\noindent [18] J.L. Lebowitz, Ph. Mouniax and W.-M. Wang, Approach to 
equilibrium for the stochastic NLS, {Comm. Math. Phys.} {321} (2013), 69--84.\par 

\noindent [19] A. Lytova and L Pastur, Central Limit Theorem for linear eigenvalue 
statistics of random matrices with independent entries, {Ann. Probab.} {37} (2009), 1778--1840.\par

\noindent [20] W. Magnus and S. Winkler, {Hill's Equation}, second edition (Dover Publications, New York, 1979).\par
\noindent [21] H.P. McKean and P. van Moerbeke, The spectrum of Hill's
equation, {Invent. Math.} {30} (1975), 217-274.\par
\noindent [22] H.P. McKean and E. Trubowitz, Hill's operator and hyperelliptic 
function theory in the presence of infinitely many branch points, {Comm. Pure Appl. 
Math.} {29} (1976), 143--226.\par
\noindent [23] J. Moser, Various aspects of integrable Hamiltonian systems, 233--289,
 {Dynamical Systems}, (Birkhauser, Boston, 1980).\par
\noindent [24] N.K. Nikolski, {Operators, Functions and Systems: 
An Easy Reading Vol. 2: 
Model Operators and Systems}, (American Mathematical Society, Providence, R.I., 2002).\par
\noindent [25] J. Ortega--Cerd\`a and K. Seip, {Fourier frames}, {Ann. of  Math. (2)} 
{155} (2002), 789--806.\par
\noindent [26] J.R. Partington, {Interpolation, identification and sampling}
 (Oxford University Press, 1997).\par
\noindent [27] H.L. Pedersen, Entire functions having small logarithmic sums over 
certain discrete subsets, {Ark. Mat.} {36} (1998), 119--130.\par
\noindent [28] B. Simon, {Trace ideals and their applications}, (Cambridge University Press,
1979).\par
\noindent [29] A Soshnikov, The Central Limit Theorem for local linear statistics in
 classical compact groups and related combinatorial identities, {Ann. Probab.} {28} (2000), 1353--1370.

\noindent [30] C.A. Tracy and H. Widom, Correlation functions, cluster functions and 
spacing distributions for random matrices, {J. Statist. Phys.} {92} (1998), 809--835.\par
\noindent [31] E. Trubowitz, The inverse problem for periodic
potentials, {Comm. Pure Appl. Math.} {30} (1977), 321-337.\par
\noindent [32] C. Villani, {Topics in Optimal Transportation} 
(American Mathematical Society, Providence R.I., 2003).\par
\noindent [33] C. Villani, {Optimal transport, old and new},
(Springer-Verlag, 2009).\par
\vfill
\eject
\end

\indent Let $\tau =(\tau_j)_{j=-\infty}^\infty$ where
$$\tau_j =\cases{ \sqrt{\lambda_{2j}}& for $j=1,2,\dots $;\cr
0, & for $j=0$;\cr
-\sqrt{ \lambda_{-2j}}& for $j=-1,-2, \dots $,\cr}\eqno(5.20)$$
\noindent and for $g\in RPW(\pi )$ introduce the linear statistic
$$F(\tau )=\sum_{j=-\infty}^\infty (g(\tau_j)-g(j)).\eqno(5.21)$$

\noindent {\sl Furthermore, $(s_j)_{j=-\infty}^\infty$
is also a sampling sequence for constants $A/2$ and $2B$ whenever
$(s_j-t_j)\in \ell^2$ satisfies $\Vert (s_j-t_j)\Vert_{\ell^2}\leq
(\sqrt{2}-1)\sqrt{3A}/\pi$.}\par
\indent There is an important case in which we can compute the linear statistics and 
their mean values asymptotically. Consider the tied spectrum $(\mu_j)_{j=1}^\infty$ for 
$-f''+qf=\lambda f$ with $f(0)=f(2\pi )=0$. There are constants
 $c_j$, depending upon $q$ such that $\mu_j=j^2/4+c_0+4c_2/j^2+\dots ,$ and 
one can invert this asymptotic series to contain $j^2/4=\psi (\mu )$. Gelfand [15] showed that there is 
 an asymptotic expansion
$$\sum_{j=1}^\infty \bigl(\mu_j+\zeta \bigr)^{-1}=\sum_{k=0}^\infty 
{{(-1)^kA_k(q)}\over{\zeta^{k+1}}}+{{d}\over{d\zeta}}\Bigl( -\psi (-\zeta )\Bigr)^{1/2}
\qquad (\zeta\rightarrow\infty) ;\eqno(4.1)$$
 \noindent where for each $k$ there exists a real polynomial $p_k$ with coefficients depending upon $q$,
 such that $A_k(q)=\sum_{j=0}^{\infty} (\mu _j^{k}-p_k(j^2))$ is a sum of powers of the eigenvalues, modified
 to ensure convergence. In particular, the first two modified moments are
$$\eqalignno{A_1(q)&=\sum_{j=1}^\infty\Bigl(\mu_j-{{j^2}\over{4}}-c_0\Bigr)-{{c_0}\over{2}},\cr
A_2(q)&=\sum_{j=1}^\infty\Bigl(\mu_j^2-{{j^4}\over{16}}-{{j^2c_0}\over{2}}-2c_1-c_0^2\Bigr) -
{{2c_2+c_0^2}\over{2}}.&(4.2)\cr}$$
\indent One can view this as an asymptotic expansion of linear statistics for $g_\zeta (x)=1/(x^2+\zeta )$.
 Brownian loop has almost surely continuous sample paths, so we can form pointwise products of the random $q(x_j)$
 and introduce the family of random variables $\Omega_N\rightarrow {\bf R}$ $q\mapsto q(x)$  for $x\in {\bf T}$,
 which have moments of order $k$ with respect to $\nu_N^\beta$ by
$$Q_k(x_1, \dots , x_k) =\int_{\Omega_N} q(x_1)\dots q(x_k)\nu_N^\beta (dq).\eqno(4.3)$$
\noindent It is easy to see that $Q_k$ reduces to a function of $k-1$ variables.\par 
\vskip.05in
\noindent {\bf Proposition 4.1} {\sl The moments of $q$ determine the mean values of the modified moments of the
 tied eigenvalues via the asymptotic series}
$$\eqalignno{ \int_{\Omega_N} \sum_{j=1}^\infty\bigl(\mu_j+\zeta \bigr)^{-1}\nu_N^\beta (dq)
&=\Bigl({{\pi\coth 2\pi\sqrt{\zeta}}\over{\sqrt{\zeta}}}-{{1}\over{2\zeta}}\Bigr)-
{{Q_1}\over{(2\sqrt{\zeta})^3}}\Bigl( 4\pi -{{1}\over{\sqrt{\zeta}}}\Bigr)&(4.4)\cr
&+{{1}\over{(2\sqrt{\zeta})^3}}\int_0^{2\pi} Q_2(v,0)e^{-2v\sqrt{\zeta}}\Bigl( 
{{2\pi -v}\over{\sqrt{\zeta}}}+(2\pi -v)v-{{1}\over{2\zeta}}\Bigr)dv+\dots \cr}$$ 
\noindent {\sl where $Q_1$ is the constant that satisfies} 
$$\int_{\Omega_N} A_1(q)\nu_N^\beta (dq)=-{{Q_1}\over {8}},\eqno(4.5)$$
\noindent {\sl and $Q_2(v,0)=Q_2(0,0)+O(v)$ as $v\rightarrow 0+.$}
\vskip.05in
\noindent {\bf Proof.} 
\noindent so that $G_\zeta (x,y)=(2\sqrt{\zeta} )^{-1}e^{-\vert x-y\vert\sqrt{\zeta}}+O(e^{-2\pi \sqrt{\zeta}})$ 
as $\zeta\rightarrow\infty$. Let $M_q$ be the operator in $L^2$ of multiplication by $q$, and note that $-d^2/dx^2+q+
\zeta$ can be viewed as a perturbation of $-d^2/dx^2+\zeta$, and that both operators have trace class inverses. Hence 
there is an expansion
$$\eqalignno{\sum_{j=1}^\infty (\zeta +\mu_j)^{-1}&={\hbox{trace}}\Bigl( \bigl( -{{d^2}\over{dx^2}} +\zeta +q\bigr)^{-1}
\Bigr)\cr
&={{\pi}\over{\sqrt{\zeta }}}\cosh 2\pi \sqrt{\zeta} +\sum_{k=1}^\infty (-1)^{k} {\hbox{trace}}
\Bigl(G_\zeta (M_qG_\zeta )^k\Bigr), &(4.7)\cr}$$
\noindent which converges for large $\zeta$ since $\Vert
M_qG_\zeta\Vert^2_{HS}\leq CN/\zeta$ for some $C>0$ and 
$${\hbox{trace}}\Bigl( G_\zeta(M_qG_\zeta )^k\Bigr) 
=O\Bigl( {{4^kN^{k/2}}\over{\zeta^{k/2}}}\Bigr)\qquad (\zeta \rightarrow\infty ).\eqno(4.8)$$ 
\indent Each term in the series is a Lipschitz function of $q\in
\Omega_N$. Indeed, the function $\Omega_N\rightarrow \Omega_N^k:$
$q\mapsto (q,\dots ,q)$ is Lipschitz, as is the function
$\Omega_N^k\rightarrow {\bf R}$  given by
$$(q_1, \dots , q_k)\mapsto \int_{{\bf T}^{k+1}} G_\zeta
(x_0,x_1)q_1(x_1)G_\zeta (x_2,x_3)q_2(x_2)\dots q(x_k)G_\zeta
(x_k,x_0)dx_0\dots dx_k.\eqno(4.9)$$
\indent We can take the expectation of this series by integrating with respect to $\nu_N^\beta$ and changing order of the integration.  The only random term in the $k^{th}$ summand is $q$, so we combine terms into the moment $Q_k$. Now $Q_k$ is symmetrical with respect to permutation of its variables and satisfies 
$$\int_{{\bf T}^k} Q_k(x_1, \dots , x_k)^2 dx_1\dots dx_k\leq N^k\eqno(4.10)$$
\noindent and the mean value of the $k^{th}$ coefficient in the series expansion is
$$\eqalignno{&(-1)^k\int_{\Omega_N} {\hbox{trace}}\bigl( G_\zeta (M_qG_\zeta )^k\bigr) \nu_N^\beta (dq) \cr
&=(-1)^k\int_{{\bf T}^{k+1}} G_\zeta (x_0, x_1)G_\zeta (x_1,x_2) 
\dots G_\zeta (x_k, x_0)Q_k(x_1, \dots , x_k) dx_0\dots dx_k,&(4.11)\cr}$$
\noindent and this term is asymptotic to
$${{1}\over{(2\sqrt{\zeta})^{k+1}}}\int_{{\bf T}^{k+1}}Q_k(x_1, \dots ,x_k) e^{-\sqrt{\zeta}
\vert x_0-x_1\vert -\dots -\sqrt{\zeta}\vert x_k-x_0\vert} dx_0\dots dx_k.
\eqno(4.12)$$
\noindent Gelfand shows that the series expansion and the asymptotic expansion determine one another, since they differ by terms that decay exponentially as $\zeta\rightarrow\infty$.\par
\indent The measure $\nu_N^\beta$ is invariant with respect to rotation of the circle, so in particular, for $k=1$ we obtain
$${{-1}\over{(2\sqrt{\zeta})^3}}\int_{\Omega_N}q(e^{0})\nu_N^\beta (dq)
 \bigl(4\pi -{{1}\over{\sqrt{\zeta}}}+e^{-4\pi \sqrt{\zeta}} \bigr);
\eqno(4.13)$$ 
\noindent also, for $k=2$ we can reduce (4.12) to a one-dimensional integral
$${{1}\over{(2\sqrt{\zeta})^3}}\int_0^{2\pi}Q_2(v,0)\Bigl( 
(2\pi -v){{e^{-2 v\sqrt{\zeta}}}\over{\sqrt{\zeta}}}+(2\pi -v) 
ve^{-2v\sqrt{\zeta}}+{{e^{-4\pi \sqrt{\zeta}}}\over{2\zeta}}-
{{e^{-2v\sqrt{\zeta}}}\over{2\zeta}}\Bigr)dv.\eqno(4.14)$$
\noindent By comparing the coefficients of $\zeta^{-2}$ in the series, we deduce that $\int_{\Omega} A_1\nu_N^\beta (dq)=-Q_1/8$.\par
\noindent Note that $Q_1$ is a constant equal to the mean value of $\int_0^{2\pi} q(x)dx/(2\pi)$ 
with respect to $\nu_N^\beta (dq).$ 
When $\beta =0$, one can use simple symmetry arguments to show that $Q_1=0$. 
However, for $\beta >0$, the value of $Q_1$ is not evident.\par 
\indent To obtain bounds on $v\mapsto Q_2(v,0)$, we observe that 
$$Q_2(v,0)-Q_2(0,0)=\int_{\Omega_N} \int_0^{2\pi} (q(s+v)-q(s))q(s)
{{ds}\over{2\pi}} \nu_N^\beta (dq)\eqno(4.15)$$
\noindent so by the Cauchy--Schwarz inequality
$$\eqalignno{\bigl\vert Q_2(v,0)-Q_2(0,0)\bigr\vert&\cr
\leq\Bigl(\int_{\Omega_N}
\Bigl\vert \int_0^{2\pi} &\bigl(q(s+v)-q(s)\bigr)q(s)
{{ds}\over{2\pi}}\Bigr\vert^2 \nu_N^0 (dq)\Bigr)^{1/2}  \Bigl(\int_{\Omega_N}
\Bigl(
{{d\nu_N^\beta}\over{d\nu_N^0}}\Bigr)^2d\nu_N^0\Bigr)^{1/2}.&(4.16)\cr}$$
\noindent Now the inner integral reduces to
$$\int_0^{2\pi} \bigl(q(s+v)-q(s)\bigr)q(s)
{{ds}\over{2\pi}}=\sum_{n=1}^\infty {{2-2\cos nv}\over{n^2}}
\gamma_n\gamma_{-n},\eqno(4.17)$$
\noindent and so by independence
$$\eqalignno{\int_{\Omega_N}\Bigl\vert\int_0^{2\pi} \bigl(q(s+v)-q(s)\bigr)q(s)
{{ds}\over{2\pi}}\Bigr\vert^2 \nu_N^0 (dq)&\leq C_\beta \Bigl(
\sum_{n=1}^\infty {{2-2\cos nv}\over{n^2}}\Bigr)^2\cr
&\leq C_\beta v^2&(4.18)\cr}$$
\noindent for some constant $C_\beta$.\par

\indent We can also deal with integrals of the form (5.2) over Hilbert space, 
under suitable extra hypotheses regarding differentiability. 
By making a change of variables, we can express the left-hand side of (5.2) in terms of $F$ and its derivatives. 
The connection with transportation inequalities is suggested by Lemma 4.2 of [1]. We also require a basic definition.\par
\vskip.05in

\indent (ii) Let $(H, \langle \, ,\,\rangle )$ and $(E, \langle\, ,\,\rangle_E)$ be real and 
separable Hilbert spaces and $\iota :H\rightarrow E$ a linear and injective map which is 
Hilbert--Schmidt. A subset $C$ of $H$ is said to be cylindrical if there exist $m\in {\bf N}$, vectors 
$y_1, \dots , y_m\in H$ and a Borel subset $B$ of ${\bf R}^m$ such that
 $C=\{ x\in H: (\langle x,y_j\rangle )_{j=1}^m \in B\}$. Let ${\cal A}$ be the algebra of subsets of $H$
 that is generated by the cylindrical sets and  let $\gamma_0 :{\cal A}\rightarrow {\bf R}$ satisfy
$$\gamma_0 \bigl( \{ x\in H: (\langle x, e_j\rangle )_{j=1}^m \in B\}\bigr) =
{{1}\over{(2\pi )^{n/2}}}\int_{B}e^{-\sum_{j=1}^m x_j^2/2} dx_1\dots dx_m\qquad (B\in {\cal A})\eqno(5.5)$$
\noindent for some orthonormal basis $(e_j)$ of $H$. Now let $\mu$ be the measure on $E$ that is 
induced from $\gamma_0$ by the map $\iota :H\rightarrow E$.\par
\indent (iii) The norm of $E$ is said to be measurable for the Gaussian measure $\mu$ on the 
Borel subsets of $E$ if for all $\varepsilon >0$, 
there exists a finite rank orthogonal projection $P_\varepsilon$ on $H$ such 
that $\mu(\{ x\in H: \Vert \iota Px\Vert_E>\varepsilon \})<\varepsilon $ for all 
finite rank projections $P$ that are orthogonal to $P_\varepsilon$; see [23]. In the current
 paper, we reserve the term measure to mean a countably additive set function on a $\sigma$-algebra.\par  
\vskip.05in

\indent The next step is to choose $n$ for the given $\varepsilon >0$ by making a
uniform approximation of $F(\sqrt{\mu_1+1+N}, \dots )$.
\vskip.05in
\noindent {\bf Lemma 4.3} {\sl Let $F:\ell^2\rightarrow {\bf R}$ be
$1$-Lipschitz for the $\ell^2$ metric. Then there exist $C_1,C_2>0$ such that for all $\varepsilon >0$ 
there exists $n$, such that $F_n:\Omega_N\rightarrow {\bf R}$ given by 
$$F_n(q)=F\bigl(\sqrt{\mu_1(q)+N+1}, \dots , \sqrt{\mu_n(q)+N+1}, C_2(n+1),
C_2(n+2),\dots \bigr)\eqno(4.8)$$
\noindent is $K$-Lipschitz for the $\ell^2$ metric with
$K=c_2C(N)\max\{(3C_1N(N+1)/C_2\varepsilon )^8, 1\}$ and} 
$$\bigl\vert F_n(q)-F\bigl(\sqrt{\mu_1(q)+N+1},\sqrt{\mu_2(q)+N+1},\sqrt{\mu_3(q)+N+1}\dots  
\dots \bigr)\bigl\vert<\varepsilon /3\qquad (q\in \Omega_N) .\eqno(4.9)$$
 \vskip.05in
\noindent {\bf Proof.} As in (2.14), let $q\in \Omega_N$ be such that for some $C_1, C_2>0$ the periodic
spectrum satisfies $\lambda_{2j}-\lambda_{2j-1}\leq C_1N\vert \hat q(j)\vert$
and $C_2j^2\leq \lambda_{2j-1}$ for all $j=1, 2, \dots $. On
$$X=\prod_{j=1}^\infty [C_2^2j^2-N, c^2_4j^2+C_1N(N+1)],$$
\noindent define
$R:X\rightarrow {\bf R}^\infty :$ $(\mu_j)\mapsto (\sqrt{\mu_j +N+1})_{j=1}^\infty $ 
and observe that $R$ is Lipschitz for the $\ell^2$ metric. Moreover, we have
$$\sum_{j=n}^\infty \bigl( \sqrt{\mu_j+N+1}-C_2j\bigr)^2\leq \sum_{j=n}^\infty 
\Bigl( {{\mu_j+N+1-C_2^2j^2}\over {\sqrt{\mu_j+N+1}+C_2j}}\Bigr)^2\leq
{{(C_1N(N+1))^2}\over{C_2^2n}};\eqno(4.10)$$
\noindent so when we choose $n=\max\{ 1, (3C_1N(N+1)/C_2\varepsilon )^2\}$, and
replace $\sqrt{\mu_j(q)+N+1}$ by $C_2j$ in the $j^{th}$ place for all $j\geq n+1$, we only change $F(\sqrt{\mu_1+N+1}, \dots )$ by at
most $\varepsilon /3$. The function thus obtained is $F_n$ where $F_n(q)$ is Lipschitz with
constant $K=c_2C(N)n^4$ by Lemma 4.2. {$\square$}\par

 A function $F:{\bf
R}^\infty\rightarrow {\bf R}$ is said to cylindrical if $F$ depends on only $k$
coordinates for some $k$, and is given by some $f\in C^\infty_c 
({\bf R}^k; {\bf R}).$\par
\vskip.05in
\noindent {\bf Theorem 4.1} {\sl For all $\beta >0$ and $0<N<N_1$, there exist 
$\kappa (N, \beta ),
\kappa_0>0$ such that
all cylindrical functions $F:{\bf R}^\infty\rightarrow {\bf R}$ that are $1$
Lipschitz for the $\ell^2$ metric satisfy}
$$\eqalignno{\nu_N^\beta \Bigl\{ q\in\Omega_N:& F\bigl(\sqrt{\mu_1 (q)+N+1}, \sqrt{\mu_2 (q)+N+1}, \dots
\bigr)-\int
Fd\nu_N^\beta >\varepsilon \Bigr\}\cr
&\leq \exp\bigl(-\kappa (N, \beta )
\varepsilon^2 \min\{1, \kappa_0(\varepsilon /N(N+1))^8\}\bigr)\qquad
(\varepsilon >0).&(4.1)\cr}$$   
\vskip.05in
\noindent For large $\varepsilon >0$, this simplifies to a Gaussian upper bound; whereas for small
$\varepsilon /N^2$, there is an adjustment to the exponent. 

\noindent {\bf Proof of Theorem 4.1} By the uniform approximation (4.9),
there is a containment of events  
$$\Bigl\{ q\in \Omega_N:F-\int Fd\nu_N^\beta >\varepsilon
\Bigr\}\subseteq \Bigl\{ q\in \Omega_N : F_n-\int F_nd\nu_N^\beta>\varepsilon
/3\Bigr\};\eqno(4.12)$$
\noindent so we can use the estimate of Lemma 4.4 with $\eta =\varepsilon/3$ to obtain
the exponent in (4.11) of 
$${{\alpha (N, \beta )\varepsilon^2}\over{9K^2}}= {{\alpha (N, \beta )\varepsilon^2}\over{9c_2^2C(N)^2}}
\min \Bigl\{1, \Bigl({{3C_2\varepsilon}\over{C_1N(N+1)}}\Bigr)^8\Bigr\}\eqno(4.13)$$ 
to conclude the proof. {$\square$}\par

\noindent {\bf Lemma 5.1} {\sl Let $\mu$ be a Gaussian measure on $E$ as above. Suppose that 
$F:H\rightarrow {\bf R}$ is $L$-Lipschitz and has Frechet derivative $\nabla F$ which is continuously differentiable
and such that} ${\hbox{Hess}}(F)$ {\sl is  a bounded family of  Hilbert--Schmidt operators on $H$. Then there exists $t_0>0$ such that for $0<t<t_0$ there exists a probability measure $\mu_t$, which is absolutely continuous with respect to $\mu$, with Radon--Nikodym derivative} 
$${{d\mu_t}\over{d\mu}}=\exp\bigl({-t\langle \eta ,\nabla F(\eta )\rangle -t^2\Vert \nabla F(\eta
)\Vert^2/2+t\,{\hbox{trace (Hess)}}(F)}\bigr) \det_2\bigl (I+t{\hbox{Hess}}(F)\bigr)\eqno(5.6)$$
\noindent {\sl and such that all bounded and continuous $\Phi :E\rightarrow {\bf R}$ satisfy}
$$\int_E \exp\bigl({Q_tF(\xi )} \bigr)\Phi (\xi )\mu (d\xi )=
\int_E \exp\bigl({F(\eta )+(t/2)\Vert \nabla F (\eta )\Vert^2}\bigr)
\Phi (\eta +t\nabla F(\eta ))\mu_t(d\eta ).\eqno(5.7)$$
\vskip.05in
\noindent {\bf Remark 5.2} Whereas the summands in $\langle \eta ,\nabla F(\eta )\rangle -{\hbox{trace(Hess)}}(F)$ 
from (5.6) are not necessarily finite, their difference exists on the support of $\mu$, as discussed in
 detail in section 4 of [22].\par

\vskip.05in

\noindent {\bf Proof of Lemma 5.1.} For $0<t<1/L$,
 the map $\eta\mapsto \eta +t\nabla F(\eta )$ is open and invertible, 
by a simple application of Banach's fixed point theorem. 
Then for all $t>0$ sufficiently small, $\eta\mapsto F(\eta )
+(2t)^{-1}\Vert \xi-\eta \Vert^2$ is uniformly convex, bounded below and 
diverges to infinity as $\Vert\eta\Vert\rightarrow\infty$. 
 Hence the infimum convolution 
$Q_tF(\xi )=\inf \{ F(\eta )+(2t)^{-1} \Vert \xi -\eta\Vert^2\}$ is attained at 
$\eta$ such that $\xi=\eta +t\nabla F(\eta )$, and 
$Q_tF(\xi )=F(\eta )+(t/2)\Vert \nabla F(\eta )\Vert^2$. Then the pull-back of 
$\mu$ by $\eta\mapsto \eta+t\nabla F(\eta )$ is a measure $\mu_t$ that is 
absolutely continuous with respect to $\mu$, with Radon--Nikodym derivative 
(4.20), as computed by Ramer [22]. The formula (4.21) follows by the definition of the pull-back measure.\par
\vskip.05in
\noindent {\bf Theorem 5.3} {\sl (i) Let $F:H\rightarrow {\bf R}$ be a bounded
and continuous function. Then, for some $\rho =\rho (N,\beta )>0$, 
the exponential integral of the Hopf--Lax semigroup
satisfies 
$$\int_{\Omega_N}\exp\bigl(\rho {Q_tF(q)}\bigr)\nu_N^\beta (dq)\leq\exp
\int_{\Omega_N} \rho F(q)\nu_N(dq).$$  
\indent (ii) In particular, let $F$ be as in Lemma 5.1, and let $\omega_t$ be the pull back of $\nu_N^\beta$
by the map $\psi_t:q\mapsto q+t\nabla F(q)$ on $L^2[0, 2\pi ]$. Then}
$$\int_{\Omega_N}\exp\bigl({Q_tF(q)}\bigr)\nu_N^\beta (dq)=
\int_{\psi_t^{-1}(\Omega_N)}
\exp\bigl({F(q)+(t/2)\Vert\nabla F(q)\Vert^2}\bigr) \omega_t(dq).\eqno(5.8)$$

\vskip.05in
\noindent {\bf Proof.} (i) The logarithmic Sobolev inequality of [2] implies this
transportation inequality by Theorem 22.17 of [26].\par
\indent (ii) Let $P_n$ be the Riesz projection onto the subspace
${\hbox{span}}\{ e^{ikx}:k=-n,\dots ,n\}$ of $L^3({\bf T})$, and introduce a
sequence of continuous functions $\Phi_N:L^2\rightarrow{\bf R}$ by
$$\Phi_n(q)=\cases{ \exp\Bigl( {{\beta}\over{6}}\int_0^{2\pi} \bigl(
P_nq(x)\bigr)^3 {{dx}\over{2\pi}}\Bigr), & for $ \Vert q\Vert^2_{L^2}\leq N$;\cr
0,& for $\Vert q\Vert^2_{L^2}>N+1$;\cr}\eqno(5.9)$$
\noindent such that $\Phi_n(q)\rightarrow \Phi (q)$ as $n\rightarrow\infty$ where 
$$\Phi (q)=\cases{ \exp\Bigl( {{\beta}\over{6}}\int_0^{2\pi} \bigl(
q(x)\bigr)^3 {{dx}\over{2\pi}}\Bigr), & for $ \Vert q\Vert^2_{L^2}\leq N$;\cr
0,& for $\Vert q\Vert^2_{L^2}>N$.\cr}\eqno(5.10)$$
\noindent We can choose the $\Phi_n$ such that $\Phi_n(q)\leq \Psi (q)$ for some
integrable function of the form
$$\Psi (q)=\cases{ \exp\Bigl( {{\beta}\over{3}}\int_0^{2\pi}
\bigl(1+q(x)^4\bigr){{dx}\over{2\pi}}\Bigr), & for $ \Vert q\Vert^2_{L^2}\leq N$;\cr
0,& for $\Vert q\Vert^2_{L^2}>N+1$.\cr}\eqno(5.11)$$  
\noindent So by applying the dominated convergence theorem, we can deduce (5.8) 
from the definition of $\nu_N^\beta$ in (1.8) via Lemma 5.1.\par 
\noindent {\bf 6 Concentration of measure on linear statistics}\par
\vskip.05in
When $zg(z)\in RPW(\pi )$, the series $\sum_n g(n)$ converges absolutely. In Theorem 6.1, 
we prove that if $g(z)$ and $z^2g(z)$ are in $RPW(\pi )$ then the linear statistic
$$\Lambda_q(g)=\sum_{j=0}^\infty (g(\sqrt{\lambda_j})+g(-\sqrt{\lambda_j}))\eqno(1.10)$$
\noindent satisfies Gaussian concentration of measure for $q$ in the probability space
 $(\Omega_N, \nu_N^\beta )$. The proof uses a contour integration argument involving 
$\Delta$.

\indent In this section we impose further conditions on the test 
function $g$, and achieve a concentration inequality for infinite 
linear statistics. In [4], we gave a version for
finite sums.\par
Let $V=\{ g\in RPW(\pi ): z^2g(z)\in RPW(\pi )\}$ which is a dense linear subspace of
$RPW(\pi )$.\par 
\vskip.05in
\noindent {\bf Theorem 6.1} {\sl For $q\in \Omega_N$ and 
$g\in V$, let
$$\Lambda_q(g)=\sum_{j=0}^\infty (g(\sqrt{\lambda_j})+g(-\sqrt{\lambda_j})).\eqno(6.1)$$
\noindent Then there exists $L(N,\beta )>0$ such that}
$$\log \int_{\Omega_N} \exp\bigl( \Lambda_q(g)\bigr)\nu_N^\beta (dq)
\leq  L^2\int_{-\infty}^\infty (1+x^2)^2\vert g(x)\vert^2\,dx+\int_{\Omega_N} \Lambda_q(g)\nu_N^\beta (dq).\eqno(6.2)$$
\vskip.05in
\noindent {\bf Proof.} The measure $\nu_N^\beta$ satisfies a 
logarithmic Sobolev inequality, and hence a concentration inequality for all Lipschitz functions on
 $\Omega_N$ by [3]. Therefore, it suffices to prove that $q\mapsto \Lambda_q(g)$ is Lipschitz 
for $q\in \Omega_N$ for $g$ as in the statement.\par
\indent By Proposition 2.3(ii), we can assume that all the
periodic eigenvalues are simple. Then we note that $\lambda_j>0$ for all but finitely 
many indices $j$, and that in the exceptional cases where $\lambda_j<0$, the sum
 $g(\sqrt{\lambda_j })+g(-\sqrt{\lambda_j})$ unambiguously gives a real random variable for
 $q\in\Omega_N$ and $g\in RPW(\pi )$. Let $X,Y>\sqrt{N}$, and let $R$ be the rectangle that has
 vertices $\pm X\pm iY$  and with boundary $\partial R$ described once in the positive sense.
 For $\lambda$ in the interval of instability $(\lambda_{4j+3}, \lambda_{4j+4})$, we have 
$\Delta (\lambda )>2$. From the Cauchy Residue Theorem, we deduce
$${{1}\over{2\pi i}}\int_{\partial R} {{2z\Delta'(z^2)g(z)dz}\over {\Delta (z^2)-2}}
=\sum_j\bigl(g(\sqrt{\lambda_{3+4j}})+g(-\sqrt{\lambda_{3+4j}})+g(\sqrt{\lambda_{4+4j}}
)+g(-\sqrt{\lambda_{4+4j}})\bigr)\eqno(6.3)$$ 
\noindent where the sum in (6.3) is over those indices $j$ such that $\sqrt{\lambda_{4j+4}}<X$, 
and we adjust the 
formula when $\lambda_{4j+4}=X^2$. \par
\indent We aim to take the sum over all $j\geq 0$, and show that the right-hand side defines a 
Lipschitz function of $q$. To express the estimates, we use the following notation. 
For an operator $T_q$ on $L^\infty [0, 2\pi ]$ with kernel depending upon $q\in \Omega_N$, and $h\in L^2({\bf T}; {\bf R})$ has $\Vert h\Vert_{L^2}\leq 1$, we define the derivative of $T_q$ in the direction of $h$ to be the limit
$$\Bigl({{\partial T_q}\over{\partial q}}\Bigr)_hf(x)=\lim_{\varepsilon\rightarrow 0+} {{T_{q+\varepsilon h}-T_q}\over{\varepsilon}}f(x)\eqno(6.4)$$
\noindent when this exists. The following lemma gives the basic estimates.\par
\vskip.05in
\noindent {\bf Lemma 6.2} {\sl There exist positive constants $c_j$, depending upon $Y$ and $N$, but not upon $X$ and $q$, such that \par
\noindent (i) $\vert \Delta (z^2)\vert\leq c_1$;\par
\noindent (ii) $\vert z\Delta'(z^2)\vert \leq c_2$;
$$ \Bigl\Vert\Bigl( {{\partial \Delta (z^2)}\over{\partial q}}\Bigr)_h\Bigr\Vert 
\leq C_1;\leqno{{}\qquad }(iii)$$
$$ \Bigl\Vert z\Bigl( {{\partial \Delta' (z^2)}\over{\partial q}}\Bigr)_h\Bigr\Vert 
\leq C_2;\leqno{{}\qquad}(iv)$$
\noindent for all $z$ inside and on $R$ and all $h\in L^2({\bf T}; {\bf R})$ such that $\Vert h\Vert_{L^2}\leq 1$.}\par
\vskip.05in
\noindent {\bf Proof of Lemma 6.2} (i) For $z\in R$, we have the uniform bounds 
$\vert \cos 2\pi z\vert \leq \cosh 2\pi Y$ and $\vert \sin 2\pi z\vert\leq \cosh 2\pi Y$; 
the estimate on $\Delta (z^2)$ then follows from the series (3.7).

\indent (ii) The bounds of (i) also hold for the enlarged rectangle $R'$ with vertices $\pm (X+1)+\pm i(Y+1)$, possibly with a different constant $c_1'$. Then we deduce (ii) for the derivative of $\Delta (z^2)$ by applying the Cauchy integral formula.\par
\indent (iii) For $\vert z\vert \leq 4\sqrt{N}\cosh 2\pi Y$, we use the local bounds from the proof of 
Theorem 5.2. \par
 In particular,
$$\Bigl({{\partial V}\over{\partial q}}\Bigr)_hf(x)=\int_0^x \cos \sqrt{\lambda} (x-t)h(t)f(t)\, dt\eqno(6.5)$$
$$\Bigl({{\partial W}\over{\partial q}}\Bigr)_hf(x)=\int_0^x \sin \sqrt{\lambda} (x-t) h(t)f(t)\, dt\eqno(6.6)$$
\noindent satisfy
$$\Bigl\Vert \Bigl({{\partial V}\over{\partial q}}\Bigr)_h\Bigr\Vert ,\Bigl\Vert \Bigl({{\partial W}\over{\partial q}}\Bigr)_h\Bigr\Vert\leq \cosh 2\pi Y\eqno(6.7)$$
\noindent as operators on $L^\infty [0, 2\pi ]$ when $\sqrt{\lambda}\in R$.\par 
\indent For $\vert z\vert \geq 4\sqrt{N}\cosh 2\pi Y$, we differentiate (3.7) with respect to $q$ in the direction of $h$ and obtain the series
$$\eqalignno{\Bigl( {{\partial \Delta (z^2)}\over{\partial q}}\Bigr)_h&={{1}\over{z}}\int_0^{2\pi} (h(t)+h(2\pi -t))\sin z(2\pi -t)\cos zt\, dt\cr
&\quad+{{1}\over{z^2}}\Bigl( \bigl({{\partial W}\over{\partial q}}\bigr)_hV(\cos sz)(2\pi)+W\bigl({{\partial V}\over{\partial q}}\bigr)_h\cr
&\quad +\bigl({{\partial W}\over{\partial q}}\bigr)_hV(\sin sz)(2\pi)
 +W\bigl( {{\partial V}\over{\partial q}}\bigr)_h(\sin sz)(2\pi )\Bigr)+\dots&(6.8)\cr}$$
\noindent where the coefficient of $z^{-k}$ is a sum of $2^kk$ terms, each of which is a 
product of $k-1$ factors of $V$ and $W$, and a factor of $\partial V/\partial q$ or 
$\partial W/\partial q$. Using the bounds (3.4) on the summands in this series, we establish 
convergence and obtain the stated bound (iii).\par
\indent (iv) By uniform convergence, the series (6.8) gives a 
holomorphic function of $z$ on $R'$. Hence (iv) follows from (iii)
 by the Cauchy integral formula. {$\square$}\par
\vskip.05in
\indent To extend the summation to all indices $j$, we need to choose  special contour 
$\Gamma^-$ on which we can control the denominator $\Delta (z^2)-2$ in (6.3). The basic
 estimates for choosing $\Gamma^-$ are in the following lemma.\par
\vskip.05in
\noindent {\bf Lemma 6.3} {\sl Let $C_1$ be the curve $z=u+iu^{-1/3}$ for $u\geq X$. Then there exists $c_5$ depending upon $N$ and $Y$ such that\par
\noindent (i) $\vert \Delta (z^2)-2\cos2\pi z\vert \leq c_5/u$, and\par
\noindent  (ii) $\vert \Delta (z^2)-2\vert \geq 4\pi^2 u^{-2/3} -c_5u^{-1}$  for all $z$ on $C_1$.}
\vskip.05in
\noindent {\bf Proof.} (i) This follows from the series (3.7).\par
\indent (ii) From (i), we have 
$$\eqalignno{\vert \Delta (z^2)-2\vert &\geq 4\vert \sin^2\pi z\vert-c_5/u\cr
&=4\sin^2\pi u+4\sinh^2(\pi u^{-1/3})-c_5/u.&(6.9)\cr}$$
\rightline{$\square$}\par
\noindent {\bf Conclusion of the proof of Theorem 6.1} Now we construct the contour $\Gamma$ to 
replace $\partial R$. In the first quadrant, we let 
$$\Gamma_1=[iY,X+ iY]\oplus [X+iY, X+i/X^{1/3}]\oplus C_1,\eqno(6.10)$$
\noindent where $X-1/2$ is an integer and $X>\max\{ c_5^3, 4\sqrt{N}\cosh 2\pi Y \}$. Then we 
reflect $\Gamma_1$ in the imaginary axis, the origin and the real axis to produce curves
 $\Gamma_1, \Gamma_2, \Gamma_3$ and $\Gamma_4$ in each of the four quadrants; finally we 
combine and re-orient the $\Gamma_j$ to produce a contour $\Gamma^-$ that is described once in 
the positive sense and resembles $\partial R$ with whiskers that are asymptotic to the real axis. 
This $\Gamma^-$ has been so chosen that there exists $c_7>0$, depending only upon $N$ and $Y$ such that 
$$\vert \Delta (z^2)-2\vert \geq c_7\vert z\vert^{-2/3}\eqno(6.11)$$

\noindent for all $z$ in $\Gamma$. To see this, it clearly suffices to obtain the bounds for
 $z$ on $\Gamma_1$, since one can then deduce the other cases by symmetry. For $z$ on $C_1$, we 
use Lemma 6.3(ii). Then on the vertical line segment $[X+iY, X+iX^{-1/3}]$, we have a better 
lower estimate $\vert \Delta (z^2)-2\vert \geq 4-c_5/u$ by (6.9) and the choice of $X-1/2$. 
On the horizontal line segment 
$[iY, X+iY]$, we use the local estimates as in Theorem 4.1.\par
\indent With this $\Gamma^-$, we express the linear statistic as a contour integral
$$\eqalignno{ \Lambda^-_q(g)&= \sum_{j=0}^\infty \bigl( g(\sqrt{\lambda_{3+4j}})
 +g(-\sqrt{\lambda_{3+4j}})+g(\sqrt{\lambda_{4+4j}}) +g(-\sqrt{\lambda_{4+4j}})\bigr)\cr
&={{1}\over{2\pi i}}\int_{\Gamma^-} {{z\Delta' (z^2)g(z)dz}\over{\Delta (z^2)-2}}&(6.12)\cr}$$
\noindent which has derivative with respect to $q$ in the direction of $h$ given by

$$\Bigl({{\partial \Lambda_q(g)}\over{\partial q}}\Bigr)_h={{1}\over{2\pi i}}\int_{\Gamma^-}
 \Bigl( {{2z }\over{\Delta (z^2)-2}} \bigl( {{\partial \Delta'(z^2)}\over{\partial q}}\bigr)_h - 
{{4z\Delta'(z^2)}\over{(\Delta (z^2)-2)^2}}\bigl( {{\partial \Delta (z^2)}\over{\partial q}}\bigr)_h
\Bigr) g(z)\, dz.\eqno(6.13)$$
\noindent Exploiting our previous estimates from Lemma 6.2 and (6.7), 
we note that 
$$  \Bigl\vert  {{2z }\over{\Delta (z^2)-2}} \bigl( {{\partial \Delta'(z^2)}\over{\partial q}}
\bigr)_h\Bigr\vert =O (\vert z\vert^{2/3} )\eqno(6.14)$$
\noindent and

$$\Bigl\vert {{4z\Delta'(z^2)}\over{(\Delta (z^2)-2)^2}}\bigl( 
{{\partial \Delta (z^2)}\over{\partial q}}\bigr)_h\Bigr\vert =O (\vert z\vert^{4/3}),\eqno(6.15)$$
\noindent where the implied constants depend upon $N$ but are uniform for $q\in \Omega_N$. From the Cauchy--Schwarz inequality we have 
$$\int_{\Gamma^-} (1+\vert z\vert^{4/3})\vert g(z)\vert \vert dz\vert\leq \Bigl( \int_{\Gamma^-} (1+\vert z\vert^{4/3})^3\vert g(z) 
\vert^2\vert dz\vert\Bigr)^{1/2}\Bigl(\int_{\Gamma^-} {{\vert dz\vert}\over{1+\vert z \vert^{4/3}}}\Bigr)^{1/2}\eqno(6.16)$$
\noindent and the right-hand side converges since  
$$\int_{-\infty}^\infty \vert (1+(x^2+y^2)^2)\vert g(x+iy)\vert^2dx\leq e^{2\pi Y}\int_{-\infty}^\infty (1+x^4)\vert g(x)\vert^2 dx\qquad (-Y<y<Y)\eqno(6.17)$$
\noindent converges when $z^2g(z)$ and $g(z)$ are in $RPW(\pi )$. Hence $q\mapsto \Lambda_q(g)$ is Lipschitz for $q\in \Omega_N$.\par
\indent Now for a slightly different contour $\Gamma^+$, we can obtain similar estimates on
$$\eqalignno{ \Lambda_q^+(g)&=\sum_j\bigl( g(\sqrt{\lambda_{1+4j}})+g(-\sqrt{\lambda_{1+4j}}) +g(\sqrt{\lambda_{2+4j}})+g(-\sqrt{\lambda_{2+4j}})\bigr)\cr
&={{1}\over{2\pi i}}\int_{\Gamma^+} {{z\Delta'(z^2) g(z)dz}\over{\Delta (z^2)+2}},&(6.18)\cr}$$
\noindent and show that this too is a Lipschitz function of 
$q\in\Omega_N$. The result follows. {$\square$}\par
\vskip.05in

\noindent {\bf Lemma 7.1} {\sl The function $t\mapsto G_n(t)$ is Lipschitz with respect to the metric}
$$d(t,t')=\Bigl( \int_{-\pi}^\pi \sum_{j=1}^n 4\sin^2 {{s(t_j-t_j')}\over {2}} {{ds}\over{2\pi}}\Bigr)^{1/2}\qquad (t,t'\in {\bf R}^n).\eqno(7.2)$$
\vskip.05in
\noindent {\bf Proof.}  For $t=(t_j)_{j=1}^n$ and $s=(s_j)_{j=1}^n$, we also introduce the linear statistic $\ell_n (t;s)=\sum_{j=1}^n \cos s_jt_j$, and the local deviations
$$\sigma_n(t;s)=\Bigl( 4\sum_{j,k=1}^n \sin^2
(s_kt_j/2)\Bigr)^{1/2},\eqno(7.3)$$
$$\rho_n(t;s)=\Bigl( 4\sum_{j=1}^n \sin^2(s_jt_j/2)\Bigr)^{1/2}.\eqno(7.4)$$
\vskip.05in
\indent  The Gram determinant satisfies
$$ G_n(t)={{1}\over{n!}}\int_{[-\pi ,\pi ]^n} \bigl\vert \det_2[ e^{it_js_k}]_{j,k=1}^n \bigr\vert^2 e^{2\ell_n(t;s)-2n} {{ds_1}\over{2\pi}}\dots {{ds_n}\over{2\pi}}\eqno(7.5)$$
\noindent and in particular $G_n(1,2,\dots ,n)=1$. By And\'eief's identity, we express the determinant of 
inner products in $L^2[-\pi ,\pi]$ and as a multiple integral
$$\eqalignno{ G_n(t)&=\det\Bigl[ \int_{-\pi}^\pi e^{i(t_j-t_k)s}{{ds}\over{2\pi}}\Bigr]_{j,k=1}^n\cr
&={{1}\over{n!}}\int_{[-\pi ,\pi ]^n}\Bigl\vert \det \bigl[e^{it_js_k}\bigr]_{j,k=1}^n \Bigr\vert^2 {{ds_1}\over{2\pi}}\dots {{ds_n}\over{2\pi}}\cr
&={{1}\over{n!}}\int_{[-\pi ,\pi ]^n}\Bigl\vert\det_2 \bigl[e^{it_js_k}\bigr]_{j,k=1}^n \Bigr\vert^2e^{2\sum_{j=1}^n \cos s_jt_j -2n} {{ds_1}\over{2\pi}}\dots {{ds_n}\over{2\pi}}&(7.6)\cr}$$
\noindent where we have formed Carleman's determinant for $A=[e^{it_js_k}]-I_n$. \par
\indent The functions ${\hbox{sinc}}(t-j))_{j=1}^n$ give an 
orthonormal sequence in $PW(\pi )$, and hence their
 Gram determinant is one. Alternatively, one can reduce 
(7.4) to a Vandermonde determinant which can be evaluated using
 standard computations from the theory of the classical compact 
Lie groups. In [11], Katz and Sarnak consider the minors $G_n(t)$ in 
relation to the Weyl denominator formula for the classical groups. 
The version that we require relates to $U(n)$. \par 

We introduce the sampling maps
$$\eqalignno{R:RPW(\pi )\rightarrow \ell^2:& f\mapsto (\dots ,0, f(t_{-n_0+1}),f(t_{-n_0+2}), \dots , 
f(t_{n_0-1}), 0,0 \dots ),\cr 
T: RPW(\pi )\rightarrow \ell^2:& f\mapsto (\dots , f(t_{-n_0}), 
0,\dots , 0, f(t_{n_0}), f(t_{n_0+1})\dots )&(7.10)\cr}$$ 
\noindent so that $R$ has finite rank and samples at the $t_j$ for $-n_0+1\leq  j\leq n_0-1$ and $T$ samples
at the $t_j$ for $j\leq -n_0$ and $j\geq n_0$. Then by 
re-indexing the standard basis of $\ell^2$, we can write $(R+T)(R+T)^\dagger$ as the matrix 
$$\left[\matrix{TT^\dagger& RT^\dagger\cr TR^\dagger
&RR^\dagger\cr}\right]=\left[\matrix{I_\infty +A&B\cr B^\dagger
&D+I_{2n_0-1}\cr}\right]\eqno(7.11)$$
\noindent where $D$ is a $(2n_0-1)\times (2n_0-1)$ matrix and $B$ is $\infty\times (2n_0-1)$.
Since ${\hbox{sinc}}(0)=1$, the matrix in (3.25) has ones on the leading diagonal, and we express
this in terms of the identity
matrices $I_{2n_0-1}$ and $I_\infty$.\par
\vskip.05in

\noindent {\bf Lemma 7.2} {\sl There exists $n_0=n_0(N, \beta )$ such that $I_\infty +A$ 
is invertible.}\par 
\vskip.05in

\noindent {\bf Proof.} A real sequence $(s_n)$ is said to be interpolating for $PW(\pi )$ 
if for all $(c_n)\in \ell^2$, there exists $f\in PW(\pi )$ such that $f(s_j)=c_j$ for all $j\in {\bf Z}$. 
Seip observed that this condition is equivalent to the linear operator 
$S:PW(\pi )\rightarrow \ell^2$ $f\mapsto (f(s_j))_{j=-\infty}^\infty$ being bounded and having an 
adjoint $S^\dagger : \ell^2\rightarrow PW(\pi )$ which is a linear embedding; 
equivalently, $({\hbox{sinc}}(t-s_j))_{j=-\infty}^\infty$ is a Riesz basis for its closed 
linear span. Note that these properties are inherited by subsequences of the original sequence $(s_j)$. \par
\indent By Kadet's $1/4$ theorem [20; 22, p227], a real sequence $(s_j)_{j=-\infty}^\infty$ such that 
$\sup\{ \vert s_j-j\vert : j\in {\bf Z}\}<1/4$  gives a Riesz basis $({\hbox{sinc}}(t-s_j))_{j=-\infty}^\infty $ 
for $PW(\pi )$. 
 For almost all $q$ in the support of $\nu_N^\beta$, we can choose $n_0=n_0(q)$ so large as to 
ensure that $\vert \sqrt{\lambda_n-\lambda_0}-n\vert <1/5$ for all 
$n\geq n_0$; hence $({\hbox{sinc}}(t-t_j))_{\vert j\vert \geq n}$ is a Riesz basis for its closed linear span. 
It follows that $R^\dagger$  is an embedding from a closed linear subspace of $\ell^2$ into $RPW(\pi )$. 
Hence $I+A$ is invertible.\par

\noindent {\bf 7 Concentration for the Gram determinant}\par
\vskip.05in
\noindent We let $t=(t_j)_{j=-\infty}^\infty$ be the modified eigenvalues of Hill's equation as in 
(1.10), let $S:PW(\pi )\rightarrow\ell^2:$ $f\mapsto (f(t_j))_{j=-\infty}^\infty$ be the sampling 
operator and consider the Gram determinant
$$G_n(t)=\det \bigl[{\hbox{sinc}}(t_j-t_k)\bigr]_{j,k=1}^n.\eqno(7.1)$$
\noindent This is a leading minor of the matrix that represents $SS^\dagger$ with respect to the 
standard basis on $\ell^2$, and thus describes some properties of the frame $({\hbox{sinc}}(t-t_k))_{k=1}^n$. In the trivial case of $t_j=j$ for all $j\in {\bf Z}$, the operator $SS^\dagger$ reduces to the identity and $G_n=1$.\par
\vskip.05in
\noindent {\bf Definition}(Carleman determinant) Let $A$ by a Hilbert--Schmidt operator with norm $\Vert A\Vert_{HS}$ and
 eigenvalues 
$(\kappa_j)_{j=0}^\infty $ listed according to multiplicity. Then the Carleman determinant of $I+A$ is defined by the convergent infinite product
$$\det_2(I+A)=\prod_{j=0}^\infty (I+\kappa_j)e^{-\kappa_j}.\eqno(5.4)$$
\indent The main result of this section is the following concentration inequality for finite Gram 
determinants.\par
\vskip.05in
\noindent {\bf Theorem 7.1} {\sl Let $\mu_n$ be the probability measure 
that is induced on ${\bf R}^n$ from $(\Omega_N, \nu_N^\beta )$ by the map 
$q\mapsto t$ where $t=(t_1, \dots , t_n)$. Then there exists $K_n>0$ such that}
$$\int_{{\bf R}^n} \exp\bigl({sG_n(t)}\bigr)\mu_n (dt)\leq
 \exp\Bigl({{{K_ns^2}\over{2}}+s\int G_n(t)\mu_n(dt)}\Bigr) \qquad (s\in {\bf R})\eqno(7.2)$$
\vskip.05in
\noindent The proof will occupy the remainder of this section, and starts with a basic definition.\par
\vskip.05in
\indent For $t=(t_j)_{j=1}^n$ and $s=(s_j)_{j=1}^n$, we also introduce the linear statistic $\ell_n (t;s)=\sum_{j=1}^n \cos s_jt_j$, and the local deviations
$$\sigma_n(t;s)=\Bigl( 4\sum_{j,k=1}^n \sin^2
(s_kt_j/2)\Bigr)^{1/2},\eqno(7.3)$$
$$\rho_n(t;s)=\Bigl( 4\sum_{j=1}^n \sin^2(s_jt_j/2)\Bigr)^{1/2}.\eqno(7.4)$$
\vskip.05in
\indent In [11], Katz and Sarnak consider the minors $G_n(t)$ in relation to the Weyl denominator formula for the classical groups. The version that we require relates to $U(n)$.\par
\vskip.05in 
\noindent {\bf Lemma 7.2} {\sl (i) The Gram determinant satisfies}
$$ G_n(t)={{1}\over{n!}}\int_{[-\pi ,\pi ]^n} \bigl\vert \det_2[ e^{it_js_k}]_{j,k=1}^n \bigr\vert^2 e^{2\ell_n(t;s)-2n} {{ds_1}\over{2\pi}}\dots {{ds_n}\over{2\pi}}\eqno(7.5)$$
\noindent {\sl and in particular $G_n(1,2,\dots ,n)=1$.}\par
\indent {\sl (ii) There exists an absolute constant such that} 
$$\int_{[-\pi ,\pi ]^n} e^{z\ell_n(t;s)} {{ds_1}\over{2\pi}}\dots {{ds_n}\over{2\pi}}\leq \min \Bigl\{ e^{zn}, C J_0(iz)^n\prod_{k=1}^n {{k}\over{t_k}}\Bigr\} \qquad (z>0).\eqno(7.6)$$
\noindent {\sl where $J_0$ is Bessel's function of the first kind of order zero.}\par
\vskip.05in
\noindent {\bf Proof of Lemma 7.2} (i) By And\'eief's identity, we express the determinant of 
inner products in $L^2[-\pi ,\pi]$ and as a multiple integral
$$\eqalignno{ G_n(t)&=\det\Bigl[ \int_{-\pi}^\pi e^{i(t_j-t_k)s}{{ds}\over{2\pi}}\Bigr]_{j,k=1}^n\cr
&={{1}\over{n!}}\int_{[-\pi ,\pi ]^n}\Bigl\vert \det \bigl[e^{it_js_k}\bigr]_{j,k=1}^n \Bigr\vert^2 {{ds_1}\over{2\pi}}\dots {{ds_n}\over{2\pi}}\cr
&={{1}\over{n!}}\int_{[-\pi ,\pi ]^n}\Bigl\vert\det_2 \bigl[e^{it_js_k}\bigr]_{j,k=1}^n \Bigr\vert^2e^{2\sum_{j=1}^n \cos s_jt_j -2n} {{ds_1}\over{2\pi}}\dots {{ds_n}\over{2\pi}}&(7.7)\cr}$$
\noindent where we have formed Carleman's determinant for $A=[e^{it_js_k}]-I_n$. \par
\indent The functions ${\hbox{sinc}}(t-j))_{j=1}^n$ give an orthonormal sequence in $PW(\pi )$, and hence their Gram determinant is one. Alternatively, one can reduce (7.4) to a Vandermonde determinant which can be evaluated using standard computations from the theory of the classical compact  Lie groups.\par 
\indent (ii) To obtain the estimate independent of $t$, one uses the simple bound $\ell_n (t;s)\leq n$. For the more refined estimate, one uses a computation which is familiar from the theory of random flight. One factors the integral factors as a product, and notes that Bessel's function of order zero satisfies
$$J_0(iz)=\int_{-\pi}^\pi e^{z\cos s}{{ds}\over{2\pi}}.\eqno(7.8)$$
In particular, we can take $M=n$ and $L_n=e^{c(2n+1)^2}$ and obtain the bounds
$$\Bigl\vert \det_2\bigl[ e^{it_js_k}\bigr]_{j,k=1}^n\Bigr\vert \leq (1+L_n),\eqno(7.7)$$
\noindent and
$$\eqalignno{\Bigl\vert \det_2\bigl[
e^{it_js_k}\bigr]_{j,k=1}^n-\det_2\bigl[e^{it'_js_k}\bigr]_{j,k=1}^n\Bigr\vert&\leq L_n\Bigl\Vert 
\bigl[ e^{it_js_k}\bigr]_{j,k=1}^n-\bigl[ e^{it'_js_k}\bigr]_{j,k=1}^n\Bigr\Vert_{HS}\cr
&=L_n\sigma_n(t-t',s).&(7.8)\cr}$$
\noindent The linear statistic $\ell_n$ satisfies the trivial bound $\ell_n(t;s)\leq n$, and by the Cauchy--Schwarz inequality
$$\eqalignno{\bigl\vert \ell_n(t;s) -\ell_n(t';s)\bigr\vert &= \Bigl\vert \sum_{j=1}^n 2\sin \bigl( (t_j-t_j')s_j/2\bigr) \sin \bigl( (t_j+t_j')s_j/2\bigr)\Bigr\vert\cr
&\leq \sqrt {n} \rho_n(t-t';s).&(7.9)\cr}$$
\noindent The stated result  follows from this.\par
\vskip.05in
\indent Passing from the case of finite to infinite sampling sets,

\vskip.05in
\noindent {\bf Proof of Theorem 7.1}  We prove that $t\mapsto G_n(t)$ is Lipschitz with respect to the metric
$$d(t,t')=\Bigl( \int_{-\pi}^\pi \sum_{j=1}^n 4\sin^2 {{s(t_j-t_j')}\over {2}} {{ds}\over{2\pi}}\Bigr)^{1/2}\qquad (t,t'\in {\bf R}^n).\eqno(7.9)$$
\noindent  On the convex level set $\{ A: \Vert A\Vert_{HS}\leq M\}$, the function $A\mapsto \det_2(I+A)$ is Lipschitz with constant $L\leq e^{c(2M+1)^2}$, where $c>0$ is some universal constant. In particular, we can take $M=n$ and $L_n=e^{c(2n+1)^2}$ and obtain the bounds
$$\Bigl\vert \det_2\bigl[ e^{it_js_k}\bigr]_{j,k=1}^n\Bigr\vert \leq (1+L_n),\eqno(7.10)$$
\noindent and
$$\eqalignno{\Bigl\vert \det_2\bigl[
e^{it_js_k}\bigr]_{j,k=1}^n-\det_2\bigl[e^{it'_js_k}\bigr]_{j,k=1}^n\Bigr\vert&\leq L_n\Bigl\Vert 
\bigl[ e^{it_js_k}\bigr]_{j,k=1}^n-\bigl[ e^{it'_js_k}\bigr]_{j,k=1}^n\Bigr\Vert_{HS}\cr
&=L_n\sigma_n(t-t',s).&(7.11)\cr}$$
\noindent The linear statistic $\ell_n$ satisfies the trivial bound $\ell_n(t;s)\leq n$, and by the Cauchy--Schwarz inequality
$$\eqalignno{\bigl\vert \ell_n(t;s) -\ell_n(t';s)\bigr\vert &= \Bigl\vert \sum_{j=1}^n 2\sin \bigl( (t_j-t_j')s_j/2\bigr) \sin \bigl( (t_j+t_j')s_j/2\bigr)\Bigr\vert\cr
&\leq \sqrt {n} \rho_n(t-t';s).&(7.12)\cr}$$
\indent By Lemma 7.2, we have
$$\eqalignno{ G_n(t)-G_n(t')&={{1}\over{n!}}\int_{[-\pi ,\pi ]^n} \Bigl( \bigl\vert \det_2\bigl[ e^{it_js_k}\bigr]_{j,k=1}^n\bigr\vert -
\bigl\vert \det_2\bigl[ e^{it'_js_k}\bigr]_{j,k=1}^n\bigr\vert\Bigr)\cr
&\times \Bigl(  \bigl\vert \det_2\bigl[ e^{it_js_k}\bigr]_{j,k=1}^n\bigr\vert + \bigl\vert \det_2\bigl[ e^{it_js_k}\bigr]_{j,k=1}^n\bigr\vert\Bigr) e^{2\ell_n(t;s)-2n} {{ds_1}\over{2\pi}}\dots {{ds_n}\over{2\pi}}\cr
&\quad + {{2}\over{n!}}\int_{[-\pi ,\pi ]^n}  \bigl\vert \det_2\bigl[ e^{it'_js_k}\bigr]_{j,k=1}^n\bigr\vert^2\bigl( \ell_n(t;s)-\ell_n(t';s)\bigr)\cr
&\qquad \times \int_0^1 
e^{-2n+u\ell_n(t;s)+(1-u)\ell_n(t';s)} du {{ds_1}\over{2\pi}}\dots {{ds_n}\over{2\pi}}&(7.13)\cr}$$
\noindent where the inner integral in the final expression emerges from DuHamel's formula. By the Cauchy--Schwarz inequality and Lemma 7.2(ii), we have
$$\int_{[-\pi ,\pi ]} \sigma_n (t-t';s)e^{2\ell (t;s)-2n} {{ds_1}\over{2\pi}}\dots {{ds_n}\over{2\pi}}$$
$$\leq {{C\sqrt{n} J_0(4i)^{n/2}}\over{e^{2n}}}\Bigl( \prod_{j=1}^n {{k}\over{t_k}}\Bigr)^{1/2}\Bigl(\int_{-\pi}^\pi\sum_{j=1}^n 4\sin^2 {{s(t_j-t_j')}\over{2}} {{ds}\over{2\pi}}\Bigr)^{1/2};\eqno(7.14)$$
\noindent and by H\"older's inequality we have
$$\int_{[-\pi ,\pi ]} \rho_n(t-t';s) \int_0^1 e^{-2n+u\ell_n(t;s)+(1-u)\ell_n(t';s)} {{ds_1}\over{2\pi}}\dots {{ds_n}\over{2\pi}}$$
$$\leq {{CJ_0(2i)^{n/2}}\over{e^{2n}}}\Bigl( \int_0^1 \prod_{k=1}^n {{k}\over{ t_k^u(t_k)^{1-u}}} du\Bigr)^{1/2}\Bigl( \int_{-\pi}^\pi \sum_{j=1}^n 4\sin^2 {{s(t_j-t_j')}\over{2}} {{ds}\over{2\pi}}\Bigr)^{1/2}.\eqno(7.15)$$ 
\noindent We deduce that 
$$\vert G_n(t)-G_n(t')\vert \leq M_n d(t,t')\eqno(7.16)$$
\noindent where $M_n$ is a constant given by the integrals. \par
\indent In section 4, we proved that the measure $\mu_n$ satisfies a logarithmic Sobolev inequality, and hence a concentration inequality for all Lipschitz functions. In particular, we obtain the theorem.\par
\indent We decompose the Gram matrix as 
$$ [{\hbox{sinc}}(t_j-t_k)]_{j,k=-\infty}^\infty =\left[\matrix{I_\infty +A&B\cr B^\dagger
&D+I_{2n_0-1}\cr}\right]\eqno(7.4)$$
\noindent where the non-zero entries are according to the ranges of indices:\par
\indent $\bullet$ $I$ for indices $j=k$;\par
\indent $\bullet$ $D$ for indices $j\neq k$ such that $\max\{ \vert j\vert, \vert k\vert\}\leq {n_0}-1$;\par
\indent $\bullet$ $A$ for indices $j\neq k$ such that $\min\{\vert j\vert, \vert k\vert\}\geq {n_0}$ and\par
\indent $\bullet$ $B$ for indices $j\neq \pm k$ such that 
$\min\{\vert j\vert ,\vert k\vert \}\leq {n_0}<\max\{ \vert j\vert , \vert k\vert \};$ \par
\noindent so that each of the blocks $D$, $A$ and $B$ are Hilbert--Schmidt as 
operators on $\ell^2$. \par
\indent Indeed by the Cauchy--Schwarz inequality, the bound on the central block is
$$\eqalignno{ \Vert D\Vert_{HS}^2&=\sum_{j,k=-{n_0}+1}^{{n_0}-1}\vert {\hbox{sinc}}(t_j-t_k)-\delta_{jk}\vert^2\cr
&=\sum_{j,k=-{n_0}+1}^{{n_0}-1}\Bigl\vert\int_{-1}^1 \bigl( e^{i\pi s(t_j-t_k)}-e^{i\pi s(j-k)}\bigr) {{ds}\over{2}}\Bigr\vert^2\cr
&\leq \sum_{j,k=-{n_0}+1}^{{n_0}-1}\int_{-1}^1 2\sin^2 {{s\pi (t_j-j+k-t_k)}\over{2}} ds\cr
&\leq\sum_{j,k=-{n_0}+1}^{{n_0}-1} {{\pi^2}\over{3}} (t_j-j+k-t_k)^2\cr
&\leq {{4\pi^2}\over{3}}(2{n_0}-1)\sum_{j=-{n_0}+1}^{{n_0}-1} (t_j-j)^2.&(7.5)\cr}$$
\noindent By simple estimates based upon Proposition 3.1, we have $\alpha >1/2$ and constants depending upon $\beta$ and $N$ such that 
$$\eqalignno{\Vert A\Vert^2_{HS}&=\sum_{j,k: j\neq k; \vert j\vert, \vert k\vert \geq {n_0}}
 {{\sin^2\pi (t_j-t_k)}\over{\pi^2 (t_j-t_k)^2}}\cr
&\leq\sum_{j,k: j\neq k; \vert j\vert, \vert k\vert \geq {n_0}} {{4\sin^2\pi
(t_j-j)+4\sin^2\pi (t_k-k)}\over{\pi^2 (t_j-t_k)^2}}\cr
&\leq \sum_{j,k: j\neq k;\vert j\vert, \vert k\vert \geq {n_0}}
 {{C(N, \beta )}\over{\vert j-k\vert^2}}\Bigl( {{1}\over{\vert j\vert^{2\alpha}}}
+{{1}\over{\vert k\vert^{2\alpha}}}\Bigr)\cr
&\leq 2C(N, \beta )\sum_{k={n_0}}^\infty {{1}\over{k^{2\alpha}}}\Bigl(\sum_{j=k+1}^\infty {{1}\over
{(j-k)^2}}\Bigr)\cr
&\leq {{C_1(N, \beta )}\over{{n_0}^{2\alpha -1}}};&(7.6)\cr}$$
\noindent and by a similar argument, the off diagonal block satisfies
$$\Vert B\Vert_{HS}^2\leq C_4(N, \beta ) \Bigl( {{1}\over{{n_0}}}+
{{\log {n_0}}\over {{n_0}^{2\alpha }}}\Bigr).\eqno(7.7)$$
\indent Hence $[{\hbox{sinc}}(t_j-t_k)]_{j,k=-\infty}^\infty$ has a Carleman determinant, 
and that its value satisfies
$$\bigl\vert \det_2 (I+D+A+B)-\det_2 (I+D)\bigr\vert\leq L_{n_0} \bigl(
 \Vert A\Vert_{HS}+\Vert B\Vert_{HS}\bigr),$$
\noindent where in this case the Lipschitz constant may be taken to be
 $L_{n_0}\leq Ce^{c\Vert D\Vert^2_{HS}}.$ Hence we have
$$\eqalignno{\int_{\Omega_N} \bigl( \det_2 (I+D+A+&B)-\det_2 (I+D)\bigr)^2 \nu_N^\beta (dq)\cr
&\leq {{C(N, \beta )}\over {{n_0}^{2\alpha }}}\int_{\Omega_N} e^{c(2{n_0}-1)\sum_{j=-{n_0}+1}^{{n_0}-1} 
(t_j-j)^2} \nu_N^\beta (dq).&(7.9)\cr}$$
\noindent Also, $D$ is zero on its leading diagonal, so we can replace $\det_2$ by $\det$ and hence by Lemma 7.1, we have 
$$\eqalignno{\int_{\Omega_N} \bigl( \det_2 (I+D)-1\bigr)^2 \nu_N^\beta (dq)&=\int_{\Omega_N} 
\bigl( G_{n_0}(t_{-{n_0}+1}, \dots ,t_{{n_0}-1}) -1\bigr)^2 \nu_N^\beta (dq).\cr
&\leq M_{n_0}\int_{\Omega_N} \sum_{j=-{n_0}+1}^{{n_0}-1} (t_j-j)^2\nu_N^\beta (dq)&(7.10)\cr}$$
\noindent So the condition
$$\int_{\Omega_N} \sum_{j=-{n_0}+1}^{{n_0}-1} (t_j-j)^2 \nu_N^\beta (dq)\rightarrow 0\eqno(7.11)$$
\noindent  as $\beta, N\rightarrow 0$ implies that 
$$\int_{\Omega_N} \det_2 \bigl[{\hbox{sinc}}(t_j-t_k)\bigr]_{j,k=-\infty}^\infty \nu_N^\beta 
(dq)\rightarrow 1,\eqno(7.12)$$
\noindent so the Gram matrix is invertible on a set of positive measure.\par

\indent  Proposition 8.1 is largely a summary of classical
results, and is included for completeness.\par
\vskip.05in
\noindent {\bf Definition} Let ${\cal K}({\bf C})$ be the set of all nonempty
compact subsets of ${\bf C}$, and for nonempty $A\subseteq {\bf C}$, let ${\cal
K}(A)=\{ K\in {\cal K}({\bf C}): K\subseteq A\}$. Now for a nonempty
open subset $U$ of ${\bf C}$, we define ${\cal K}(U)$ to be an open
subset of ${\cal K}({\bf C})$. Also, for a nonempty open subset $V$ of
${\bf C}$, we define ${\cal M}(V)=\{K\in {\cal K}({\bf C}): V\cap K\neq
\emptyset\}$, and say that ${\cal M}(V)$ is an open subset of 
${\cal K}({\bf C})$. The Michael topology is the coarsest topology 
that has basic open sets 
$${\cal K}(U)\cap {\cal M}(V_1)\cap \dots \cap {\cal M}(V_r)$$
\noindent for nonempty open subsets $U, V_1, \dots ,
V_r$ of ${\bf C}$.\par
\vskip.05in
\noindent {\bf Proposition 8.1} {\sl There exists a sequence of
polynomials $T_n(x)$ and hyperelliptic algebraic curves 
${\cal E}_n=\{ (x,y): y^2=4-T_n(x)^2\}$ such that:\par
\indent (i) the projection $\Sigma_n$ to ${\bf R}$ of the real points 
on ${\cal E}_n$
is a compact set;\par
\indent (ii) $D(0,R)\cap \Sigma_n\rightarrow D(0,R)\cap \Sigma$
in the Michael topology as $n\rightarrow\infty$ for all $R>0$;\par
\indent (iii the equilibrium measure on $\Sigma_n$ determines 
${\cal E}_n$; \par
\indent (iv) for all $g$ rational, $\sum_{j=1}^n g(\mu_j)$ is given by a meromorphic function on the
Jacobian of ${\cal E}_n$.}\par 
\noindent {\bf Proof.} (i) Totik has shown that for all $n$ and $\varepsilon>0$, there exist $a_1, \dots ,a_n\in {\bf R}$ such
that $\vert \lambda_{2j-1}-a_j\vert<\varepsilon/n$ and  
 a real polynomial $T_n$ of degree $N$ such that
all the zeros of $T_n'$ satisfy $T_n(y)^2\geq 4$ and $\Sigma_n=\cup_{j=1}^n
[a_j,\lambda_{2j}]$ satisfies 
$$\Sigma_n=\{ x\in {\bf C}: -2\leq T_n(x)\leq 2\}.\eqno(8.2)$$
\noindent hence ${\cal E}_n=\{ (x,y): y^2=4-T_n(x)^2\}$ is a
hyperelliptic curve of genus less than or equal to $N-1$, and for each $x\in \Sigma_n$ there
exists a real point on ${\cal E}_n$. \par

\indent (ii) For all $m$, we can choose a
$2^{-m-1}$ net $x_1, \dots ,x_r$ of the compact set $\Sigma\cap D(0,R)$
and introduce the open neighbourhood $U_m=\{ z\in {\bf C}:d(\Sigma\cap
D(0,R), z)<1/m\}$ of $\Sigma\cap D(0,R)$ and the open discs
$V_k=D(x_k,2^{-m})$ which are contained in $U_m$ and together cover
$\Sigma\cap D(0,R)$; then consider 
${\cal U}_m={\cal K}\cap {\cal M}(V_1)\cap \dots \cap {\cal M}(V_r)$.
Clearly $\Sigma\cap D(0,R)\in {\cal U}_m$, and we observe that any $\lambda\in {\bf
C}\setminus\Sigma$ will have $\lambda\in {\bf C}\setminus U_m$ for all
sufficiently large $m$, hence $\{\lambda\}$ will not be an element of
${\cal U}_m$. Likewise, if $K\in {\cal K}({\bf C})$ does not intersect
with $\Sigma$, then $K$ does not intersect with $V_j$ for all $j$ and
all sufficiently large $m$. However, for all $m$, there exists
$n_0$ such that $\Sigma_n\cap D(0,R)\in {\cal U}_m$ for all $n\geq n_0$.\par

\indent (iii) Now suppose that $\Sigma_n$ is made of conducting material and
supports a unit of positive electrical charge. Then the charge on $\Sigma_n$ has the equilibrium distribution
$$\omega_n(dx)={{\vert T_n'(x)\vert dx}\over{\pi N\sqrt{4-T_n(x)^2}}},
\eqno(8.3)$$  
\noindent and
$$T_n(x)=2\cos\Bigl( N\int_x^\infty \omega_n(dt)\Bigr)\qquad (x\in
\Sigma_n).\eqno(8.4)$$
\noindent Then $T_n$ is a polynomial of degree $N$, with leading coefficient
$\tau_N$, where $(1/\tau_N)^{1/N}$ equals the analytic capacity of
$\Sigma_n$.\par
\indent  To introduce the real Jacobi variety, also introduce the coefficients
$\sum_{k=0}^N\alpha_{k,n}x^k=4-T_n(x)^2$ and the generating function
$$S_n=\sum_{j=1}^n 2\int_{\mu_j}^{\lambda_{2j}}\sqrt{4-T_n(x)^2}\, dx.
\eqno(8.6)$$
\noindent Then the corresponding phases are
$$\varphi_{k,n}={{\partial S_n}\over{\partial \alpha_{k,n}}}=\sum_{j=1}^n
\int_{\mu_j}^{\lambda_{2j}}{{x^kdx}\over{\sqrt{4-T_n(x)^2}}},\qquad (k=0, \dots ,
n-1)\eqno(8.7)$$
\noindent and the real Jacobian map is $J_n:\times_{j=1}^n [a_j,
\lambda_{2j}]\rightarrow {\bf R}^n:$ 
$(\mu_j)_{j=1}^n\mapsto (\varphi_{k,n})_{k=0}^{n-1}$. The change of variables $J_n$ has derivative with determinant
$$\det\Bigl[
-{{\partial\varphi_{k-1}}\over{\partial\mu_j}}\Bigr]_{j,k=1}^n
={{\prod_{1\leq k<\ell\leq n} (\mu_\ell-\mu_j)}\over{\prod_{1\leq k\leq
n} \sqrt{4-T_n (\mu_k )^2}}},\eqno(8.8)$$
\noindent which is of constant sign on $\Sigma_n$.\par
\indent (iv) The linear statistic is symmetric and rational in each
variable, hence defines a meromorphic function on the symmetric
Cartesian product ${\cal E}_n^{n}/\sim$, which is equivalent to the
Jacobian variety.\par

\indent Suppose that $\Delta (\lambda )^2-4$ has only simple real zeros. Then for all 
$n\in {\bf N}$, $R>0$ and $\varepsilon >0$, we can introduce a Taylor
polynomial $T_n (\lambda )$ of $\Delta (\lambda )$ about $\lambda =0$ with degree $N_n$ so
large that $\vert T_n(\lambda )-\Delta (\lambda )\vert \leq \varepsilon $ for all
$\lambda $ such that $\vert \lambda \vert \leq R$, and for all $j=1, \dots ,2n$ there
exists $\hat\lambda_j\in {\bf R}$ such that $T_n(\hat\lambda_j)=\Delta (\lambda_j)$ and
$\vert \hat\lambda_j-\lambda_j\vert <\varepsilon$. Then ${\cal E}_n=\{ (z, \lambda ):
z^2=4-T_n(\lambda )^2\}$ defines a hyperelliptic curve of genus less than or equal to 
$N_n-1$ with branch points
$\hat\lambda_j$ for $j=1,\dots , 2n$ and additionally $j=2n+1, \dots , 2N_n$. Let
$4-T_n(x)^2=\sum_{k=0}^{2N_n} \alpha_kx^k$ and regard $\alpha_0,\dots \alpha_{N_n-2}$ as new
real variables. The
generating function is
$$S_n=\sum_{j=1}^{N_n-1} \int^{\mu_j}_{\hat\lambda_{2j}} \sqrt{4-T_n(x)^2}\, dx.\eqno(8.1)$$   
\noindent and the corresponding angle variables are
$$\varphi_{k}={{\partial S_n}\over{\partial \alpha_{k}}}=\sum_{j=1}^{N_n}
\int_{\mu_j}^{\lambda_{2j}}{{x^kdx}\over{\sqrt{4-T_n(x)^2}}},\qquad (k=0, \dots ,
N_n-2)\eqno(8.2)$$
\noindent so the real Jacobian map is $J_{N_n}:\times_{j=1}^{N_n-1} [\hat\lambda_{2j},
\hat\lambda_{2j+1}]\rightarrow {\bf R}^{N_n-1}:$ 
$(\mu_j)_{j=1}^{N_n-1}\mapsto (\varphi_{k})_{k=0}^{N_n-2}$. By results of Jacobi
discussed in [], the real Jacobian of ${\cal E}_n$ gives the action and angle variables of an 
integrable dynamical system.\par

\indent We now extend the definitions to the context of ${\cal E}$. First we
motivate the sampling functions. Let $U_n$ be the Chebyshev polynomial of the second kind of degree $n$, 
and $h(z)$ a rational function. We observe that the 
rational integral on ${\cal E}_n$ 
$$\sum_{j=1}^n
\int_{a_j}^{\lambda_{2j}}h\Bigl(
{{T_n(s_j)}\over{2}}
\cos{{\pi t_j}\over{n+1}}-{{\sqrt{4-T_n
(s_j)^2}}\over{2}}\sin{{\pi t_j}\over{n+1}}\Bigr)
U_{n}\bigl({{T_n(s_j)}\over{2}}\bigr)
{{T_n'(s_j)\,ds_j}\over{\sqrt{4-T_n(s_j)^2}}}\eqno(8.2)$$
\noindent transforms via
$T_n(s_j)=2\cos(\pi\theta_j/(n+1))$ to an integral of $g(\theta )=h(\cos (\pi\theta
/(n+1))$, namely
$$\sum_{j=1}^n \int_0^{n+1}g(\theta_j+t_j) {{\sin\pi\theta_j }\over{(n+1)\sin
{{\pi\theta_j}\over{n+1}}}}\, d\theta_j,\eqno(8.3)$$
\noindent where 
$$\int_0^{n+1} g(\theta)^2{{d\theta}\over{\pi
(n+1)}}={{1}\over{\pi}}\int_{-1}^1 h(x)^2{{dx}\over{\sqrt{1-x^2}}}.\eqno(8.4)$$
\vskip.05in 
\noindent {\bf Proposition 8.1} {\sl Let 
$V=\{ g\in RPW(\pi ): z^2g(z)\in RPW(\pi )\}$. Then For all $g\in V$,} 
$$\sum_{j=1}^n \int_0^{n+1}g(\theta_j+t_j) {{\sin\pi\theta_j }\over{(n+1)\sin
{{\pi\theta_j}\over{n+1}}}}\, d\theta_j\rightarrow\sum_{j=1}^\infty
\int_0^\infty g(\theta +t_j)
{\hbox{sinc}}(\theta ) d\theta\qquad (n\rightarrow\infty ).\eqno(8.5)$$
\vskip.05in
\noindent {\bf Proof.} For $g\in V$, the function $xg(x)$ is integrable so the
result follows by the dominated
convergence theorem.\par

\noindent {\bf Proposition 9.1} {\sl Let $\Delta
(\mu_j)^2-4=(\lambda_{2j}-\mu_j)(\mu_j-\lambda_{2j-1})h_j(\mu_j)$. Then
$h_j$ is continuous and strictly positive and there exists $\mu_j'\in (\lambda_{2j-1}, \lambda_{2j})$ such that the
image of $\times_{j=1}^n [\lambda_{2j-1}, \lambda_{2j}]$ under the
transformation  $J:(\mu_j)_{j=1}^n\mapsto (\alpha_j)_{j=1}^n$ has volume}
$${{\pi^n\prod_{1\leq j<k\leq n}(\mu'_k-\mu'_j)}\over{\prod_{j=1}^n
\sqrt{h_j(\mu'_j)}}}.$$ 
\vskip.05in
\noindent {\bf Proof.} 
\noindent and the image of $\times_{j=1}^n [\lambda_{2j-1}, \mu_j]$
under this transformation has volume
$$\int_{\times_{j=1}^n [\lambda_{2j-1}, \mu_j]} \det\Bigl[
{{\partial\varphi_{k-1}}\over{\partial\mu_j}}\Bigr]_{j,k=1}^nds_1\dots
ds_n=\det\Bigl[ \int_{\lambda_{2j-1}}^{\mu_j}
{{s^{k-1}}\over{\sqrt{\Delta (s)^2-4}}}\Bigr]_{j,k=1}^n.$$
\noindent We introduce $\theta_j$ by
$\mu_j=2^{-1}(\lambda_{2j}+\lambda_{2j-1})+2^{-1}(\lambda_{2j}-\lambda_
{2j-1})\cos \theta_j$ and then apply the mean value theorem to the
resulting integral.\par
\indent The inverse spectral problem involves recovering $q$ from the 
spectral data consisting of the periodic spectrum 
$\{ \lambda_j: j\in {\bf N}\}$ and the family of tied spectra 
$\{ \mu_j(s):j\in {\bf N};
s\in [0, 2\pi ]\}$ of the translated potentials $q(x+s)$ 
for $s\in [0, 2\pi ]$. This gives $\delta (s)=\sum_j({\bf p}_j-{\bf
q}_j(s))\in {\hbox{Div}}$, where 
${\bf q}_j(s)=(\mu_j(s), \varepsilon\sqrt{\Delta (\mu_j(s))^2-4})$ is a
real point on ${\cal E}$, and the classical inversion formulas suggest that we
consider derivatives of $\omega_\infty (\delta (s))$. In the 
present setting, $q\in \Omega_N$ is typically not differentiable, 
and so the classical inversion formulas do not apply directly. \par
\indent By the construction of $\nu_N^\beta$ from random Fourier series, there exists for almost all $q$ in the support of $\nu_N^\beta$ a sequence $q_n\in\Omega_N\cap C^\infty$ such that 
$q_n\rightarrow q$ in $L^2$ as $n\rightarrow\infty$. Let $q_n(x)$ have periodic spectrum $\{ \lambda_{j,n}\}$ and let $q_n(s+x)$ have tied spectrum $\{ \mu_{j,n}(s)\}$ where $\lambda_{2j-1,n}\leq \mu_{j,n}(s)\leq \lambda_{2j,n}$. Then we introduce the discriminant $\Delta_{n}$ for Hill's equation with potential $q_n$, and the phase function, which is determined by the spectral data,
$$\varphi_n (s)=\sum_{j=1}^\infty 2\int_{\mu_{j,n}(s)}^{\lambda_{2j,n}} 
{{\Delta_{n}'(x)dx}\over{\sqrt{\Delta_{n}(x)^2-4}}}.\eqno(8.24)$$
\vskip.05in
\noindent {\bf Proposition 8.5} {\sl The limit of
the phase functions determines the potential $q$ by}
$$\varphi_n(x)-\varphi_n(0)\rightarrow \int_0^xq(s)ds -x\int_0^{2\pi} q(s) 
{{ds}\over{2\pi}}\qquad (n\rightarrow\infty ).\eqno(8.25)$$
\vskip.05in
\noindent {\bf Proof.} Trubowitz has shown that the smooth potentials $q_n$ give a differentiable phase function $\varphi_n$ such that 
$$q_n(x)-\int_0^{2\pi}q_n(s) {{ds}\over{2\pi}}={{d\varphi_n}\over{dx}},\eqno(8.26)$$
\noindent and we can integrate this formula, then take the limit of 
the left-hand side. {$\square$}\par